\input amstex
\documentstyle{amsppt}
\magnification=\magstep1
 \NoRunningHeads
 \loadbold
\topmatter
\title
 Rank-one nonsingular actions of countable groups and their odometer factors 
\endtitle

\author
Alexandre  I. Danylenko  and Mykyta I. Vieprik
\endauthor

\address
 B. Verkin Institute for Low Temperature Physics and Engineering
of the  National Academy of Sciences of Ukraine,
47 Nauky Ave.,
 Kharkiv, 61164, UKRAINE
  \newline 
 \newline
 Mathematical Institute of the Polish Academy of Sciences,
ul. {\'S}niadeckich 8,
 Warszawa, 00-656,
 POLAND
\endaddress
\email            alexandre.danilenko\@gmail.com
\endemail

\address
V. Karazin Kharkiv National University, 
4 Svobody sq.,  Kharkiv, 61077, Kharkiv,
UKRAINE\newline
\newline
 Mathematical Institute of the Polish Academy of Sciences,
ul. {\'S}niadeckich 8,
 Warszawa, 00-656,
 Poland
\endaddress
\email
nikita.veprik\@gmail.com
\endemail

\abstract
For an arbitrary countable discrete infinite group $G$, nonsingular rank-one actions are introduced.
It is shown that the class of nonsingular rank-one actions coincides  with the class of nonsingular $(C,F)$-actions.
Given a decreasing sequence $\Gamma_1\supsetneq\Gamma_2\supsetneq\cdots$ of cofinite subgroups in $G$ with $\bigcap_{n=1}^\infty\bigcap_{g\in G}g\Gamma_ng^{-1}=\{1_G\}$, the projective limit of the homogeneous $G$-spaces $G/\Gamma_n$ as $n\to\infty$ is a $G$-space.
Endowing this  $G$-space with an ergodic nonsingular nonatomic measure we obtain a dynamical system which is called a nonsingular odometer.  
Necessary and sufficient conditions are found for a rank-one nonsingular $G$-action to have a finite factor and a nonsingular  odometer factor in terms of the underlying $(C,F)$-parameters.
Similar conditions are also found for a rank-one nonsingular $G$-action to be isomorphic to an odometer.
Minimal Radon uniquely ergodic locally compact Cantor models are constructed for the nonsingular rank-one extensions of odometers.
Several concrete examples are constructed and several facts are proved that illustrate a sharp difference of the nonsingular noncommutative case from the classical finite measure preserving one: odometer actions which are not of rank one, factors of rank-one systems which are not of rank-one, however each probability preserving odometer is a factor of an infinite measure preserving rank-one system, etc.
\endabstract


\thanks
This work  was  supported in part by the 
``Long-term program of
support of the Ukrainian research teams at the 
Polish Academy of Sciences carried out in
collaboration with the U.S. National Academy of
Sciences with the financial support of external
partners'' and
by the Akhiezer Foundation. 
 \endthanks
\endtopmatter

\document

\head{0. Introduction}\endhead
This work is motivated by a recent paper \cite{Fo--We},
where  M.~Foreman, S.~Gao, A.~Hill, C.E.~Silva and B.~Weiss describe odometer factors of rank-one transformations in terms of the underlying cutting-and-stacking parameters. 
This description is considered as a step towards classification of the rank-one transformations up to isomorphism relation.
Our purpose here is to generalize the main results of \cite{Fo--We}   in the following 3 directions:
\roster
\item"---" to consider actions of arbitrary countable infinite groups including noname\-nable ones (\cite{Fo--We} deals only with $\Bbb Z$-actions),
\item"---" to consider arbitrary nonsingular group actions (\cite{Fo--We} deals only with  probability preserving actions) and 
\item"---" to consider rank-one actions along arbitrary sequence of shapes (\cite{Fo--We} deals only with  the classical rank-one, i.e. rank-one along a sequence of intervals in $\Bbb Z$. 
In particular,  our  results hold for the {\it funny} rank-one probability preserving $\Bbb Z$-actions which were not studied in \cite{Fo--We}).
\endroster

We now briefly outline the content of the paper, 
which consists of 6  sections.
  Section~1 is divided into 7 subsections.
In \S1.1 we define, for an arbitrary countable group $G$,  nonsingular $G$-actions of rank one.
According to this definition, a nonsingular $G$-action $T$ is of rank one if $T$ is free and $T$ admits a refining sequence of Rokhlin towers that approximate both the entire  $\sigma$-algebra of Borel subsets and the $G$-orbits and, in addition, the Radon-Nikodym derivative of $T$  
is constant
on each transposition of the levels  within  every tower 
(see Definition~1.1).
Definition~1.1 can be considered as an {\it abstract} definition of rank one.
In the case of probability preserving $\Bbb Z$-actions, there exist several  equivalent  {\it constructive} definitions of this concept \cite{Fe}. 
One of the most useful  of these is the  {\it cutting-and-stacking construction}, which explicitly 
associates a rank-one transformation to a sequence of integer-valued parameters.\footnote{Thus, the class of rank-one transformations is parametrized with a nice Polish space of integer parameters. However, different sequences of parameters can define isomorphic rank-one maps. 
A challenging open problem in this field is to find necessary and sufficient conditions for the parameters that determine isomorphic transformations.}
This transformation  is defined on the unit interval. 
It preserves the 
 Lebesgue measure.
 A natural generalization of this construction for general countable groups was suggested in \cite{dJ2} and 
\cite{Da1} in similar but non-equivalent versions. 
We call it the $(C,F)$-construction. 
The most general version of  the $(C,F)$-construction, including  the versions from \cite{dJ2} and \cite{Da1} as particular cases, appeared in \cite{Da3}.
However, \cite{Da3} deals only with {\it measure preserving} actions.
In \S1.2--1.3 here we 
define {\it nonsingular} $(C,F)$-actions.    
\S1.2 is preliminary: we define  $(C,F)$-equivalence relations and related quasi-invariant  $(C,F)$-measures.
The nonsingular $(C,F)$-actions  related to the $(C,F)$-equivalence relations and $(C,F)$-measures appear in~\S1.3.
They  include all the nonsingular rank-one transformations (and actions of Abelian groups) that have been studied earlier in the literature: see \cite{Aa}, \cite{AdFrSi}, \cite{HaSi}, \cite{Da1}, \cite{Da2}, \cite{DaSi} and references therein.
The main result of \S1 is the following (see Theorem~1.13).

\proclaim{Theorem A} Each nonsingular $(C,F)$-action of $G$ is of rank one and each rank-one nonsingular action of $G$ is isomorphic to a $(C,F)$-action.
 \endproclaim

It is worth noting that if a probability  preserving $G$-action is of rank one along a sequence $(F_n)_{n=1}^\infty$ of subsets in $G$ then $G$ is amenable and  $(F_n)_{n=1}^\infty$ is left F{\o}lner (see Corollary~1.11(ii)).

Important concepts of telescoping and reduction for the parameters of $(C,F)$-actions are introduced in \S1.4.
They are used in  \S1.5 to construct  continuous models of the nonsingular $(C,F)$-actions.
We remind that the famous Jewett-Krieger theorem provides strictly ergodic models for the ergodic probability preserving $\Bbb Z$-actions. 
In \cite{Yu}, an analogue of this theorem was proved for the {\it infinite} measure preserving ergodic transformations.
In the present paper, we prove
the existence of Radon uniquely ergodic minimal topological models for the {\it rank-one \it nonsingular} actions (see Theorem~1.19).

\proclaim{Theorem B} If $(X,\mu,(T_g)_{g\in R})$ is a nonsingular $G$-action of rank one then there are
a Radon uniquely ergodic minimal free continuous $G$-action $(R_g)_{g\in G}$ on a locally compact Cantor space $Y$, an $R$-quasi-invariant Radon measure $\nu$ on $Y$ and a measure preserving  isomorphism $\phi$ of $(X,\mu)$ onto $(Y,\nu)$ such that 
 $\phi T_g\phi^{-1}=R_g$ a.e. 
 and the Radon-Nikodym derivative $\rho_g:=\frac{d\nu\circ R_g}{d\nu}$
 is a continuous mapping from $Y$ to $\Bbb R_+^*$
 for each $g\in G$.
 Moreover, $\nu$ is the only (up to scaling) $R$-quasiinvariant Radon measure on $Y$ whose
 Radon-Nikodym cocycle equals $(\rho_g)_{g\in G}$.
\endproclaim

We note that the   continuity of the Radon-Nikodym derivatives was used essentially in \cite{DadJ} for the almost continuous orbit classification of nonsingular homeomorphisms of Krieger type $III$.

In \S1.6 and \S1.7 we discuss the case of  nonsingular $\Bbb Z$-actions of rank one along intervals in more detail.
It is shown in \S1.6 that the $(C,F)$-construction in this case is equivalent to the classical cutting-and-stacking with a single tower at every step of the construction.
However, in contrast with the  measure preserving case,  the towers are now divided into subtowers of different width. 
It is  shown in~\S1.7 that each  rank-one nonsingular transformation is isomorphic to  a transformation built over a classical nonsingular odometer of product type (called also Krieger's adding machine) and under a piecewise constant function.

Section~2 is devoted to description of finite factors of rank-one nonsingular actions. 
We remind that a {\it factor} of a dynamical system is an invariant sub-$\sigma$-algebra of measurable subsets.
Equivalently, a factor of a system is a dynamical system which appears as the image of the original system under a nonsingular  equivariant mapping.
Hence,  the finite factors of an ergodic $G$-action correspond to  the $G$-equivariant mappings onto homegeneous  $G$-spaces $G/\Gamma$, where $\Gamma$ is a cofinite subgroup in $G$. 
Each nonsingular $(C,F)$-action is parametrized by  an underlying sequence $\Cal T=(C_n,F_{n-1},\kappa_n,\nu_{n-1})_{n=1}^\infty$ of $(C,F)$-parameters, where $C_n$ and $F_{n-1}$ are finite subsets of $G$,
 $\kappa_n$ is a probability on $C_n$ and $\nu_{n-1}$ is measure on $F_{n-1}$ for each $n\in\Bbb N$.
 These parameters have to satisfy some conditions listed in \S1.3.
The following is the main result of \S2 (a stronger version of it is proved as Theorem~2.3; see 
also~Remark~2.4).

\proclaim{Theorem C} A  nonsingular $(C,F)$-action action $T$ of $G$ has a finite factor $G/\Gamma$ if and only  there is a telescoping $\Cal T=(C_n,F_{n-1},\kappa_n,\nu_{n-1})_{n=1}^\infty$ of the $(C,F)$-para\-meters of $T$ and a coset $g\Gamma\in G/\Gamma$ such that
$$
 \sum_{n=1}^\infty\kappa_n\big(\{c\in C_n\mid c\not\in g\Gamma g^{-1}\}\big)<\infty.
$$
An explicit formula for the factor mapping is obtained.
\endproclaim

We note that if  $G$  is Abelian, $\Gamma\subset G$ is a cofinite subgroup  and  
 the homogeneous $G$-space $G/\Gamma$ is a factor-space of an ergodic $G$-action $T$
  then the corresponding factor-algebra of $T$ is defined uniquely.
This is no longer true for non-Abelain $G$:
we provide an example of a rank-one $G$-action $T$ and two $T$-invariant sub-$\sigma$-algebras $\goth F_1\ne\goth F_2$ such that $T\restriction\goth F_1$ and $T\restriction\goth F_2$
are isomorphic $G$-actions on finite spaces (Example~2.5).

A criterion of total ergodicity for a nonsingular $(C,F)$-action in terms of the underlying $(C,F)$-parameters is obtained as a corollary from Theorem~C (see Corollary~2.6).

Starting from Section~3 we assume that $G$ is residually finite.
\S3 consists of 2 subsections. 
In \S3.1 we consider  {\it topological $G$-odometers} as the projective limits of  homogeneous $G$-spaces $G/\Gamma_n$ for a decreasing  sequence $\Gamma_1\supsetneq\Gamma_2\supsetneq\cdots$ of cofinite subgroups  in $G$ such that $\bigcap_{n=1}^\infty\bigcap_{g\in G}g\Gamma_ng^{-1}=\{1_G\}$.
By a {\it nonsingular $G$-odometer} we mean a topological $G$-odometer endowed with a $G$-quasiinvariant measure.
Topological properties of odometers  are not of our primary interest in the present work.
Measure theoretical    odometers (for general groups) were under study in \cite{LiSaUg} and \cite{DaLe}, but only in the finite measure preserving case.
In this paper we study nonsingular  $G$-odometers.
Some sufficient conditions for a nonsingular odometer to be of rank one are found in Propositions~3.2.
These conditions are satisfied for all known nonsingular odometers (see
Example~3.3).
It is worth noting  that there exist  odometers which are not of rank one.
Examples of non-rank-one probability preserving $G$-odometers  for non-amenable $G$ are given in  Example~3.4 and for amenable $G$ (including the Grigorchuk group) in Examples~3.5 and 3.6.
On the other hand, each probability preserving $G$-odometer is a factor of an infinite measure preserving rank-one $G$-action  (see Theorem~3.9 for a slightly stronger result):

\proclaim{Theorem D} For a topological $G$-odometer  $O$  defined on a compact space $Y$, there exist 
\roster
\item"---" a rank-one  measure preserving continuous $G$-action  $T$ on  a locally compact Cantor space $X$ equipped with a $\sigma$-finite measure  $\mu$ and
\item"---" a $G$-equivariant continuous mapping $\pi:X\to Y$
\endroster
 such that  $O$ is a factor of $T$ and the measure $\mu\circ\pi^{-1}$ is equivalent (i.e. has the same ideal of subsets of zero measure) to the Haar measure on $Y$.
\endproclaim

Thus, a factor of a rank-one nonsingular action is not necessarily of rank one.
This is in contrast with the classical case of rank-one finite measure preserving  $\Bbb Z$-actions \cite{Fe}.
Theorem~3.9 is about an interplay between odometer factors and an ``unordered'' sequence of finite factors for an ergodic $G$-action. 
This theorem is trivial in the case where $G$ is Abelian.

Nonsingular normal covers for nonsingular odometers are introduced  in \S3.2.
The existence of nonsingular normal covers is proved in Proposition~3.11.

In Section~4  we study  odometer factors of nonsingular $(C,F)$-actions.
The main result of the paper
 is the following (see Theorem~4.4).

\proclaim{Theorem E} Let $(X,\mu,T)$ be the  nonsingular $(C,F)$-action of $G$ associated with a sequence of $(C,F)$-parameters $\Cal T$.
Let  $O$ be the  topological $G$-odometer defined on the projective space $Y=\projlim_{n\to\infty}G/\Gamma_n$ corresponding to a nested sequence $\Gamma_1\supsetneq\Gamma_2\supsetneq\cdots$ of cofinite subgroups  in $G$ such that $\bigcap_{n=1}^\infty\bigcap_{g\in G}g\Gamma_ng^{-1}=\{1_G\}$.
A measurable $G$-equivariant mapping  $\pi:X\to Y$ exists if and only if there are  a telescoping 
$\Cal T'=(C_n,F_{n-1},\kappa_n,\nu_{n-1})_{n=1}^\infty$ of $\Cal T$ and an element $(g_n\Gamma_n)_{n=1}^\infty\in Y$ such that
$$
 \sum_{n=1}^\infty\kappa_n\big(\{c\in C_n\mid c\not\in g_n\Gamma_n g_n^{-1}\}\big)<\infty.
 $$ 
 An explicit formula for $\pi$ is obtained.
 Necessary and sufficient conditions for $\pi$ to be an isomorphism of $(X,\mu, T)$ onto
 $(Y,O,\mu\circ\pi^{-1})$ are given in terms of $\Cal T'$.
\endproclaim

It is worth noting that each rank-one nonsingular action $T$ is parametrized by the $(C,F)$-parameters $\Cal T$ (see Theorem~A) in a highly non-unique way. 
However, the properties of $\Cal T$ specified in the statement of Theorem~E (to determine an odometer factor $O$ or an isomorphism of $T$ with $O$) are independent on the choice of  $\Cal T$.
Hence Theorem~E can be considered as a contribution to the classification problem for the rank-one nonsingular systems.

As a corollary from Theorem~E,  criteria for a $(C,F)$-action to have an odometer factor
 or to be isomorphic to an odometer factor in terms of the underlying  $(C,F)$-parameters are obtained in Corollaries~4.6 and  4.7 respectively.
 Corollary~4.8 provides  minimal Radon uniquely ergodic models for the rank-one nonsingular extensions of nonsingular odometers.
 This corollary can be interpreted as a ``relative'' counterpart of Theorem~B:

 \proclaim{Theorem F} Let $(X,\mu,T)$ be a rank-one nonsingular $G$-action, $(Y,\nu,O)$ a nonsingular $G$-odometer and $\pi:X\to Y$ a $G$-equivariant mapping with $\mu\circ\pi^{-1}=\nu$.
 Then there exist a locally compact Cantor space $\widetilde X$, a minimal Radon uniquely ergodic free continuous $G$-action $\widetilde T$ on $\widetilde X$, a continuous $G$-equivariant mapping $\widetilde\pi:X\to Y$ and a Borel isomorphism $R:X\to\widetilde X$ such that
 $\mu\circ R^{-1}$ is a Radon measure on $\widetilde X$,
 $RT_gR^{-1}=\widetilde T_g$ for each $g\in G$,
  the Radon-Nikodym derivative 
 $\rho_g:=\frac{d(\mu\circ R^{-1})\circ \widetilde T_g}{d(\mu\circ R^{-1})}$
 is a continuous mapping from
 $\widetilde X$ to $\Bbb R^*_+$ for each $g\in G$
 and $\widetilde\pi\circ R=\pi$.
 Moreover, $\widetilde T$ is also Radon $(\rho_g)_{g\in G}$-uniquely ergodic.
 \endproclaim

It follows from  the Glimm-Effros theorem (see \cite{Ef}, \cite{DaSi}) that each topological odometer $(Y,O)$ (in fact, each topological $G$-action with a recurrent point) has uncountably many ergodic quasiinvariant measures.
However, the space of these measures is huge and   ``wild'' to describe it in good parameters.
Utilizing Theorems~E and F, we can isolate a good class  of ergodic finite quasi-invariant measures that admits a good parametrization. 
This is the  class of factor-measures on $Y$ for all rank-one nonsingular $G$-actions for which $Y$ is a factor. 
Every such measure can be parametrized by the $(C,F)$-parameters (see Corollary 4.9).

Section~5 is devoted completely to construction of 5 concrete rank-one actions with odometer factors and interesting properties.
In \S5.1 we continue study  the example of non-odometer rank-one probability preserving $\Bbb Z$-action $(X,\mu,T)$
 from  \cite{Fo--We}.
 It was shown there that the maximal odometer factor $\goth F$ of $T$  is non-trivial and isomorphic to the classical 2-adic odometer.
We prove  that $\goth F$ is the Kronecker factor of $T$ and that 
 $T$ is an uncountable-to-one extension of $\goth F$.
 It follows, in particular, that the spectrum of $T$ has a continuous component.
 In \S5.2 we
 consider nonsingular counterparts of the aforementioned system $(X,\mu,T)$.
 In particular, for each $\lambda\in[0,1]$, we construct a measure $\mu_\lambda$ on $X$
 such that 
 \roster
 \item"---"
 the triple  $(X,\mu_\lambda, T)$ is a rank-one nonsingular system of Krieger type $III_\lambda$, 
  \item"---" 
  $(X,\mu_\lambda, T)$ has a factor  $\goth F$ which is isomorphic to the probability preserving 2-adic odometer,
 \item"---"
 $\goth F$ is the maximal  (in the class of nonsingular odometers)  factor of $(X,\mu_\lambda, T)$,
  \item"---" the extension $X\to\goth F$ is uncountable-to-one (mod $\mu_\lambda$).  
  \endroster
In the  $III_0$-case, we extend this result   to systems whose associated flow is an arbitrary finitary AT in the sense of Connes-Woods \cite{CoWo}.
In \S5.3 we
provide an example of rank-one $\Bbb Z^2$-action $T=(T_{g})_{g\in \Bbb Z^2}$ such that the generators $T_{(0,1)}$ and $T_{(1,0)}$ have $\Bbb Z$-odometer factors but $T$ has no $\Bbb Z^2$-odometer factor.
Another  construction of such an action has appeared earlier in \cite{JoMc, \S6} but our example is much simpler.
In \S5.4 we construct a rank-one action $T$ of the Heisenberg group $H_3(\Bbb Z)$ which has an odometer factor $\goth F$ but which is not isomorphic to any odometer action.
We show there that $\goth F$ is the maximal odometer factor of $T$ and the extension $T\to\goth F$ is uncountable-to-one.
In \S5.5 we provide an example of non-normal $H_3(\Bbb Z)$-odometer which is canonically isomorphic to a normal odometer.

 The final \S6  is devoted to the  article \cite{JoMc} which  appeared in the course of  our work on the present paper.
The purpose of \cite{JoMc} is the same as ours: to generalize \cite{FoWe}.
However only finite measure preserving $\Bbb Z^d$-actions are studied in \cite{JoMc}.
 Therefore  in  \S6 we discuss  the results of  \cite{JoMc} and compare them with results of the present paper.

 {\sl Acknowledgement.} We thank R.~Grigorchuk and A.~Dudko for drawing our attention to Examples~3.5 and 3.6.

\head{1. Rank-one nonsingular actions  of countable groups and  $(C,F)$-construction}\endhead

\subhead 1.1. Nonsingular actions of rank one\endsubhead
Let $G$ be a discrete infinite countable  group. 
Let $T=(T_g)_{g\in G}$ be a free nonsingular action of $G$ on a standard $\sigma$-finite non-atomic measure space $(X,\goth B,\mu)$.
By a {\it  Rokhlin tower for $T$} we mean a pair $(B,F)$, where $B\in\goth B$ with $0<\mu(B)<\infty$ and $F$ is a finite subset of $G$ with $1_G\in F$ such that
\roster
\item"---" the subsets $T_fB$, $f\in F$, are mutually disjoint,
\item"---" the Radon-Nikodym derivative $\frac{d\mu\circ T_f}{d\mu}$ is constant on $B$ for each $f\in F$.
\endroster

Given a Rokhlin tower $(B,F)$, we let $X_{B,F}:=\bigsqcup_{f\in F}T_fB\in\goth B$.
Of course,  $\mu(X_{B,F})<\infty$.
By $\xi_{B,F}$ we mean the finite partition of  $X_{B,F}$ into the subsets $T_fB$, $f\in F$.
If $x\in T_fB$ then we set $O_{B,F}(x):=\{T_gx\mid g\in Ff^{-1}\}$.

\definition{Definition 1.1} Let $\{1_G\}=F_0\subset F_1\subset F_2\subset\cdots$ be an increasing sequence
of finite subsets in $G$.
We say that $T$ is of {\it rank-one along $(F_n)_{n=0}^\infty$} if there is a decreasing sequence 
$B_0\supset B_1\supset\cdots$ of subsets of positive measure in $X$ such that
$(B_n,F_n)$ is a
 Rokhlin tower for $T$ for each $n\in\Bbb N$ and
\roster
\item"(i)" $\xi_{B_n,F_n}\prec\xi_{B_{n+1},F_{n+1}}$ for each $n\ge 0$ and $\bigvee_{n=0}^\infty \xi_{B_n,F_n}$ is the partition of $X$ into singletons (mod 0),
\item"(ii)" $\{T_gx\mid g\in G\}=\bigcup_{n=1}^\infty O_{B_n,F_n}(x)$ for a.e. $x\in X$.
\endroster
\enddefinition

It follows from (i) that $X_{B_0,F_0}\subset X_{B_1,F_1}\subset X_{B_2,F_2}\subset\cdots$ and $\bigcup_{n=0}^\infty X_{B_n,F_n}=X$.
The piecewise constant property of the Radon-Nikodym derivative on the Rokhlin towers yields that:
\roster
\item"(iii)"
if $T_cB_{n+1}\subset B_{n}$ for some $c\in F_{n+1}$ then
$$
\frac{\mu(T_{fc}B_{n+1})}{\mu(T_fB_n)}=\frac{\mu(T_cB_{n+1})}{\mu(B_n)}\quad\text{for each $f\in F_n$ and $n\ge 0$}.
$$
\endroster

\proclaim{Proposition 1.2}
Let  $T$ satisfy (i)  from Definition~1.1.
Then  $T$ is ergodic. 
In particular, every rank-one nonsingular action is ergodic.
\endproclaim

\demo{Proof} Let two subsets $A_1,A_2\in\goth B$ be of positive measure.
It follows from (i) that there are $n>0$ and $f_1,f_2\in F_n$ such that
$$
\mu(A_1\cap T_{f_1}B_n)>0.9\mu(T_{f_1}B_n)\quad\text{and}\quad 
\mu(A_2\cap T_{f_2}B_n)>0.9\mu(T_{f_2}B_n).
$$
As $(B_n,F_n)$ is a Rokhlin tower,  $T_{f_2f_1^{-1}}T_{f_1}B_n=T_{f_2}B_n$ and
 $$
 \frac{d\mu\circ T_{f_2f_1^{-1}}    } {d\mu}(x)=
 \frac{\mu(T_{f_2}B_n)}{\mu(T_{f_1}B_n)}\quad\text{at a.e. $x\in A_1$.}
 $$
It follows that
$$
\align
\mu\big(T_{f_2f_1^{-1}}A_1\cap T_{f_2}B_n\big)&=\mu\big(T_{f_2f_1^{-1}}(A_1\cap T_{f_1}B_n)\big)\\
&=
\mu(A_1\cap T_{f_1}B_n) \frac{\mu(T_{f_2}B_n)}{\mu(T_{f_1}B_n)}\\
&>0.9\mu(T_{f_2B_n}).
\endalign
$$
Therefore $\mu(T_{f_2f_1^{-1}}A_1\cap T_{f_2}B_n\cap A_2)>0.8\mu(T_{f_2}B_n)$.
Hence, $\mu(T_{f_2f_1^{-1}}A_1\cap A_2)>0$, as desired. \qed
\enddemo

\subhead{1.2.  $(C,F)$-equivalence relations and nonsingular $(C,F)$-measures}
\endsubhead
Fix two sequences $(F_n)_{n\geq 0}$ and $(C_n)_{n\geq 1}$ of finite subsets in $G$ such that
$F_{0} = \{1_G\}$ and for 
 each $n>0$,
$$
\aligned
	  &  1_G\in F_n\cap C_n, \  \#C_{n} > 1, \\ 
      &F_{n} C_{n+1}\subset F_{n+1},\\ 
      &F_{n} c\cap F_{n} c' = \emptyset\text{ if $c, c'\in C_{n+1}$ and $c \neq c'$. }
\endaligned
\tag1-1
$$
We let $X_n := F_{n} \times C_{n+1} \times C_{n+2} \times\ldots$ and endow this set with the infinite product topology. 
Then $X_n$ is a compact Cantor space. The mapping 
$$
	X_n \ni (f_n,c_{n+1},c_{n+2}\ldots) \mapsto (f_n c_{n+1},c_{n+2},\ldots) \in X_{n+1}
$$
is a continuous embedding of $X_n$ into $X_{n+1}$. 
Therefore the topological  inductive limit $X$ of the sequence $(X_n)_{n\geq 0}$  is  well defined.
Moreover, $X$ is a locally compact Cantor space. 
Given a subset $A\subset F_n$, we let
$$
	[A]_n := \{x=(f_n,c_{n+1},\ldots)\in X_n, f_n\in A\}
$$
and call this set an $n$-{\it cylinder} in $X$.
It is open and compact in $X$. 
For brevity, we will write $[f]_n$ for $[\{f\}]_n$ for an  element $f\in F_n$.

We remind that two points $x=(f_n,c_{n+1},\ldots)$ and $x'=(f_n',c_{n+1}',\ldots)$ of $X_n$ are {\it tail equivalent}
if there is $N>n$ such that $c_l=c_{l}'$ for  each $l>N$.
We thus obtain the tail equivalence relation on $X_n$.

\definition{Definition 1.3 \cite{Da3}}
The {\it $(C,F)$-equivalence relation (or the tail equivalence relation)}  $\Cal R$ on $X$ is defined as follows: for each $n\ge 0$, the restriction of $\Cal R$ to $X_n$ is 
the tail equivalence relation on $X_n$. 
\enddefinition

The following properties of $\Cal R$ are easy to check:
\roster
\item"---"
Each $\Cal R$-class is countable.
\item"---"
 $\Cal R$ is {\it minimal}, i.e. the $\Cal R$-class of every point is dense in $X$. 
\item"---"
 $\Cal R$ is {\it hyperfinite}, i.e. there is a sequence $(\Cal S_n)_{n=1}^\infty$ of subrelations of $\Cal R$ such that  $\Cal S_1\subset \Cal S_2\subset\cdots$, $\bigcup_{n=1}^\infty\Cal S_n=\Cal R$ and $\#\Cal S_n(x)<\infty$ for each $x\in X$ and $n\ge 0$.
Indeed, we can define $\Cal S_n$ by the following:  $(x,y)\in\Cal S_n$ if either $x,y\not\in X_n$ and $x=y$ or
 $x=(f_n,c_{n+1},\dots)\in X_n,y=(f_n',c_{n+1}',\dots)\in X_n$ and $c_m=c_m'$ for all $m>n$.
 \endroster
 
 We remind that the {\it full group $[\Cal R]$ of $\Cal R$} is the group of all Borel bijections $\gamma:X\to X$ such that  $(x,\gamma x)\in\Cal R$
 for each $x\in X$.
A
Borel measure $\mu$ on $X$ is called $\Cal R$-quasi-invariant if $\mu\circ\gamma\sim\mu$ for each  $\gamma\in[\Cal R]$.
  Then there is a Borel mapping $\rho_\mu:\Cal R\to\Bbb R^*_+$ such that
 $$
 \rho_\mu(x,y)\rho_\mu(y,z)=\rho_\mu(x,z)\quad\text{ for  all $(x,y), (y,z)\in\Cal R$ }
 $$
 and
 $\rho_\mu(\gamma x,x)=\frac{d\mu\circ\gamma}{d\mu}(x)$ at a.e. $x\in X$ for each
 $\gamma\in[\Cal R]$.
 The mapping $\rho_\mu$ is called {\it the Radon-Nikodym cocycle of $(\Cal R,\mu)$}.

 \comment

\definition{Definition 1.4}  A  non-atomic $\sigma$-finite Radon measure $\mu $ on $X$ is called
a {\it  $(C,F)$-measure} if  $\mu$ is $\Cal R$-quasi-invariant and there exists a sequence $(a_n)_{n=0}^\infty$ of functions $a_n:F_n\to\Bbb R^*_+$ such that 
$$
\rho_\mu(x,y)=\frac{a_n(f_n)}{a_n(f_n')}
$$
whenever $x=(f_n,c_{n+1},\dots)\in X_n, y=(f_n',c_{n+1}',\dots)\in X_n$ and $(x,y)\in\Cal S_n$ for some $n\ge 0$.
\enddefinition

If $\mu$ is  a nonsingular $(C,F)$-measure then
$$
\frac{a_n(f)}{a_n(f')}=\frac{a_{n+1}(fc)}{a_{n+1}(f'c)}\quad\text{for all $f,f'\in F_n,\, c\in C_{n+1}$, $n\ge 0$.}\tag1-2
$$
In a sense, \thetag{1-2} can be resolved.
For that, we define a mapping $b_{n+1}:C_{n+1}\to\Bbb R_+^*$ by setting
 $b_{n+1}(c):=a_{n+1}(c)/a_{n}(1_G)$.
 Then  \thetag{1-2} yields that 
 $$
 a_{n+1}(fc)=a_n(f)b_{n+1}(c)\quad\text{for each $f\in F_n$, $c\in C_{n+1}$, $n\ge 0$.}
 $$
 It follows that if 
$x,y\in X_n$ for some $n\ge 0$  and $x$ and $y$ are tail equivalent then 
$$
\rho_\mu(x,y)=\frac{a_n(f_n)}{a_n(f_n')}\prod_{m>n}\frac{b_m(c_m)}{b_m(c_m')}.
$$

We now show how to construct explicitly  $(C,F)$-measures on $X$.

\endcomment

Suppose  that for each $n\in\Bbb N$, a non-degenerated probability measure $\kappa_n$ on $C_n$
is given.
We now let $\mu_0:=\nu_0\otimes\kappa_1\otimes\kappa_2\otimes\cdots$, where $\nu_0$ is the Dirac measure supported at $1_G$.
Then $\mu_0$ is an $(\Cal R\restriction X_0)$-quasi-invariant probability on $X_0$. 
Of course, $\mu_0$ is non-atomic if and only if 
$$
\prod_{n>0}\max_{c\in C_n}\kappa_n(c)=0.\tag1-2
$$
By the Kolmogorov 0-1 law, $(\Cal R\restriction X_0)$ is ergodic on the probability space $(X_0,\mu_0)$.
There are many ways to extend $\mu_0$ to an $\Cal R$-quasiinvariant measure on $X$.
However all such measures will be mutually equivalent.
Select for each $n\in\Bbb N$, a non-degenerated finite measure $\nu_n$ on $F_n$ such that 
$$
\nu_{n+1}(fc)=\nu_n(f)\kappa_{n+1}(c)\quad\text{ for each $f\in F_n$ and $c\in C_{n+1}$.}\tag1-3
$$
It is often convenient to consider $\nu_n$ and $\kappa_n$ as finite measures on $G$ supported on $F_n$ and $C_{n}$ respectively.
Then \thetag{1-3}
 means that $\nu_{n+1}\restriction F_nC_{n+1}=\nu_n*\kappa_{n+1}$, where the symbol $*$ means the convolution.
We now define a Borel measure $\mu$ on $X$ by setting 
$$
\mu([f]_n):=\nu_n(f), \quad\text{for each  $g\in F_n$ and every $n\in\Bbb N$.}
$$
It is  straightforward to verify that  $\mu$ is a well defined $\sigma$-finite Radon measure.
Moreover, $\mu$ is $\Cal R$-quasi-invariant
and
$$
\rho_\mu(x,y)=\frac{\nu_n(f_n)}{\nu_n(f'_n)}\prod_{m>n}\frac{\kappa_m(c_m)}{\kappa_m(c_m')},
$$
whenever $x=(f_n,c_{n+1},\dots)$ and $y=(f_n',c_{n+1}',\dots)$ are $\Cal R$-equivalent points that belong to $X_n=F_n\times C_{n+1}\times C_{n+2}\times\cdots$ for some $n>0$.
\comment

However, if we choose another sequence  $(\nu_n')_{n=1}^\infty$ satisfying ( ) and define the corresponding measure $\mu'$ on $X$ via ( ) with $\nu_n'$ in place of $\nu_n$ then $\mu'\sim\mu$.
Indeed, if $x\in X_n$ then $\frac{d\mu'}{d\mu}(x)=\frac{\nu_n'(f)}{\nu_n(f)}$, where $f$ is the element of $F_n$ such that $x\in[f]_n$.
\endcomment

The following definition extends \cite{Da1, Definition~4.2},  where the case of Abelian $G$ was considered.

\definition{Definition 1.4} If \thetag{1-2} and \thetag{1-3} hold then we call $\mu$  {\it the  $(C,F)$-measure on $X$ determined by $(\kappa_n)_{n=1}^\infty$ and $(\nu_n)_{n=0}^\infty$.}
\enddefinition

Consider another sequence $(\nu_n')_{n=0}^\infty$ of non-degenerated measures on $(F_n)_{n=0}^\infty$ respectively (in $n$) such that $\nu_0'$ is the Dirac measure supported at $1_G$ and $\nu_{n+1}'(fc)=\nu_n'(f)\kappa_{n+1}(c)$  for each $f\in F_n$ and $c\in C_{n+1}$ for each $n>0$.
Then the  $(C,F)$-measure $\mu'$ determined by $(\kappa_n)_{n=1}^\infty$ and $(\nu_n')_{n=0}^\infty$ is equivalent to $\mu$ and 
$$
\frac{d\mu'}{d\mu}(x)=\frac{\nu_n'(f_n)}{\nu_n(f_n)}\quad\text{if $x=(f_n,\dots)\in X_n$.}
$$
Another useful observation is that given $(\kappa_n)_{n=1}^\infty$, we can always find  $(\nu_n)_{n=0}^\infty$ satisfying \thetag{1-3}.
Thus, the equivalence class of a nonsingular $(C,F)$-measure  is  completely determined by $(\kappa_n)_{n=1}^\infty$ alone.
In particular, 
we may always replace a $\sigma$-finite nonsingular $(C,F)$-measure with an equivalent finite nonsingular $(C,F)$-measure.

\remark{Remark 1.5} We note that $\Cal R$ is {\it Radon uniquely ergodic}, i.e. there is a unique  $\Cal R$-invariant Radon measure $\xi$ on $X$ such that $\xi(X_0)=1$.
We call it  the {\it  Haar measure  for $\Cal R$}.
It is $\sigma$-finite.
Let $k_n$ be the equidistribution on $C_n$ and let $\nu_n(f)=\prod_{k=1}^n\kappa_k(1_G)$ for each  $f\in F_n$ and $n\ge 0$.
Then \thetag{1-2} and \thetag{1-3} holds for $(\kappa_n)_{n=1}^\infty$ and $(\nu_n)_{n=0}^\infty$.
Of course, the Haar measure for $\Cal R$ is
a $(C,F)$-measure determined by $(\kappa_n)_{n=1}^\infty$ and $(\nu_n)_{n=0}^\infty$.
The Haar measure is finite if and only if 
$$
\prod_{n=1}^\infty\frac{\# F_{n+1}}{\# F_n\# C_{n+1}}<\infty.
$$
\endremark

It is easy to verify that  $\Cal R$ is conservative and ergodic on the $\sigma$-finite measure space $(X,\mu)$.
This means that for each $\Cal R$-invariant subset $A\subset X$,  either  $\mu(A)=0$ or $\mu(X\setminus A)=0$.

Since the set of quasi-invariant probability measures with a fixed Radon-Nikodym derivative is a simplex \cite{GrSc}, it makes sense to introduce the following definition.

\definition{Definition 1.6} Let $\Cal S$ be a Borel countable equivalence relation on a locally compact Polish space $Z$. 
Given a Borel cocycle $\rho:\Cal S\to\Bbb R^*_+$, we say that $\Cal S$ is {\it Radon
$\rho$-uniquely ergodic} if there is a unique (up to scaling) Radon $\Cal S$-quasi-invariant measure $\lambda$ on $Z$ such that $\rho_\lambda=\rho$.
\enddefinition

\proclaim{Proposition 1.7} Let $\mu$ be a $(C,F)$-measure on $X$ determined by two sequences $(\kappa_n)_{n=1}^\infty$ and $(\nu_n)_{n=0}^\infty$ of finite measures satisfying \thetag{1-2} and \thetag{1-3}.
Then $\Cal R$ is Radon
$\rho_\mu$-uniquely ergodic.
\endproclaim
\demo{Proof} Let $\lambda$ be a Radon measure on $X$ such that $\rho_\lambda=\rho_\mu$ and  
   $\lambda(X_0)=1$.
We will prove that $\lambda=\mu$.
For that, it suffices to show that $\lambda([f]_n)=\mu([f]_n)$ for all $f\in F_n$ and $n\ge 0$.
As
$$
\mu([f]_n)=\frac{\nu_n(f)}{\nu_n(1_G)}\mu([1_G]_n)\qquad\text{and}\qquad
\lambda([f]_n)=\frac{\nu_n(f)}{\nu_n(1_G)}\lambda([1_G]_n),
$$
it is enough to prove that $\mu([1_G]_n)=\lambda([1_G]_n)$ for each $n\ge 0$.
This will be done inductively.
Of course, $\mu(X_0)=\mu([1_G]_0)=\lambda([1_G]_0)=\lambda(X_0)=1$.
Suppose that $\mu([1_G]_n)=\lambda([1_G]_n)$ for some $n$.
Then for each $c\in C_{n+1}$,
$$
\lambda([c]_{n+1})=\frac{\nu_n(c)}{\nu_n(1_G)}\lambda([1_G]_{n+1})
=\frac{\nu_n(1_G)\kappa_{n+1}(c)}{\nu_n(1_G)\kappa_{n+1}(1_G)}\lambda([1_G]_{n+1}).
$$
Since $[1_G]_n=\bigsqcup_{c\in C_{n+1}}[c]_{n+1}$, we obtain that
$$
\frac{\lambda([1_G]_{n})}{\lambda([1_G]_{n+1})}
=\frac{\sum_{c\in C_{n+1}}\lambda([c]_{n+1})}{\lambda([1_G]_{n+1})}
=\frac{\sum_{c\in C_{n+1}}\kappa_{n+1}(c)}{\kappa_{n+1}(1_G)}
=\frac{1}{\kappa_{n+1}(1_G)}
=\frac{\mu([1_G]_n)}{\mu([1_G]_{n+1})}.
$$
Hence $\lambda([1_G]_{n+1})=\mu([1_G]_{n+1})$, as desired.
\qed
\enddemo

\subhead{1.3. Nonsingular $(C,F)$-actions}
\endsubhead
Nonsingular $(C,F)$-actions were defined  in \cite{Da1}  and \cite{Da2} for Abelian groups only.
We  extend this definition to arbitrary countable  groups.
Given $g\in G$,  let 
$$
X_n^g:=\{(f_n, c_{n+1}, c_{n+2},\dots)\in X_n\mid gf_n\in F_n\}.
$$
Then $X_n^g$ is a compact open subset of $X_n$ and $X_n^g\subset X^g_{n+1}$.
  Hence the union
  $X^g:=\bigcup_{n\ge 0}X_n^g$ is an open subset of $X$.
Let $X^G:=\bigcap_{g\in G}X^g$.
 Then $X^G$ is a $G_\delta$-subset of $X$.
  Hence $X^G$ is Polish and totally disconnected  in the induced topology. 
  Given $g\in G$ and  $x\in X_G$, there is $n >0$ such that
  $x= (f_n, c_{n+1},\dots )\in X_n$ and $gf_n\in F_n$. 
  We now let $T_gx:= (gf_n, c_{n+1}, \dots)\in X_n\subset X$.
  It is straightforward  to verify that
  \roster
  \item"(i)" $T_gx\in X^G$,
  \item"(ii)" the mapping $T_g:X^G\ni x\mapsto T_gx\in X^G$ is a homeomorphism of $X^G$ and
  \item"(iii)" $T_gT_{g'}=T_{gg'}$ for all $g, g'\in G$.
\endroster
Hence $T:=(T_g)_{g\in G}$
is a continuous $G$-action on $X^G$.

\definition{Definition 1.8 \cite{Da3}} The action $T$ is called {\it the topological $(C,F)$-action of $G$ associated with $(C_n,F_{n-1})_{n=1}^\infty$.}
\enddefinition

This action is free.
The subset $X^G$ is $\Cal R$-invariant.
The $T$-orbit equivalence relation 
coincides with the restriction of $\Cal R$ to $X^G$.

\proclaim{Proposition 1.9 \cite{Da3, Proposition 1.2}} $X^G=X$ if and only if for each $g\in G$ and $n >0$, there is
$m > n$ such that
$$
gF_nC_{n+1}C_{n+2}\cdots C_m\subset F_m.
\tag1-4
$$
\endproclaim
Thus, if \thetag{1-4} holds then $T$ is a minimal continuous $G$-action on a locally compact Cantor space $X$.
Moreover, $T$ is {\it Radon uniquely ergodic}, i.e. there exists a unique  $T$-invariant Radon measure $\xi$ on $X$ such that $\xi(X_0)=1$.

From now on, $T$ is a topological $(C,F)$-action of $G$ on $X^G$ and
 $\mu$ is the nonsingular $(C,F)$-measure on $X$ determined by $(\kappa_n)_{n=1}^\infty$ and $(\nu_n)_{n=0}^\infty$ satisfying \thetag{1-2} and~\thetag{1-3}.
Since $X^G$ is $\Cal R$-invariant, we obtain that either $\mu(X^G)=0$ or $\mu(X\setminus X^G)=0$.
In the latter case $T$ is $\mu$-nonsingular, conservative and ergodic.

\proclaim{Proposition 1.10}  
The following are equivalent.
\roster
\item"(i)"  $\mu(X\setminus X^G)=0$.
\item"(ii)" 
For each $g\in G$ and every $n \ge0$,
$$
\lim_{m\to\infty}\nu_m\big((F_nC_{n+1}C_{n+2}\cdots C_m)\cap g^{-1}F_m\big)=\nu_n(F_n).
$$
\item"(iii)"
For each $g\in G$,
$$
\lim_{m\to\infty}\kappa_1*\cdots*\kappa_m( g^{-1}F_m)=1.
$$
\endroster
 \endproclaim
\demo{Proof}
(i)$\Leftrightarrow$(ii) Since $\mu(X\setminus X^G)=0$ if and only if 
$
\mu(X_n\cap X^g_m)\to\mu(X_n)$
 as $m\to\infty$ for each $g\in G$ and $n \ge 0$, it suffices to note that
 $$
 \align
 \mu(X_n\cap X^g_m)&=\mu([F_n]_n\cap[F_m\cap g^{-1}F_m]_m)\\
 &=\mu([F_nC_{n+1}\cdots C_m]_m\cap[F_m\cap g^{-1}F_m]_m)\\
 &=\mu([F_nC_{n+1}\cdots C_m\cap F_m\cap g^{-1}F_m]_m)\\
 &=\nu_m((F_nC_{n+1}\cdots C_m)\cap g^{-1}F_m)
 \endalign
 $$
 and $\mu(X_n)=\mu([F_n]_n)=\nu_n(F_n)$.
 
 (ii)$\Rightarrow$(iii)
 We set $\kappa_{1,m}:=\kappa_1*\cdots*\kappa_m$.
 Then
$$
\align
\kappa_{1,m}\big((C_1\cdots C_m)\setminus g^{-1}F_m\big) &=
\nu_m\big((F_0 C_{1} \cdots C_m) \setminus g^{-1} F_m\big)  \\
&=\nu_m(F_0 C_{1} \cdots C_m) - \nu_m\big((F_0 C_{1} \cdots C_m) \cap g^{-1} F_m\big) \\
&=\nu_0(F_0) - \nu_m\big((F_0 C_{1} \cdots C_m) \cap g^{-1} F_m\big).
\endalign
$$
Hence, 
$
\lim_{m\to\infty}\kappa_{1,m}\big((C_1\cdots C_m)\setminus g^{-1}F_m\big) = 0
$
according to (ii).
As  $\kappa_{1,m}$ is supported on $C_1\cdots C_m$, it follows that
$$
\align
\lim_{m\to\infty}\kappa_{1,m}(g^{-1}F_m) &=\lim_{m\to\infty}\kappa_{1,m}\big((C_1\cdots C_m)\cap  g^{-1}F_m\big)\\
&=\lim_{m\to\infty}\kappa_{1,m}(C_1\cdots C_m)=1,
\endalign
$$
as desired.

(iii)$\Rightarrow$(i)
Fix $g\in G$.
Take  arbitrary $n\ge 0$ and $f\in F_n$.
Then it follows from~(iii) that 
for $\mu$-a.e.  $x= (1_G, c_{n+1}, c_{n+2},)\in [1_G]_n$, 
there exists $m>0$ such that
$gfc_{n+1}\cdots c_m \in F_m$. 
This means that for $\mu$-a.e.  $y= (f_n, c_{n+1}, c_{n+2},\dots)\in X_n$, 
$$
gf_nc_{n+1}\cdots c_m \in F_m\quad\text{eventually in $m$}, 
$$
i.e. $y\in X^g$.
Hence $\mu(X\setminus X^g)=0$ and (i) follows.
 \qed
\enddemo

In the case where $\mu$ is the Haar measure for $\Cal R$, the equivalence (i)$\Leftrightarrow$(ii)
of~Proposition~1.10 was proved in~\cite{Da3}.

\proclaim{Corollary 1.11}  \roster
\item"(i)" If $\mu(X\setminus X^G)=0$ and $\mu(X)<\infty$ then $\nu_n(F_n\triangle gF_n)\to 0$ as $n\to\infty$ for each $g\in G$.
\item"(ii)" If  $\mu(X\setminus X^G)=0$, $\mu(X)<\infty$  and $\mu$ is the Haar measure for $\Cal R$ then
$G$ is amenable and 
$(F_n)_{n=1}^\infty$ is a left F{\o}lner sequence in $G$.
\item"(iii)"  If  $\mu(X\setminus X^G)=0$, $\mu(X)<\infty$, $\mu$ is the Haar measure for $\Cal R$
and 
there exists a 
subgroup $H$ of $G$ such that
$C_n\subset H$ eventually in $n$  then $H$ is of finite index in $G$.
\endroster
\endproclaim
\demo{Proof} (i) We note that $\nu_n(F_n)=\mu([F_n]_n)=\mu(X_n)\to\mu(X)$ as $n\to\infty$.
Hence it follows from Proposition~1.10(ii) that for each $\epsilon>0$, there is $n>0$ such that if $m>n$ then
$$
\aligned
\nu_m\big(F_nC_{n+1}C_{n+2}\cdots C_m\big)&>(1-\epsilon)\nu_m(F_m)\quad\text{and}\\
\nu_m\big((F_nC_{n+1}C_{n+2}\cdots C_m)\cap gF_m\big)&>(1-\epsilon)\nu_m(F_m).
\endaligned
$$
Hence $\nu_m(F_m\cap gF_m)>(1-2\epsilon)\nu_m(F_m)$.
It follows that $\lim_{m\to\infty}\nu_m(F_m\triangle gF_m)=0$, as desired.

(ii) Since $\mu$ is the Haar measure for $\Cal R$, it follows that $\nu_n (A)=\frac{\# A}{\#C_1\cdots\# C_n}$
for each subset $A\subset F_n$.
Since $\mu$ is finite, there exists a limit 
$$
\lim_{n\to\infty}\frac{\# F_n}{\#C_1\cdots\# C_n}=\mu(X).
$$
This fact and  (i) yield that for each $g\in G$,
$$
0=\lim_{m\to\infty}\nu_n(F_n\triangle gF_n)=\lim_{m\to\infty}\frac{\#(F_n\triangle gF_n)}{\#C_1\cdots\# C_n}
=\mu(X)\lim_{m\to\infty}\frac{\#(F_n\triangle gF_n)}{\#F_n}.
$$
Hence $(F_n)_{n=1}^\infty$ is a left F{\o}lner sequence in $G$.
Therefore $G$ is amenable.

(iii) Suppose that 
$H$ is  of infinite index in $G$.
We first prove an auxiliary claim.

{\it Claim A.} For each finite subset  $S\subset G$, there exists an element $g\in G$ such that $g\not\in \bigcup_{a,b\in S}aHb^{-1}$.

By~(ii), $G$ is amenable.
Hence there exists a left-invariant finitely additive measure $\xi$ on the $\sigma$-algebra of all subsets of $G$ such that $\xi(G)=1$.
We first observe that since $H$ is of infinite index, $\xi(H)=0$.
Indeed, for each $n>0$, there are elements $g_1,\dots,g_n\in G$ such that the cosets $g_1H,\dots,g_nH$ are mutually disjoint.
Hence 
$$
1\ge \xi\bigg(\bigsqcup_{j=1}^ng_jH\bigg)=\sum_{j=1}^n\xi(g_jH)=\sum_{j=1}^n\xi(H)=n\xi(H).
$$
This yields that $\xi(H)=0$, as desired.
As $g^{-1}Hg$ is also a subgroup of infinite index in $G$, 
it follows that  $\xi(g^{-1}Hg)=0$  for each $g\in G$.
Since $\xi$ is left-invariant, $\xi(kHg)=0$ for all $k,g\in G$.
This implies that  
$
\xi\big(\bigcup_{a,b\in S}aHb^{-1}\big)=0.
$
Therefore, $G\ne\bigcup_{a,b\in S}aHb^{-1}$.
Thus, Claim~A is proved.

\comment
Indeed, since
$H$ is virtually normal, there is a cofinite subgroup $N\subset H$ which is normal in $G$.
Select a finite subset $B\subset H$ such that $BN=H$.
Then 
$$
\bigcup_{a,b\in S}aHb^{-1}=\bigcup_{a,b\in S}aBNb^{-1}=\Bigg(\bigcup_{a,b\in S}aBb^{-1}\Bigg)N.
$$
Since $H$ is of infinite index in $G$, it follows that $N$ is of infinite index in $G$.
As the subset in parenthesis is finite,  we obtain that  $\bigcup_{a,b\in S}aHb^{-1}\ne G$.
Thus, the claim is proved.

\endcomment

Since $\mu(X)<\infty$, 
there exists $n$ such that  $\nu_n(F_n)>0.5\nu_m(F_m)$  and $C_m\subset H$   for each $m\ge n$.
Hence
$$
\nu_m\big((F_nH)\cap F_m\big)\ge\nu_m\big((F_nC_{n+1}C_{n+2}\cdots C_m)\cap F_m\big)=\nu_n(F_n)>0.5\nu_m(F_m)
$$
for each $m\ge n$.
By Claim~A, there is $g\in G$ such that $gF_nH\cap F_nH=\emptyset$.
Since $\mu$ is the Haar measure, it follows that
$$
\align
\nu_m\big((gF_nH)\cap F_m\big)&\ge 
\nu_m\big((gF_nC_{n+1}C_{n+2}\cdots C_m)\cap F_m\big)\\
&=
\nu_m\big((F_nC_{n+1}C_{n+2}\cdots C_m)\cap g^{-1}F_m\big).
\endalign
$$
This inequality and (i) yield that $\nu_m\big((gF_nH)\cap F_m\big)>0.5\nu_m(F_m)$
eventually in $m$.
Therefore $\nu_m(F_nH\cap gF_nH)>0$ eventually in $m$, a contradiction.
 \qed
\enddemo

\comment

By Proposition~1.6, 
$$
\align
&\nu_m\big((F_nC_{n+1}C_{n+2}\cdots C_m)\cap F_m\big)\to \nu_n(F_n)\quad\text{and}\\
&\nu_m\big((F_nC_{n+1}C_{n+2}\cdots C_m)\cap g^{-1}F_m\big)\to\nu_n(F_n)
\endalign
$$
as $m\to\infty$.
Since $T$ preserves $\mu$, we obtain that
$$
\nu_m\big((F_nC_{n+1}C_{n+2}\cdots C_m)\cap g^{-1}F_m\big)=\nu_m\big((gF_nC_{n+1}C_{n+2}\cdots C_m)\cap F_m\big)
$$
for every $m$.
Hence there is $m>n$ such that
$$
\align
&\nu_m\big((F_nH)\cap F_m\big)\ge\nu_m\big((F_nC_{n+1}C_{n+2}\cdots C_m)\cap F_m\big)>0.5\nu_m(F_m)\quad\text{and}\\
&\nu_m\big((gF_nH)\cap F_m\big)\ge 
\nu_m\big((gF_nC_{n+1}C_{n+2}\cdots C_m)\cap F_m\big)
>0.5\nu_m(F_m).
\endalign
$$
It follows that $\nu_m(F_nH\cap gF_nH )>0$, a contradiction.

$$
\lim_{m\to\infty}\frac{\#(gF_nC_{n+1}\cdots C_m\cap F_m)}{\#(F_nC_{n+1}\cdots C_m)}=1.
$$
Since $T$ is finite measure preserving,  there is $N$ such that if $m>n>N$ then
$$
\frac{\#(F_m\setminus (F_nC_{n+1}\cdots C_m))}{\#(F_nC_{n+1}\cdots C_m)}<\frac1{10}.
$$
Then there is $g\in G$ such that  $gF_nC_{n+1}\cdots C_m\cap (F_nC_{n+1}\cdots C_m)=\emptyset$
for all $m>n$, a contradiction.

\qed

\endcomment

\definition{Definition 1.12} If $\mu(X\setminus X^G)=0$ then the dynamical system $(X,\mu, T)$ (or simply $T$) is called
{\it the nonsingular $(C,F)$-action associated with  $(C_n,F_{n-1},\kappa_n,\nu_{n-1})_{n=1}^\infty$.}
\enddefinition

From now on we consider only the case where $\mu(X\setminus X^G)=0$.
As $X=X^G$ mod~0,  we will assume that $T$ is defined on the entire space $X$.
Then for each $n$ and every two elements $g,h\in F_n$, 
we have that $T_{hg^{-1}}[g]_n=[h]_n$ and the Radon-Nikodym derivative of the transformation $T_{hg^{-1}}$ is constant on the subset $[g]_n$.
More precisely, this constant equals $\frac{\nu_n(h)}{\nu_n(g)}$.

We now prove the main result of this section.

\proclaim{Theorem 1.13} Each nonsingular $(C,F)$-action  is of rank one.
Conversely, each rank-one nonsingular $G$-action is  isomorphic (via a   measure preserving isomorphism) to a $(C,F)$-action.
\endproclaim
\demo{Proof} 
Let a sequence $(C_n,F_{n-1},\kappa_n,\nu_{n-1})_{n=1}^\infty$ satisfy
\thetag{1-1}--\thetag{1-3} and Proposit\-ion~1.10(ii).
We claim that 
 the $(C,F)$-action $T=(T_g)_{g\in G}$ associated with this sequence
is of rank one along  $(F_n)_{n=0}^\infty$.
Let $X$ be the space of this action and let $\mu$ be the $(C,F)$-measure on $X$ determined by
$(\kappa_n)_{n=1}^\infty$ and
$(\nu_n)_{n=0}^\infty$.
Then $X=\bigcup_{n\ge 0}X_n$ and $\Cal R=\bigcup_{n\ge 0}\Cal S_n$, where 
$X_n$ and $\Cal S_n$ were introduced in \S1.2.
Of course, for  each $n\in\Bbb N$, the pair $([1_G]_n,F_n)$ is  a Rokhlin tower for $T$.
Moreover,  
\roster
\item"a)" $X_{[1_G]_n,F_n}=X_n$, 
\item"b)" $\xi_{[1_G]_n,F_n}$ is the partition of $X_n$ into cylinders $[f]_n$, $f\in F_n$, and
\item"c)" if $x=(f_n,c_{n+1},\dots)\in X_n\cap X^G$ then $O_{[1_G]_n,F_n}(x)=\{T_gx\mid g\in {F_nf_n^{-1}}\}=\Cal S_n(x)$.
\endroster
We note that a) and b) imply  that  Definition~1.1(i) holds.
It follows from Proposition~1.10 that for a.e. $x\in X$ (or, more precisely, for each $x\in X^G$), the $T$-orbit of $x$ equals $\Cal R(x)$.
As $\Cal R(x)=\bigcup_{n=1}^\infty\Cal S_n(x)$, 
it follows that  c) implies (ii) from  Definition~1.1.
Hence $T$ is of rank one along $(F_n)_{n=1}^\infty$.

Conversely, suppose that $T$ is a nonsingular $G$-action of rank one along an increasing sequence
$(Q_n)_{n=0}^\infty$ of finite subsets in $G$ with $Q_0=\{1_G\}$.
Let $(B_n,Q_n)_{n=0}^\infty$ be the corresponding generating sequence of Rokhlin towers such that (i) and (ii) of Definition~1.1 hold.
We have to define a sequence $(C_n,F_{n-1},\kappa_n,\nu_{n-1})_{n=1}^\infty$, 
satisfying~\thetag{1-1}--\thetag{1-3} and Proposition~1.10(ii)
such that the associated $(C,F)$-action   is isomorphic to $T$.
We first set 
$F_n:=Q_n$ for each $n\ge 0$.
By Definition~1.1(i), for each $n\ge 0$, there is a subset $R_{n+1}\subset Q_{n+1}$ such that
$B_n=\bigsqcup_{f\in R_{n+1}}T_f B_{n+1}$.
Without loss of generality we may assume that $1_G\in R_{n+1}$.
Indeed, if this is not the case, we replace $(B_{n+1},Q_{n+1})$ with another Rokhlin tower
$(T_sB_{n+1}, Q_{n+1}s^{-1})$ for an element $s\in R_{n+1}$.
Then $\xi_{B_{n+1},Q_{n+1}}=\xi_{T_sB_{n+1}, Q_{n+1}s^{-1}}$ and $O_{B_{n+1},Q_{n+1}}(x)=
O_{T_sB_{n+1}, Q_{n+1}s^{-1}}(x)$ for each $x\in X$.
Hence such replacements will not affect  (i) and (ii) of Definition~1.1.
We now set $C_{n+1}:=R_{n+1}$.
Thus, we defined the entire sequence $(C_n,F_{n-1})_{n\ge 1}$.
It is straightforward to verify that ~\thetag{1-1} holds.
Let $\nu_0$ be the  Dirac measure  supported at $1_G$.
Next, for each $n> 0$ and $f\in F_n$, we let $\nu_n(f):=\mu(T_fB_n)$.
Thus, we obtain a non-degenerated measure  $\nu_n$  on $F_n$.
Finally, we define a probability $\kappa_{n+1}$ on $C_{n+1}$ by setting
$$
\kappa_{n+1}(c):=\frac{\nu_{n+1}(c)}{\nu_n(1_G)}\quad\text{for each $c\in C_{n+1}$ and $n\ge 0$.}
$$
Thus,  the entire sequence of measures  $(\kappa_n,\nu_{n-1})_{n\ge 1}$ is defined.
It follows from  the property (iii)  which is below Definition~1.1  that
$$
\frac{\nu_{n+1}(fc)}{\nu_n(f)}=\frac{\nu_{n+1}(c)}{\nu_n(1_G)}=\kappa_{n+1}(c)\quad\text{for each $c\in C_{n+1}$ and $f\in F_n$,}
$$
i.e. \thetag{1-3} holds.
We note that for each $n>0$, the restrictions of $\xi_{B_n,Q_n}$ to the subset $X_{B_0,Q_0}=B_0$
is the finite partition of $B_0$ into subsets $T_{c_1\cdots c_n}B_{n}$, where $(c_1,\dots,c_n)$ runs the subset $C_1\times\cdots \times C_n$.
As $\xi_{B_n,Q_n}\restriction B_0$ converges to the  partition into singletons, we obtain that
$$
\max_{c_1\in C_1,\dots, c_n\in C_n}\mu(T_{c_1\cdots c_n}B_{n})\to 0\quad\text{as $n\to\infty$.}
$$ 
Since
$
\mu(T_{c_1\cdots c_n}B_{n})=\mu(X_{B_0})\prod_{j=1}^n\kappa_j(c_j),
$
 \thetag{1-2} follows.
 
 Fix $n>0$, $g\in G$ and $\epsilon>0$.
It follows from Definition~1.1(ii) that there exists $M>n$ such that for each $m>M$ there is a subset $A\subset X_{B_n,Q_n}$  such that
$\mu( X_{B_n,Q_n}\setminus A)<\epsilon$ and $T_gx\in O_{B_m,Q_m}(x)$ for each $x\in A$.
Hence there exist $f_1,f_2\in Q_m$ such that $T_gx=T_{f_1f_2^{-1}}x$ and $T_{f_2^{-1}}x\in B_m$.
As $T$ is free, $gf_2=f_1\in Q_m=F_m$.
It follows  that $T_{gf_2}B_m=T_{f_1}B_m\subset X_{B_m,Q_m}$.
 Since $T_{f_2}B_m\ni x$ and $x\in X_{B_n,Q_n}$, we obtain that  $T_{f_2}B_m\subset X_{B_n,Q_n}$ because $\xi_{B_m,Q_m}$ is finer than $\xi_{B_n,Q_n}$.
 Thus, without loss of generality we may assume that $A$ is measurable with respect to the partition $\xi_{B_m,Q_m}$.
 Since 
 $$
 X_{B_n,Q_n}=\bigsqcup_{f\in F_n} T_fB_n=\bigsqcup_{f\in F_n} T_f\left(\bigsqcup_{c\in C_{n+1}\cdots C_m}T_cB_m\right),
 $$
it follows  that $T_{f_2}B_m\subset X_{B_n,Q_n}$ if and only if $f_2\in F_nC_{n+1}\cdots C_m$.
 Hence
 $$
\nu_m(\{f_2\in F_m\mid gf_2\in F_m, f_2\in F_nC_{n+1}\cdots C_m\})\ge \mu(A)>\mu(X_{B_n,Q_n})-\epsilon.
 $$
This implies Proposition~1.10(ii).

Thus, $(C_n,F_{n-1},\kappa_n,\nu_{n-1})_{n=1}^\infty$, 
satisfies~\thetag{1-1}--\thetag{1-3} and Proposition~1.10(ii).
Denote by  $R$ the nonsingular $(C,F)$-action associated with $(C_n,F_{n-1},\kappa_n,\nu_{n-1})_{n=1}^\infty$.
Let $(Y,\nu)$ be the space of this action.
The correspondence 
$$
T_fB_n \longleftrightarrow R_f[1_G]_n,\quad\text{where $f$ runs $F_n$ and $n$ runs $\Bbb N$,}
$$
gives rise to a Boolean measure preserving isomorphism of the underlying algebras of measurable subsets
on $X$ and $Y$.
The Boolean isomorphism is generated by a
certain  pointwise measure preserving isomorphism $\theta$ of $(X,\mu)$ onto $(Y,\nu)$.
We claim that $\theta$ intertwines $T$ with $R$.
Indeed, take $g\in G$.
As was shown above, for each $n>0$ and $\epsilon>0$, there exists $M>n$ such that for each $m>M$  there is a subset  $Q'\subset Q_m$ such that $\bigsqcup_{f\in Q'}T_fB_m\subset X_{B_n,Q_n}$, $\mu(X_{B_n,Q_n}\setminus \bigsqcup_{f\in Q'}T_fB_m)<\epsilon$ and $gQ'\subset Q_m$.
Hence
$$
\theta(T_gT_fB_m)=R_g\theta(T_fB_m)\quad\text{for all $f\in Q'$.}
$$
Passing to the limit as $\epsilon\to0$ and using the fact that $\nu\circ\theta=\mu$, we obtain that
$\theta(T_gx)=R_g\theta x$ for a.e. $x\in X_n$.
Since $n$ is arbitrary, $\theta T_g=R_g\theta$, as claimed.
\qed
\enddemo

\remark{Remark 1.14}  
\roster
\item"(a)" Theorem~1.13 corrects \cite{Da3, Theorem~2.6}, where  the particular case of  $\sigma$-finite measure preserving rank-one actions was under consideration:
the condition~(ii) (see  Definition~1.1) is missing in the definition of rank one in \cite{Da3}.
However, this condition can not be omitted:
 counterexamples of non-rank-one action satisfying (i) (and hence not satisfying (ii)) is provided in  Examples~3.4--3.6 below.
\item"(b)" It is worth noting that the condition on the Radon-Nikodym derivatives in the definition of Rokhlin tower in \S1.1 is important and can not be omitted either.
Indeed, the associated flow of each nonsingular $(C,F)$-system is AT in the sense of Connes-Woods \cite{CoWo} (see also \cite{Ha}).
In \cite{DoHa}, A.H.~Dooley and T.~Hamachi constructed explicitly a Markov nonsingular odometer  ($\Bbb Z$-action) whose associated
flow is non-AT. 
Hence, this Markov odometer is not isomorphic (if fact, it is not even orbit equivalent) to any
rank-one nonsingular $\Bbb Z$-action.\footnote{We do not provide definitions of orbit equivalence, associated flow and AT-flow  because we will not use it anywhere below in this paper. 
Instead,  we refer the interested reader to the survey \cite{DaSi}.}
On the other hand, it is easy to see that this odometer satisfies a ``relaxed version'' of Definition~1.1 in which we  drop only the condition on the Radon-Nikodym derivatives.
\endroster
\endremark

\subhead{1.4.} Telescopings and reductions
\endsubhead
Let a sequence $\Cal T=(C_n,F_{n-1},\kappa_n,\nu_{n-1})_{n=1}^\infty$ satisfy   \thetag{1-1}--\thetag{1-3} and Proposition~1.10(ii).
 Denote by $T=(T_g)_{g\in G}$ the $(C,F)$-action of $G$ on $X$ associated with $\Cal T$.
 Let $\mu$ stand for the 
 the nonsingular $(C,F)$-measure on $X$ determined by $(\kappa_n)_{n=1}^\infty$ and $(\nu_n)_{n=0}^\infty$.

Given a strictly increasing  infinite sequence of integers $\boldsymbol l=(l_n)_{n=0}^\infty$ such that
$l_0=0$, we let 
$$
\widetilde F_n:=F_{l_n},\quad \widetilde C_{n+1}:=C_{l_n+1}\cdots C_{l_{n+1}},\quad \widetilde \nu_n:=\nu_n,\quad\widetilde\kappa_{n+1}:=\kappa_{l_n+1}*\cdots*\kappa_{l_{n+1}}.
$$
for each $n\ge 0$.

\definition{Definition  1.15}
We call the sequence $\widetilde{\Cal T}:=(\widetilde C_n,\widetilde F_{n-1},\widetilde\kappa_n,\widetilde\nu_{n-1})_{n=1}^\infty$ the {\it $\boldsymbol l$-tele\-scoping} of $\Cal T$.
\enddefinition

It is easy to check that $\widetilde{\Cal T}$   satisfies  \thetag{1-1}--\thetag{1-3} and Proposition~1.10(ii).
Hence a nonsingular $(C,F)$-action $\widetilde T=(\widetilde T_g)_{g\in G}$ of $G$ 
 associated with $\widetilde{\Cal T}$
is well defined.
Let $\widetilde X$ denote the  space of $\widetilde T$ and let $\widetilde\mu$
denote the nonsingular $(C,F)$-measure on $\widetilde X$ determined by $(\widetilde\kappa_n)_{n=1}^\infty$ and $(\widetilde\nu_n)_{n=0}^\infty$.
There is a {\it canonical} measure preserving isomorphism $\iota_{\boldsymbol l}$ of  $(X,\mu)$ onto $(\widetilde X,\widetilde\mu)$ that intertwines
$T$ with $\widetilde T$.
Indeed, if $x\in X$ then we select the smallest  $n\ge 0$ such that $x=(f_{l_n},c_{l_n+1},c_{l_n+1},\dots)\in X_{l_n}$.
Let
$$
\iota_{\boldsymbol l}(x):=(f_{l_n},c_{l_n+1}\cdots c_{l_{n+1}}, c_{l_{n+1}+1}\cdots c_{l_{n+2}},\dots)\in \widetilde X_n\subset\widetilde X,
$$
where $\widetilde X_n=\widetilde F_n\times\widetilde C_{n+1}\times \widetilde C_{n+2}\times\cdots$.
It is a routine to verify that $\iota_{\boldsymbol l}$ is a homeomorphism of $X$ onto $\widetilde X$ such that  $\iota_{\boldsymbol l} T_g=\widetilde T_g\iota_{\boldsymbol l}$ for each $g\in G$, as desired.

  Let  $\boldsymbol l=(l_n)_{n=0}^\infty$  and  $\boldsymbol m=(m_n)_{n=0}^\infty$ be two strictly increasing sequences of integers such that
$l_0=m_0=0$.
If $\widetilde {\Cal T}$ is the   $\boldsymbol l$-telescoping of  ${\Cal T}$
and $\Cal S$ is the $\boldsymbol m$-telescoping of $\widetilde{\Cal T}$ then $\Cal S$ is the $\boldsymbol l\circ\boldsymbol m$-telescoping of $\Cal T$, where $\boldsymbol l\circ\boldsymbol m:=(l_{m_n})_{n=1}^\infty$
and $\iota_{\boldsymbol m}\circ \iota_{\boldsymbol l}=\iota_{\boldsymbol l\circ\boldsymbol m}$.

Given a sequence $\boldsymbol A=(A_n)_{n=1}^\infty$ of subsets $A_n\subset C_{n}$ such that
$1_G\in A_n$ for each $n\in\Bbb N$ and $\sum_{n=1}^\infty(1- \kappa_n(A_n))<\infty$, we let
$$ 
 \kappa'_n(a):=\frac{\kappa_n(a)}{\kappa_n(A_n)},\quad a\in A_n.
$$
Then $ \kappa'_n$ is a non-degenerated probability on $A_n$ for each $n\in\Bbb N$.
We also define a measure $\nu_n'$ on $F_n$ by setting 
$$
\nu_n'=\frac{1}{\prod_{j=1}^n\kappa_j(A_j)}\cdot \nu_n
$$
if $n>0$ and $\nu_0':=\nu_0$.
Let $\Cal T':=(A_n,F_{n-1},\kappa_n',\nu_{n-1}')_{n=1}^\infty$.

\definition{Definition 1.16}
We call $\Cal T'$   an  {\it $\boldsymbol A$-reduction} of $\Cal T$.
\enddefinition
It is easy to check that $\Cal T'$   satisfies  \thetag{1-1}--\thetag{1-3}.
We note that 
$$
(\kappa_1'*\cdots*\kappa'_m)(G\setminus g^{-1}F_m)\le
\frac {(\kappa_1*\cdots*\kappa_m)(G\setminus g^{-1}F_m)}{\prod_{j=1}^m\kappa(A_j)}
$$
for each $g\in G$.
Passing to the limit and using 
 Proposition~1.10(iii) for $\Cal T$,  
 we obtain that $(\kappa_1'*\cdots*\kappa'_m)(G\setminus g^{-1}F_m)\to 0$ as $m\to\infty$.
 In other words, 
 $$
 \lim_{m\to\infty}(\kappa_1'*\cdots*\kappa'_m)( g^{-1}F_m)=1, 
 $$
i.e.  
Proposition~1.10(iii) holds for $\Cal T'$.
\comment

 \thetag{1-6} is equivalent to 
$$
\lim_{m\to\infty}\nu_m\big((F_nC_{n+1}\cdots C_m)\setminus g^{-1}F_m\big)=0\qquad\text{for each $n>0$.}\tag1-8
$$
Since \thetag{1-6} and hence \thetag{1-8} hold for $\Cal T$, it follows from the inequality
$$
\nu_m'\big((F_nA_{n+1}\cdots A_m)\setminus g^{-1}F_m\big)\le \frac{\nu_m\big((F_nC_{n+1}\cdots C_m)\setminus g^{-1}F_m\big)}{\prod_{j=1}^m\kappa_j(A_j)}
$$
that \thetag{1-8} holds for $\Cal T'$.
Therefore,  \thetag{1-6} holds for $\Cal T'$.

\endcomment
Hence, a nonsingular $(C,F)$-action $ T'=(T_g')_{g\in G}$ of $G$ 
 associated with $\Cal T'$
is well defined.
Let $X'$ denote the  space of $T'$ and let $\mu'$
denote the nonsingular $(C,F)$-measure on $X'$ determined by $(\kappa_n')_{n=1}^\infty$ and $(\nu_n')_{n=0}^\infty$.

\proclaim{Proposition 1.17}
There is a canonical measure scaling isomorphism $\iota_{\boldsymbol A}$ of  $(X,\mu)$ onto $(X',\mu')$ that intertwines
$T$ with $T'$
and  $\mu\circ \iota_{\boldsymbol A}^{-1}=\prod_{m>0}{\kappa_m(A_m)}\cdot\mu'$.
\endproclaim
\demo{Proof}
Indeed,
fix $n\in\Bbb N$.
Since $\mu\restriction X_n=\nu_n\otimes\kappa_{n+1}\otimes\kappa_{n+2}\otimes\cdots$, it follows from the Borel-Cantelli lemma that for a.e. $x=(f_n,c_{n+1},c_{n+2},\dots)\in X_n$, there is
$N=N_x>n$ such that $c_m\in A_m$ for each $m>N$.
We then let
$$
\iota_{\boldsymbol A,n}(x):=(f_nc_{n+1}\cdots c_N, c_{N+1},c_{N+2},\dots)\in X_N'\subset X'.
$$
It is routine to verify that $\iota_{\boldsymbol A,n}:X_n\ni x\mapsto \iota_{\boldsymbol A,n}(x)\in X'$ is a well defined nonsingular 
 mapping and
$$
\frac{d\mu'\circ \iota_{\boldsymbol A,n}}{d\mu}(x)=\frac{\nu_N'(f_nc_{n+1}\cdots c_N)}{\nu_N(f_nc_{n+1}\cdots c_N)}\prod_{m>N}\frac{1}{\kappa_m(A_m)}=
\prod_{m>0}\frac{1}{\kappa_m(A_m)}.
$$
Moreover,
$\iota_{\boldsymbol A,n+1}\restriction X_n=\iota_{\boldsymbol A,n}$ for each $n\in\Bbb N$.
Hence a measurable mapping $\iota_{\boldsymbol A}:X\to X'$ is a well defined
by the restrictions $\iota_{\boldsymbol A}\restriction X_n=\iota_{\boldsymbol A,n}$ for all $n\in\Bbb N$.
It is straightforward to verify that $\iota_{\boldsymbol A}$ is an  isomorphism of $(X,\mu)$ onto $(X',\mu')$ with  $\mu\circ \iota_{\boldsymbol A}^{-1}=\prod_{m>0}{\kappa_m(A_m)}\cdot\mu'$ and
 $\iota_{\boldsymbol A}T_g=T_g'\iota_{\boldsymbol A}$ for each $g\in G$.
 \qed
\enddemo

\subhead 1.5. Locally compact  models for rank-one nonsingular systems
\endsubhead
Let $Z$ be  a locally compact Polish $G$-space.
We remind that a Borel mapping $\rho:G\times Z\to\Bbb R^*_+$ is called a $G$-cocycle if
$$
\rho(g_2,g_1z)\rho(g_1,z)=\rho(g_2g_1,z)\qquad\text{for all $g_1,g_2\in G$ and $z\in Z$}.
$$
The following definition is a dynamical analogue of Definition~1.6.

\definition{Definition 1.18}
Fix a $G$-cocycle $\rho$.
We say that the $G$-action on $Z$ is Radon {\it $\rho$-uniquely ergodic} if there exists a unique (up to scaling) Radon $G$-quasiinvariant measure $\gamma$ on $Z$ such that
$$
\frac{d\gamma\circ g}{d\gamma}(z)=\rho(g,z)\quad\text{for all $g\in G$ and $z\in Z$}.
$$
\enddefinition

We now show that each rank-one nonsingular action has a uniquely ergodic continuous realization on a locally compact Cantor space.

\proclaim{Theorem 1.19} Let a nonsingular action $R$ of $G$ on a $\sigma$-finite standard non-atomic measure space $(Z,\eta)$ be of rank one along a sequence $(Q_n)_{n=1}^\infty$.
Then there exist
\roster
\item"(i)"
a continuous, minimal,  Radon uniquely ergodic $G$-action $T'=(T_g')_{g\in G}$
defined on a locally compact Cantor space $X'$,
\item"(ii)" 
a $T'$-quasi-invariant Radon measure $\mu'$ on $X'$ such that
the Radon-Nikodym derivative  $\frac{d\mu'\circ T_g'}{d\mu'}:X'\to\Bbb R^*_+$ 
is continuous for each $g\in G$,
\item"(iii)" 
a measure 
preserving Borel isomorphism of $(Z,\eta)$ onto $(X',\mu')$  that intertwines $R$ with $T'$,
\item"(iv)"
a sequence $\Cal T'=(C_n',F_{n-1}',\kappa_n',\nu_{n-1}')_{n=1}^\infty$ satisfying \thetag{1-1}--\thetag{1-4} such that $(X',\mu',T')$ is the the $(C,F)$-action associated with $\Cal T'$ and
\item"(v)"  a sequence  $(z_n)_{n=1}^\infty$ such that $z_n\in Q_n$ for each $n$ and 
$(F_n')_{n=1}^\infty$ is a subsequence of $(z_n^{-1}Q_n)_{n=1}^\infty$.
\endroster
Moreover, $T'$ is Radon  $\Big(\frac{d\mu'\circ T_g'}{d\mu'}\Big)_{g\in G}$-uniquely ergodic.
\endproclaim

\demo{Proof} By Theorem~1.13, there is a 
a sequence $\Cal T=(C_n,F_{n-1},\kappa_n,\nu_{n-1})_{n=1}^\infty$ satisfying \thetag{1-1}--\thetag{1-3} and Proposition~1.10(ii) such that the $(C,F)$-action $T$ of $G$ associated with $\Cal T$ is isomorphic to $R$ via a measure preserving isomorphism.
Denote by $(X,\mu)$  the space of $T$.
Let $G=\{g_j\mid j\in\Bbb N\}$.
It follows from Proposition~1.10(ii) that there are an increasing  sequence $\boldsymbol l=(l_n)_{n=0}^\infty$ of integers 
and a sequence $(D_n)_{n=1}^\infty$ of subsets in $G$ such that
$l_0=0$, $D_{n+1}:= \big(C_{l_n+1}\cdots C_{l_{n+1}}\big)\cap\bigcap_{j=1}^n g_j^{-1}F_{l_{n+1}}$ and
$$
\kappa_{l_n+1}*\cdots*\kappa_{l_{n+1}}(D_{n+1})>1-\frac1{(n+1)^2}
$$
for each $n\ge 0$.
Denote by $\widetilde{ \Cal T}=(\widetilde C_n,\widetilde F_{n-1},\widetilde \kappa_n,\widetilde \nu_{n-1})_{n=1}^\infty$ the $\boldsymbol l$-telescoping of $\Cal T$.
Let
 $(\widetilde X,\widetilde\mu, \widetilde T)$  stand for
the $(C,F)$-action of $G$ associated with $\widetilde{ \Cal T}$.
Then $D_n\subset\widetilde C_n$ and $\widetilde\kappa_n(D_n)>1-n^{-2}$ for each $n>0$.
Denote by
 $\iota_{\boldsymbol l}$  the canonical measure preserving isomorphism intertwining  $T$ with $\widetilde T$.
 In general, $1_G\not\in D_n$.
 Therefore we need to modify the $(C,F)$-parameters $\widetilde{\Cal T}$.
 First, we choose, for each $n>0$, an element $c_n\in D_n$.
 Then we let $z_0:=1_G$ and $z_n:=c_1\cdots c_n$ for each $n>0$.
 Finally, we define a new sequence  $\widehat{ \Cal T}=(\widehat C_n,\widehat F_{n-1},\widehat \kappa_n,\widehat \nu_{n-1})_{n=1}^\infty$ by setting:
 $$
 \widehat C_n:=z_{n-1}\widetilde C_nz_n^{-1}, \  \widehat F_{n-1}:=\widetilde F_{n-1}z_{n-1}^{-1},
 $$
$\widehat \kappa_n$ is the image of $\kappa_n$ under the bijection $\widetilde C_n\ni c\mapsto z_{n-1}cz_n^{-1}\in  \widehat C_n $ and $\widehat \nu_{n-1}$ is the image of $\nu_{n-1}$ under the bijection $\widetilde F_{n-1}\ni f\mapsto fz_{n-1}^{-1}\in  \widehat F_{n-1}$.
It is straightforward to verify that $\widehat{ \Cal T}$ satisfies  
\thetag{1-1}--\thetag{1-3} and Proposition~1.10(ii).
Denote by $(\widehat X,\widehat \mu,\widehat T)$ the $(C,F)$-action of $G$ associated with $\widehat{ \Cal T}$.
Then there is a canonical continuous measure preserving isomorphism $\vartheta: (\widetilde X,\widetilde\mu)\to( \widehat X,\widehat \mu)$ that intertwines $\widetilde T$ with $\widehat T$:
$$
\widetilde X\supset \widetilde X_n\ni(f_n,c_{n+1},\dots)\mapsto(f_nz_{n}^{-1}, z_{n}c_{n+1}z_{n+1}^{-1},
z_{n+1}c_{n+2}z_{n+2}^{-1},\dots)\in\widehat X_n\subset\widehat X.
$$
Let $\widehat D_n$ is the image of $D_n$ under the bijection $\widetilde C_n\ni c\mapsto z_{n-1}cz_n^{-1}\in  \widehat C_n$.
Then $1\in \widehat D_n$ and $\widehat\kappa_n(\widehat D_n)>1-n^{-2}$ 
for each $n\in \Bbb N$.
Hence $\sum_{n=1}^\infty(1-\widehat\kappa_n(\widehat D_n))<\infty$. 
We now set $\boldsymbol D:=(D_n)_{n=1}^\infty$.
Denote by $\Cal T'$ the $\boldsymbol D$-reduction of $\widehat{ \Cal T}$.
Then  $\Cal T'$ satisfies not only \thetag{1-1}--\thetag{1-3} and Proposition~1.10(ii)
but also \thetag{1-4}.
Let $( X',\mu', T')$ denote the $(C,F)$-action of $G$ associated with $\Cal T'$ and let
$\iota_{\boldsymbol D}$ stand for the canonical measure scaling isomorphism of $(\widehat X,\widehat\mu)$ onto $( X',\mu')$ that intertwines $\widehat T$
with $T'$.
Then $\iota_{\boldsymbol D}\circ \vartheta\circ \iota_{\boldsymbol l}$ is a measure scaling isomorphism
 of $(X,\mu, T)$ onto $(X',\mu',T')$.
 Replacing $\mu'$ with $a\cdot\mu'$ for an appropriate $a>0$ we obtain that 
 $\iota_{\boldsymbol D}\circ \vartheta\circ \iota_{\boldsymbol l}$ is measure preserving.
 It follows from Proposition~1.9 that $T'$ is a Radon uniquely ergodic minimal  continuous action of $G$ on the locally compact Polish space $X'$.
Thus, we proved (i), (iii), (iv) and~(v).
The property (ii) follows easily from \thetag{1-4} and the definition of $\mu'$.

The final claim of the theorem follows from (iv) and Proposition~1.7.
 \qed
\enddemo

\comment
It is worthy to note that  (v) implies the following:  if $G=\Bbb Z^d$ and $T$ is of rank one along a sequence of parallelepipeds or cubes in $\Bbb Z^d$ (this corresponds to the classical concept of rank one for $\Bbb Z^d$-actions)  then $T'$ is also of rank one along parallelepipeds or cubes respectively. ????????

\endcomment

\subhead 1.6.  Nonsingular $\Bbb Z$-actions of rank one along intervals and $(C,F)$-const\-ruction
\endsubhead
Let $G=\Bbb Z$.
Suppose that a sequence
 $\Cal T=(C_n,F_{n-1},\kappa_n,\nu_{n-1})_{n=1}^\infty$ satisfies \thetag{1-1}--\thetag{1-3} and Proposition~1.10(ii)
and there is a sequence $(h_n)_{n=0}^\infty$ of positive integers such that
$F_n=\{0,1,\dots,h_n-1\}$.
Denote by $(X,\mu,T)$ the $(C,F)$-dynamical system associated with $\Cal T$.

We now show how to obtain $(X,\mu,T)$ via the classical inductive geometric cutting-and-stacking  in the case of $\Bbb Z$-actions.
On the initial step of the  construction we define a column
 $Y_0$ consisting of a single interval $[0,1)$ equipped with Lebesgue measure.
Assume that at the $n$-th step we have a column 
$$
Y_n = \{I(i,n) \mid i=0,\dots,h_n-1 \}
$$
 consisting of 
 of
disjoint intervals $I(i,n) \subset \Bbb R$ such that  $\bigsqcup_{i=0}^{h_n-1} I(i,n) = [0,\nu_n(F_n))$.
 Then we define a continuous  {\it $n$-th column mapping}  
 $$
 T^{(n)}:\big[0,\nu_n(F_n)\big)\setminus I(h_n-1,n)\to \big[0,\nu_n(F_n)\big)\setminus I(0,n)
 $$ 
 such that  $T^{(n)}\restriction I(i,n)$ is 
 the    orientation preserving affine mapping of $I(i,n)$ onto $I(i+1,n)$ for $i=0,\dots,h_n-2$.
 It is convenient to think of $I(i,n)$ as a {\it level} of $Y_n$.
  The levels may be of different length but they are  parallel to each other and the $i$-th level is above the $j$-th level if $i>j$.
  The $n$-th column mapping moves every level, except the highest one,   one level up.
  By the move here we mean the orientation preserving affine mapping.
  On the highest level, $T^{(n)}$ is not defined.
 
  On the $(n+1)$-th step of the construction, we first cut each $I(i, n)$ into subintervals $I(i+c, n+1)$,  $c\in C_{n+1}$, such that 
 $I(i+c, n+1)$ is from the left of $I(i+c', n+1)$ whenever $c<c'$ and
 $$
\frac{\text{the length of $I(i+c,n+1)$}}{\text{the length of $I(i,n)$}} =\kappa_{n+1}(c)\quad\text{for each $c\in C_{n+1}$}.
 $$
Hence, we obtain that  $\bigsqcup_{j\in F_n+C_{n+1}}I(j,n+1)=\big[0,\nu_{n}(F_n)\big)$.
Next, we cut the interval $\big[\nu_{n}(F_n),\nu_{n+1}(F_{n+1})\big)$ into subintervals $I(j, n+1)$, $j\in F_{n+1}\setminus(F_n+C_{n+1})$, such that
the length of $I(j, n+1)$ is $\nu_{n+1}(j)$ for each $j$.
These new levels are called {\it spacers}.
Thus, we obtain a new column $Y_{n+1}=\{I(j,n+1) \mid j \in F_{n+1}\}$ with
$\bigsqcup_{i\in F_{n+1}} I(i,n+1) = [0,\nu_{n+1}(F_{n+1}))$.
If an element  $c\in C_{n+1}$  is not maximal in $C_{n+1}$ then we denote by  $c^+$  the least element of $C_{n+1}$ that is greater than $c$.
We define a {\it spacer mapping} $s_{n+1}:C_{n+1}\to\Bbb Z_+$ by setting
$$
s_{n+1}(c):=
\cases
c^+-c-h_n &\text{if $c\ne \max C_{n+1}$},\\
h_{n+1}-c-h_n&\text{if $c= \max C_{n+1}$.}
\endcases
$$
The subcolumn $Y_{n,c}:=\{I(i+c,n+1)\mid i\in F_n\}\subset Y_{n+1}$ is called {\it the $c$-copy of $Y_n$}, $c\in C_{n+1}$.
Thus, $Y_{n+1}$  consists of $\# C_{n+1}$  copies  of  $Y_{n}$ and spacers between them and above the highest copy of $Y_n$.
More precisely, 
there are exactly $s_{n+1}(c)$ spacers above the $c$-copy of $Y_n$ in $Y_{n+1}$.
The $(n+1)$-th column mapping 
 $$
 T^{(n+1)}:\big[0,\nu_{n+1}(F_{n+1})\big)\setminus I(h_{n+1}-1,n+1)\to \big[0,\nu_{n+1}(F_{n+1})\big)\setminus I(0,n+1)
 $$ 
 is defined in a similar way as $T^{(n)}$.
 Of course, 
 $$
 T^{(n+1)}\restriction\big( \big[0,\nu_n(F_n)\big)\setminus I(h_n-1,n)\big)=T^{(n)}.
 $$
 Passing to the limit as $n\to\infty$, we obtain a well-defined nonsingular (piecewise affine) transformation
 $Q$ of the interval $[0,\lim_{n\to\infty}\nu_n(F_n))\subset\Bbb R$ equipped with Lebesgue measure
 such that 
 $$
 Q\restriction \Big( \big[0,\nu_n(F_n)\big)\setminus I(h_n-1,n)\Big)=T^{(n)}\qquad\text{for each $n\in\Bbb N$.}
 $$
 It is possible that $\lim_{n\to\infty}\nu_n(F_n)=\infty$ and then $T$ is defined on $[0,+\infty)$.
 Of course, this transformation (or, more precisely, the $\Bbb Z$-action generated by $Q$)
 is isomorphic to $(X,\mu, T)$.
 The according  nonsingular isomorphism is generated by the following correspondence:
 $$
 X\supset [i]_n \longleftrightarrow I(i,n)\subset \big[0,\lim_{n\to\infty}\nu_n(F_n)\big),\quad i\in F_n,n\in\Bbb N.
 $$
Without loss of generality we may assume $s_n(\max C_n)=0$ for each $n>0$ (see, for instance, \cite{Da4}).
This means that there are no spacers over the highest copy of $Y_n$ in $Y_{n+1}$.

 
 \subhead 1.7. Nonsingular rank-one $\Bbb Z$-actions as  transformations built under function over 
 nonsingular odometer base
 \endsubhead
 Let $(X,\mu, T)$ be as in the previous subsection.
 We will assume that $\max F_{n+1}=\max F_n+\max C_{n+1}$ for each $n\ge 0$.
 In terms of the geometrical cutting-and-stacking (see \S1.6) this
 means exactly that there are no spacers on the top of the $(n+1)$-th column.
 We remind that $X_0=C_1\times C_2\times\cdots$,  $\mu\restriction X_0=\kappa_1\otimes\kappa_2\otimes\cdots$ and 
 $$
  c^+:=\min\{d\in C_n\mid d>c\}.
 $$
 for each $n>0$ and $c\in C_n$ such that $c\ne\max C_n$.
Denote by $R$  the transformation induced by $T_1$ on $(X_0,\mu\restriction X_0)$.
Since $T_1$ is conservative, $R$ is a well defined nonsingular transformation of  $X_0$.
  Take $x=(c_1,c_2,\dots)\in X_0$.
Choose $n\ge 0$ such that $c_i=\max C_i$ for each $i=1,\dots,n$ and $c_{n+1}\ne\max C_{n+1}$.
  It is straightforward to verify that
  $$
 Rx:=
 (\underbrace{0,\dots,0}_{n\text{ times}}, c_{n+1}^+, c_{n+2}, c_{n+3}\dots)
 $$
Thus, $R$ is a classical {\it nonsingular odometer  of product type} (see \cite{Aa}, \cite{DaSi} and references therein).
Let $x^{\max}:=(\max C_1,\max C_2,\dots)\in X_0$.  
We now define a function $\vartheta:X_0\setminus\{x^{\max}\}\to\Bbb Z_+$ by setting
$$
\vartheta(x):=s_n(c_n)
$$
if $x=(c_1,c_2,\dots)$, $c_i=\max C_i$ for $i=1,\dots,n-1$ and  $c_n\ne \max C_n$, where
$s_n$ is the spacer mapping (see \S1.6).
Of course, $\vartheta$ is continuous.
Then $(X,\mu,T_1)$ is isomorphic to the transformation $R^\theta$ built under $\vartheta$ over the base $R$.
Indeed, 
let 
$$
X_0^\vartheta:=\{(x,j)\mid x\in X_0\setminus \{x^{\max}\}, 0\le j\le \vartheta(x)\}
$$
be the space of $R^\vartheta$.
Furnish  $X_0^\vartheta$ with a probability Borel measure $\mu^\vartheta$  such  that 
\roster
\item"---" $\mu^\vartheta$
projects onto $\mu$ under the natural projection $X_0^\vartheta\to X$ and 
\item"---"  for each $x\in X_0$, the corresponding
conditional measure over $x$ is non-degenerated on the entire fiber $\{x\}\times \{0,1,\dots,\vartheta(x)\}$.
\endroster
We now construct a nonsingular  isomorphism $\phi$ of $(X,\mu, T)$ onto $(X_0^\vartheta,\mu^\vartheta,R^\vartheta)$.
Given $x\in X$, we find $n>0$ such that $x\in X_n\setminus X_{n-1}$.
Then $x=(f_n,c_{n+1},c_{n+2},\dots)$ with  $f_n\in F_n\setminus F_{n-1}C_n$
and $c_j\in C_j$ for each $j>n$.
Hence $f_n=j+c_n+h_{n-1}$ for some (uniquely determined) $c_n\in C_n$ and $0< j\le s_n(c_n)$.
Since there are no spacers over  $\max C_n$ in the $n$-th column, we obtain that $c_n\ne\max C_n$.
We now let
$$
\phi(x):=\big((\max C_1,\dots,\max C_{n-1}, c_n,c_{n+1},\dots), j\big)\in X_0^\vartheta.
$$
 It is straightforward to verify that $\mu\circ\phi^{-1}\sim\mu^\vartheta$ and $\phi T=R^\vartheta\phi$, as desired.

\head{2. Finite factors of nonsingular $(C,F)$-actions}\endhead

 Let a sequence $\Cal T=(C_n,F_{n-1},\kappa_n,\nu_{n-1})_{n=1}^\infty$ satisfy  \thetag{1-1}--\thetag{1-3} and Proposition~1.10(ii) and let $\lim_{n\to\infty}\nu_n(F_n)<\infty$.
  Let $\Gamma$ be a cofinite  subgroup of $G$. 
We consider the left coset space $G/\Gamma$ as a homogeneous $G$-space on which $G$ acts by left translations.
It is obvious  that  for each coset $g\Gamma\in G/\Gamma$,  the subgroup $g\Gamma g^{-1}\subset G$ is the stabilizer of  $g\Gamma$ in $G$.

\definition{Definition 2.1} 
Given  a coset $g\Gamma\in G/\Gamma$, we say that
${\Cal T}$  is {\it compatible with $g\Gamma$} if
$$
 \sum_{n=1}^\infty\kappa_n\big(\{c\in C_n\mid c\not\in g\Gamma g^{-1}\}\big)<\infty.
 $$ 
 \enddefinition

  Denote by $(X,\mu,T)$ the $(C,F)$-action of $G$  associated with $\Cal T$.
  Then $\mu(X)=\lim_{n\to\infty}\nu_n(F_n)<\infty$.
For a point $ x=(c_1,c_2,\dots)\in X_0= C_1\times   C_2\times\cdots$,  we let
$$
\pi_{({\Cal T},g\Gamma)}( x):=\lim_{n\to\infty}  c_1 c_2\cdots  c_ng\Gamma\in G/\Gamma
$$
whenever this limit exists\footnote{The quotient space $G/\Gamma$ is endowed with the discrete topology.}.
It follows from the Borel-Cantelli lemma that if $\Cal T$  is compatible with  $g\Gamma\in G/\Gamma$ then $\pi_{({\Cal T},g\Gamma)}(x)$ is well defined for $\mu$-a.e.
$x\in X_0$.
It is straightforward to verify that for each $h\in G$,
$$
\pi_{({\Cal T},g\Gamma)}( T_h x)=h\pi_{({\Cal T},g\Gamma)}( x)
$$ 
whenever $\pi_{({\Cal T},g\Gamma)}( x)$ is well defined and $ T_h x\in X_0$.
It follows that  the mapping 
$$
\pi_{({\Cal T},g\Gamma)}: X_0\ni
 x\mapsto \pi_{({\Cal T},g\Gamma)}( x)\in G/\Gamma
$$ 
extends uniquely (mod 0) to a  measurable $G$-equivariant mapping from $ X$ to $G/\Gamma$.  
We denote the extension by the same symbol $\pi_{(\Cal T,g\Gamma)}$.
It is routine to verify that if $x=(  f_n, c_{n+1}, c_{n+2},\dots)\in  X_n$ for some $n\in\Bbb N$ then
$$
\pi_{({\Cal T},g\Gamma)}( x):=\lim_{m\to\infty} f_n c_{n+1}c_{n+2}\cdots c_mg\Gamma.\tag2-1
$$
\comment

It follows that  if $\widetilde{\Cal T}$  is compatible with  $g\Gamma$ then the mapping 
$$
\pi_{(\widetilde{\Cal T},g\Gamma)}:=\widetilde\pi_{(\Cal T,g\Gamma)}\circ\iota\tag2-2
$$
is a well defined $G$-equivariant measurable mapping from $X$ to $G/\Gamma$,
where $\iota$ stands for  the canonical isomorphism of $(X,\mu)$ onto $(\widetilde X,\widetilde \mu)$ (see \thetag{1-7}).

\endcomment

\definition{Definition 2.2} We call $\pi_{({\Cal T},g\Gamma)}$ {\it the $({\Cal T},g\Gamma)$-factor mapping}
for $T$.
\enddefinition

\comment

 Let $\widetilde{\Cal T}
=(\widetilde C_n,\widetilde F_{n-1},\widetilde\kappa_n,\widetilde\nu_{n-1})_{n=1}^\infty$
be a telescoping of $\Cal T$.

Denote by $\widetilde T$ the $(C,F)$-action of $G$  associated with $(C_n,F_{n-1})_{n=1}^\infty$.
 Let $\widetilde X$ be the space of $\widetilde T$ and let $\widetilde \mu$ stand for the 
 the nonsingular $(C,F)$-measure on $\widetilde X$ determined by $(\widetilde\kappa_n)_{n=1}^\infty$ and $(\widetilde\nu_n)_{n=0}^\infty$.

 We need some notation.
 Given two positive integers $n\le m$, we denote the subset $C_n\cdots C_m$ of $G$ by $C_{n,m}$.
It is straightforward to verify that $[f]_n=\bigsqcup_{c\in C_{n,m}}[fc]_m$ for each $f\in F_n$.
 We denote by by $\kappa_{n,m}$ the probability measure $\kappa_n*\cdots*\kappa_m$ which is supported on the finite subset $C_n\cdots C_m$ of $ G$.
\endcomment

We need some notation. 
 Given $1<n<m$, we denote the subset $C_n\cdots C_m$ of $G$ by $C_{n,m}$.
 Let  $\kappa_{n,m}$ stand for the probability measure $\kappa_n*\cdots*\kappa_m$.
 It is   supported on  $C_n\cdots C_m$.
 We now  state
 the main result of the section.

\proclaim{Theorem 2.3}
	The following are equivalent.
	\roster
	\item"(i)" There is a measurable factor map $\tau: X\to G/\Gamma$, i.e. for each $g\in G$,
	$$
	 \tau(T_gx)=g\tau(x)\quad\text{at $\mu$-a.e. $x\in X$}
	$$
	\item"(ii)" There exists a sequence $(g_n)_{n>0}$ of elements of  $G$ such that
	$$
	\lim\limits_{N\to\infty}\sup\limits_{m> n\geq N}\kappa_{n+1,m}\big(\{c\in C_{n+1,m} \mid c\notin g_n^{-1}\Gamma g_m\}\big) = 0. 
	$$ 
	\item"(iii)" 
	There exist   a coset  $g_0\Gamma\in G/\Gamma$ and a $g_0\Gamma$-compatible 
	telescoping of $\Cal T$. 	
	\endroster
	It follows that $T$ has no factors isomorphic to $G/\Gamma$ if and only if there is  no telescoping of $\Cal T$  compatible with $g\Gamma g^{-1}$ for any $g\in G$.
\endproclaim

\demo{Proof} (i)\,$\Longrightarrow$\,(ii)
Let $Y_j := \tau^{-1}(j)$, for each $j\in G/\Gamma$. 
Then $X=\bigsqcup_{j\in G/\Gamma}Y_j$. 
Consider $n$ large so that $\mu(Y_j\cap X_n)>0$ for each $j\in G/\Gamma$.
Let $g,h\in F_n$.
Since $\tau$ is $G$-equivariant, it follows that   $T_{hg^{-1}} Y_j = Y_{hg^{-1}j}$ and
hence
$$
	T_{hg^{-1}}([g]_n\cap Y_j) = [h]_n\cap Y_{hg^{-1}j}.
	\tag{2-2}
$$
For each  $j\in G/\Gamma$ and $g \in F_n$, let 
$$
d_{n, g}(j) := \mu([g]_n\cap Y_j)/\mu([g]_n).
$$
 Then, the set $\{d_{n, g}(j)\mid j\in G/\Gamma\}\subset (0,1)$ does not depend on $g\in F_n$.
Indeed, for each $h\in F_n$, 
the Radon-Nikodym derivative of the transformation $T_{gh^{-1}}$ is constant on $[h]_n$ and we deduce from \thetag{2-2} that
$$
d_{n,h}(j) := \frac{\mu([h]_n\cap Y_j)}{\mu([h]_n)} =
\frac{\mu([g]_n\cap Y_{gh^{-1}j})}{\mu([g]_n)}=d_{n,g}(gh^{-1}j).
\tag2-3
$$
Hence
$
\{d_{n, h}(j): j\in G/\Gamma\} = 
\{d_{n, g}(j): j\in G/\Gamma\}.
$
Let 
$$
\delta_n := \max\limits_{j\in G/\Gamma} d_{n,g}(j). 
$$
We claim that $\delta_m \to 1$ as $m\to\infty$. 
Indeed, for $m\ge n$, let ${\Cal P}_m$ denote the  finite $\sigma$-algebra  generated 
by  the $m$-cylinders (which are compact and open subsets of $X$) that are contained 
in $X_n$. 
Then  ${\Cal P}_n\subset{\Cal P}_{n+1}\subset\cdots$ and $\bigvee_{m>n}{\Cal P}_m$ is the entire Borel $\sigma$-algebra on $X_n$.
Hence,  for each $j\in G/\Gamma$, there is $g_m\in F_m$ such that  
$$
\frac{\mu([g_m]_m\cap Y_j)}
{\mu([g_m]_m)}\to 1\quad\text{ as $m\to\infty$.}
$$
This implies that $\delta_m \to 1$ as $m\to\infty$, as claimed.
In what follows we consider $n$ large so that
$\delta_n > 0.9$.
Then for each $g\in  F_n$, there is
a unique $\Gamma$-coset $j_n(g)\in G/\Gamma$ such that $\delta_n=d_{n,g}(j_n(g))$.
It is convenient to consider $G/\Gamma$ as a set of {\it colors}.
Then  
 $j_n(g)$ is
 {\it the dominating color} on $[g]_n$. 
 It follows from \thetag{2-3} that
  $j_n(g)=gh^{-1} j_n(h)$, for all $g,h\in F_n$.
Choose $g_n\in F_n$ such that $j_n(g_n) = \Gamma$. 
Then $j_n(g)= gg_n^{-1} \Gamma$
for each $g\in F_n$.
Given $\epsilon < \frac{1}{2}$, there is $N > 0$ such that $\delta_n > 1-\epsilon^2$ for all $n > N$. Hence for all 
$m > n > N$,
$$
(1-\epsilon^2)\mu([1_G]_n)<\mu([1_G]_n\cap Y_{j_n(1_G)}).\tag2-4
$$
We remind that  $[1_G]_n=\bigsqcup_{c\in C_{n+1,m}}[c]_m$.
Let 
$$
D:=\{c\in C_{n,m}\mid \mu([c]_m\cap Y_{j_n(1_G)})>(1-\epsilon)\mu([c]_m)\}.
$$
It follows from \thetag{2-4} that
$$
\align
(1-\epsilon^2)\mu([1_G]_n)&<\sum_{c\in D}\mu([c]_m)+(1-\epsilon)\sum_{c\in C_{n+1,m}\setminus D}\mu([c]_m)\\
&=\sum_{c\in D}\mu([c]_m)+(1-\epsilon)\Bigg(\mu([1_G]_n)-\sum_{c\in D}\mu([c]_m)\Bigg).
\endalign
$$
This yields that $\sum_{c\in D}\mu([c]_m)>(1-\epsilon)\mu([1_G]_n)$ or, equivalently,
$$
\sum_{c\in D}\kappa_{n+1,m}(c)>(1-\epsilon)\kappa_n(1_G).\tag2-5
$$
By the definition of $D$, an element $c\in C_{n+1,m}$ belongs to  $D$ if and only if $j_n(1_G)$ is the dominating color on $[c]_m$, i.e. $j_m(c)=j_n(1_G)$.
Therefore, $cg_m^{-1}\Gamma=g_n^{-1}\Gamma$, i.e. $c\in g_n^{-1}\Gamma g_m$.
Hence \thetag{2-5} yields that
$$
\kappa_{n+1,m}\big(\{c\in C_{n+1,m}\mid c\not\in g_n^{-1}\Gamma g_m\}\big)\le \epsilon
$$
and (ii) follows.

 (ii)\,$\Longrightarrow$\,(iii)
As $\Gamma$ is cofinite and (ii) holds, there exist  an  increasing sequence  $0=q_0<q_1 <q_2 < \cdots$ of positive integers  and an element $g_0\in G$ such that  $g_n^{-1}\in g_0\Gamma $  and
$$
\kappa_{q_n+1,q_{n+1}}\big(
	 \{c\in C_{q_{n}+1, q_{n+1}} \mid c\notin g_0\Gamma g_0^{-1}\}\big)<2^{-n}
	 $$
for each $n\in\Bbb N$.
Let $\boldsymbol q:=(q_n)_{n=0}^\infty$.
Then $g_0\Gamma$ is compatible with the $\boldsymbol q$-telescoping of $\Cal T$.
This implies (iii).

 (iii)\,$\Longrightarrow$\,(i)  Denote by $\widetilde{\Cal T}$ the $\boldsymbol q$-telescoping of $\Cal T$.
 Then  $\pi_{(\widetilde{\Cal T},g_0\Gamma)}\circ \iota_{\boldsymbol q}$
 is a factor mapping of $T$ onto $G/\Gamma$.
 
 Thus, the first statement of the theorem is proved completely.
 The second statement follows from the first one and a simple observation that given two subgroups $\Gamma$ and $\Gamma'$ of $G$, the corresponding $G$-actions by left translations
 on $G/\Gamma$ and $G/\Gamma'$ are isomorphic if and only if $\Gamma$ and $\Gamma'$ are conjugate.
 \qed

 \comment
 For each $n\in\Bbb N$, the mapping
$$
(c_{g_n+1},\dots,c_{q_{n+1}})\mapsto c_{g_n+1}\cdots c_{q_{n+1}}
$$
 is a bijection of $C_{q_n+1}\times\cdots\times C_{q_{n+1}}$ onto $C_{q_n+1,q_{n+1}}$.
We identify the space
$$
X_{q_n}=F_{q_n}\times(C_{q_n+1}\times\cdots\times C_{q_{n+1}})\times
(C_{q_{n+1}+1}\times\cdots\times C_{q_{n+2}})\times\cdots
$$
 with the the infinite product space 
$
\widetilde{X}_{q_n}:=F_{q_n}\times C_{q_n+1,q_{n+1}}\times  C_{q_{n+1}+1,q_{n+2}}\times\cdots
$
via these mappings.
Under this identification,  the measure $\mu\restriction X_{q_n}$ corresponds to
the infinite  product measure $\widetilde\mu_{n}:=\nu_{q_n}\times\kappa_{q_n+1,q_{n+1}}\times\kappa_{q_{n+1}+1,q_{n+2}}\times\cdots$ on
$\widetilde{X}_{q_n}$.
Let 
$$
D_n:=\{d\in C_{q_{n+1}+1,q_{n+2}}\mid  d\in g_0^{-1}\Gamma g_0\}.
$$
It follows from the condition of (iii) and the Borel-Cantelli lemma that for $\widetilde\mu_n$-a.e. 
$
\widetilde x=(f_n,d_{n+1},d_{n+2},\dots)\in \widetilde{X}_{q_n},
$
 there is $N_{\widetilde x}>n$ such that
$d_m\in D_m$ for each $m>N_{\widetilde x}$. 
We now let  for each $x=(f_{q_n},c_{q_n+1},\dots)\in X_{q_n}$,
$$
\pi_n(x):=\lim_{m\to\infty}f_{q_n}c_{q_n+1}\cdots c_{q_m} g_0^{-1}\Gamma\in G/\Gamma. \tag2-5
$$
Since $c_{q_m+1}\cdots c_{q_{m+1}}g_0^{-1}\Gamma =g_0^{-1}\Gamma$ eventually in $m$ for $\mu$-a.e. $x\in X_{q_n}$,
the  limit ~\thetag{2-5} exists  for $\mu$-a.e. $x\in X_{q_n}$.
Of course, $\pi_n$ is a well defined Borel mapping from $X_n$ (or, more precisely, from a $\mu$-conull subset of $X_n$) to $G/\Gamma$.
It is straightforward to verify that  $\pi_n(T_gx)=g\pi_n(x)$ if $gf_nd_{n+1}\cdots d_m\in F_m$ for some $m>n$.
It is straightforward to verify that $\pi_{n+1}\restriction X_n=\pi_n$ for each $n\in\Bbb N$.
Hence a Borel mapping $\pi:X\to G/\Gamma$ is well defined ($\mu$-mod 0) by the formula
$\pi\restriction X_n:=\pi_n$ for each $n\in\Bbb N$.
Of course, $\pi$ is $G$-equivariant.

\endcomment
\enddemo

The following  important remark will be used essentially in the proof of the main result of \S3.

\remark{Remark 2.4}   In fact, we  obtained  more than what is stated in Theorem~2.3.
We proved indeed that given a factor mapping $\tau:X\to G/\Gamma$, there exist a coset $g_0\Gamma$ and a 
$g_0\Gamma$-compatible $\boldsymbol q$-telescoping $\widetilde{\Cal T}$ of $\Cal T$ such that $\tau=\pi_{(\widetilde{\Cal T},g_0\Gamma)}\circ \iota_{\boldsymbol q}$.
 To explain this fact, we use below the notation from the proof of Theorem~2.3.
 For $n\in\Bbb N$, let 
 $$
 X_n':=\bigsqcup_{f_n\in F_n}\big([f_n]_n\cap Y_{j_n(f_n)}\big)\subset X_n.
 $$
  We remind that $\mu(X)<\infty$.
 Since  $\delta_n\to 1$ and $\mu(X_n)\to\mu(X)$ as $n\to\infty$, it follows that  
 $\mu(X_{q_n}') \to \mu(X)$ as $n\to\infty$.
 We can assume (passing to a subsequence of $(q_n)_{n=1}^\infty$ if needed) that $\sum\limits_{n=1}^{\infty} \mu(X\setminus X_{q_n}') < \infty$.
 Then, the Borel-Cantelli lemma yields that for a.e. $x\in X$, we have that  $x\in X_{q_n}'$ eventually in $n$.
 Hence, for a.e. $x\in X$,
 $$
\tau(x) = \lim\limits_{m\to\infty} j_{q_m}(f_{q_m}) = \lim\limits_{m\to\infty} f_{q_m}g_{q_m}^{-1}\Gamma = \lim\limits_{m\to\infty} f_{q_n} g_0 \Gamma,\tag2-6
$$
 where $f_{q_m}$ is the first coordinate of $x$ in $X_{q_m}$, i.e. $x=(f_{q_m}, c_{q_m+1},\dots)\in X_{q_m}$.
 On the other hand, it follows from \thetag{2-1} and \thetag{2-2} that
 $$
 \pi_{\widetilde{(\Cal T}, g_0\Gamma)} (\iota_{\boldsymbol q}(x))=\lim_{m\to\infty} f_{q_m}g_0\Gamma\tag2-7
 $$
 at a.e. $x\in X$.
 Therefore, \thetag{2-6} and \thetag{2-7} yield that $ \tau= \pi_{\widetilde{(\Cal T}, g_0\Gamma)}\circ\iota_{\boldsymbol q}$ almost everywhere, as desired.
\endremark

\vskip 10pt

We note in this connection that if 
\roster
\item"---" $G$ is Abelian or  $G$ is arbitrary but 
 $\Gamma$ is normal in $G$ and
 \item"---" the homogeneous space $G/\Gamma$ is a factor of an ergodic  nonsingular free action 
 of $G$  on a standard measure space $(Y,\goth Y,\nu)$
 \endroster
  then this factor (considered as an invariant  sub-$\sigma$-algebra of $\goth Y$) is defined uniquely by $\Gamma$.
Indeed, if $\pi_1,\pi_2:Y\to G/\Gamma$ are two $G$-equivariant measurable maps then
the mapping $Y\ni y\mapsto \pi_1(y)\pi_2(y)^{-1}\in G/\Gamma$ is invariant under $G$.
Hence it is constant.
Therefore there is $a\in G/\Gamma$ such that $\pi_1(y)=a\pi_2(y)$ for a.e. $y\in Y$.
It follows that 
$$
\{\pi_1^{-1}(j)\in \goth Y\mid j\in G/\Gamma\}=\{\pi_2^{-1}(j)\in \goth Y\mid j\in G/\Gamma\}.
$$
Therefore the equality  $ \tau= \pi_{\widetilde{(\Cal T}, g_0\Gamma)}\circ \iota_{\boldsymbol q}$  (at least, up to a rotation of $G/\Gamma$) stated in Remark~2.4 is a 
trivial fact.
However, it  is no longer true if $\Gamma$ is not normal.

\example{Example 2.5} Let  $\Bbb Z_3:=\Bbb Z/3\Bbb Z$, $\Bbb Z_2:=\Bbb Z/2\Bbb Z$,
 $G=\Bbb Z_3\rtimes\Bbb Z_2$
and $\Gamma= \{0\}\times \Bbb Z_2$.
Then $\Gamma$ is a non-normal cofinite subgroup of $G$ of index $3$.
We consider $\Bbb Z_3$ as a quotient   $G/\Gamma$.
Then $\Bbb Z_3$ is a $G$-space.
Hence the product space $\Bbb Z_3\times\Bbb Z_3$ is also a 
$G$-space (we consider the diagonal $G$-action).
Since the diagonal $D=\{(j,j)\mid j\in G/\Gamma\}$  is an invariant subspace of  $\Bbb Z_3\times\Bbb Z_3$, the complement 
$Y:=(\Bbb Z_3\times\Bbb Z_3)\setminus D$ of $D$ in $\Bbb Z_3\times\Bbb Z_3$ is also $G$-invariant.
It is easy to verify the $G$-action on  $Y$ is  transitive and free.
Endow $Y$ with the (unique) $G$-invariant probability measure $\nu$.
Of course,  the coordinate projections $\pi_1,\pi_2:Y\to\Bbb Z_3$ are two-to-one $G$-equivariant maps. 
However, the corresponding $\sigma$-algebras  of $\pi_1$-measurable and $\pi_2$-measurable subsets in $Y$ are different.
Consider a  rank-one $\Bbb Z$-action on a standard probability space $(Z,\goth Z,\kappa)$.
Then the product $(\Bbb Z\times G)$-action  on $(Z\times Y,\kappa\otimes\nu)$ is of  rank one.
Denote it by $R$.
The subgroup $\Bbb Z\times \Gamma$ of $\Bbb Z\times G$ is non-normal.
It is  of index 3.
Hence we can consider the corresponding finite quotient space 
$\Bbb Z_3$ as a  $(\Bbb Z\times G)$-space.
The  mappings $1\otimes\pi_1$ and $1\otimes\pi_2$ from $Z\times Y$ onto $\Bbb Z_3$ are $(\Bbb Z\times G)$-equivariant.
However the corresponding factors of $R$, i.e. the invariant sub-$\sigma$-algebras,  are different. 
\endexample

A nonsingular $G$-action is totally ergodic if and only if has no non-trivial finite factors or, equivalently, each cofinite subgroup of $G$ acts ergodically.
We thus deduce from Theorem~2.3 the following criterion of total ergodicity for the rank-one nonsingular actions.

\proclaim{Corollary 2.6} Let $T$ be a $(C,F)$-action of $G$ associated with a sequence $\Cal T$ satisfying \thetag{1-1}--\thetag{1-3} and  Proposition~1.10(ii).
Then $T$ is totally ergodic if no telescoping of $\Cal T$ is compatible with any proper
 cofinite subgroup of $ G$, i.e. for each increasing sequence $n_1<n_2<\cdots$ of integers and
 each proper cofinite subgroup $\Gamma$ in $G$,
 $$
 \sum_{k=1}^\infty \kappa_{n_k+1}*\cdots *\kappa_{n_{k+1}}\big(\{c\in C_{n_k+1}\cdots C_{n_{k+1}}\mid c\not\in\Gamma\}\big)=\infty.
 $$
\endproclaim

\head{3. Nonsingular odometer  actions of residually finite groups}\endhead
 \subhead 3.1. Nonsingular odometers\endsubhead
 From now on $G$ is residually finite.
Then there is  a decreasing sequence  $\Gamma_1\supsetneq \Gamma_2\supsetneq\cdots$ of cofinite subgroups $\Gamma_n$ in $G$ such that
 $$
 \bigcap_{n=1}^\infty\bigcap_{g\in G}g\Gamma_ng^{-1}=\{1_G\}.\tag3-1
 $$
Consider the natural  inverse  sequence of homogeneous $G$-spaces and $G$-equivariant mappings intertwining them:
$$
G/\Gamma_1\longleftarrow G/\Gamma_2\longleftarrow\cdots.\tag3-2
$$
Denote by $Y$ the projective limit of this sequence.
A point of $Y$ is a sequence $(g_n\Gamma_n)_{n=1}^\infty$ such that
$g_n\Gamma_n=g_{n+1}\Gamma_n$, i.e. $g_n^{-1}g_{n+1}\in \Gamma_n$ for each $n>0$.
Endow $Y$ with the topology of projective limit.
Then $Y$ is a  compact Cantor $G$-space. 
Of course, the $G$-action on $Y$ is  minimal and uniquely ergodic.
Denote this action by $O=(O_g)_{g\in G}$.
It follows from \thetag{3-1} that $O$ is {\it faithful}, i.e.  $O_g\ne I$ if $g\ne 1_G$.
We note that a faithful action is not necessarily free.

\definition{Definition 3.1}
The  dynamical system $(Y,O)$ is called {\it the topological $G$-odometer associated with $(\Gamma_n)_{n=1}^\infty$}.
If $\nu$ is a non-atomic Borel measure on $Y$ which is quasi-invariant and ergodic under $O$ then
we call the dynamical system $(Y,\nu,O)$ a {\it nonsingular $G$-odometer}.
By the {\it Haar measure} for $(Y,O)$ we mean the unique $G$-invariant probability on $Y$.
\enddefinition

We note that \thetag{3-1} is in no way restrictive. 
Indeed, let  $\bigcap_{n=1}^\infty\bigcap_{g\in G}g\Gamma_ng^{-1}=N\ne\{1_G\}$.
Define $(Y,O)$ as above.
Then $N$ is a proper normal subgroup of $G$ and $N=\{g\in G\mid O_g=I\}$.
We now let $\widetilde G:=G/N$ and $\widetilde\Gamma_n:=\Gamma_n/N$.
Then $\widetilde\Gamma_n$ is a cofinite subgroup in $\widetilde G$ for each $n\in\Bbb N$,
$\widetilde \Gamma_1\supsetneq \widetilde \Gamma_2\supsetneq\cdots$ and
 $\bigcap_{n=1}^\infty\bigcap_{g\in \widetilde G}g\widetilde\Gamma_ng^{-1}=\{1_{\widetilde G}\}$.
 Let $(\widetilde Y,\widetilde O)$ denote the topological $\widetilde G$-odometer associated with
 the sequence $(\widetilde\Gamma_n)_{n=1}^\infty$.
 Then, of course, $Y=\widetilde Y$ and $O_g=\widetilde O_{gN}$ for each $g\in G$.

We now isolate a class of nonsingular odometers  of rank one.
For each $n>0$, we choose a finite subset $D_n\subset\Gamma_{n-1}$ such that $1_G\in D_n$ and each $\Gamma_{n}$-coset in $\Gamma_{n-1}$ intersects $D_n$ exactly once\footnote{
For consistency of the notation, we let $\Gamma_0:=G$.}.
We then call $D_n$ {\it a $\Gamma_n$-cross-section} in $\Gamma_{n-1}$.
Then the product $D_1\cdots D_n$ is a $\Gamma_n$-cross-section in $G$. 
Hence there is a unique bijection 
$$
\omega_n:G/\Gamma_n\to D_1\cdots D_n
$$ 
such that
$\omega_n(g\Gamma_n)\Gamma_n=g\Gamma_n$ for each $g\in G$ and $\omega_n(\Gamma_n)=1_G$.
It follows, in particular, that 
$$
\text{if $\omega_n(g\Gamma_n)=h\omega_n(g'\Gamma_n)$ for some
$g,g',h\in G$ then $g\Gamma_n=hg'\Gamma_n$.}\tag3-3
$$
It is straightforward to verify that the  diagram 
$$
\CD
G/\Gamma_1 @<<<G/\Gamma_2 @<<< G/\Gamma_3 @<<<\cdots\\
@V{\omega_1}VV @V{\omega_2}VV @V{\omega_3}VV \\
D_1@<\phi_1<< D_1D_2@<\phi_2<< D_1D_2D_3@<\phi_3<<\cdots\\
@V\psi_1VV @V\psi_2VV @V\psi_3VV \\
D_1@<\phi_1'<< D_1\times D_2@<\phi_2'<< D_1\times D_2\times D_3@<\phi_3'<<\cdots\\
\endCD
\tag3-4
$$
commutes.
The horizontal arrows in the upper line denote the natural projections.
The other mappings in the diagram are defined as follows: 
$$
\align
&\phi_n(d_1\cdots d_{n+1}):=d_1\cdots d_{n}, \\
&\phi_n'(d_1,\dots, d_{n+1}):=(d_1,\dots, d_{n})\quad\text{and}\\
&\psi_n(d_1\cdots d_{n}):=(d_1,\dots, d_{n})
\endalign
$$
for each $(d_1,\dots,d_{n+1})\in D_1\times\cdots\times D_{n+1}$ and $n\ge 0$.
It follows from \thetag{3-3} that there exists  a natural  homeomorphism\footnote{The homeomorphism  pushes down to a bijection between $G/\Gamma_n$ and $D_1\times\cdots\times D_n$ for each $n$.} of $Y$ onto the infinite product
space $D:=D_1\times D_2\times\cdots$.

\proclaim{Proposition 3.2} If, for each $n\in\Bbb N$, there is a $\Gamma_n$-cross-section $D_n$ in $\Gamma_{n-1}$ and a probability $\kappa_n$ on $G$ such that
\roster
\item"(i)" $\text{{\rm supp}\,} \kappa_n=D_n$ for each $n$, 
\item"(ii)" $\prod_{n=1}^\infty\max_{d\in D_n}\kappa_n(d)=0$ and 
\item"(iii)"
$\lim_{n\to\infty}(\kappa_1*\cdots*\kappa_n)(gD_1\cdots D_n)=1$ for each $g\in G$
\endroster
then  
there is a  non-atomic probability Borel measure $\mu$ on $Y$ which is quasi-invariant under $O$ and such that the nonsingular odometer $(Y,\mu,O)$ is of rank one along the sequence $(D_1\cdots D_n)_{n=1}^\infty$.
\endproclaim

\demo{Proof}
We  set $F_0:=\{1_G\}$, $F_{n}:=D_1\cdots D_{n}$, $C_n:=D_n$  and $\nu_n:=\kappa_1*\cdots*\kappa_n$ for each $n\in\Bbb N$.
Then \thetag{1-1}--\thetag{1-3} hold for the sequence $(C_n,F_{n-1},\kappa_n,\nu_{n-1})_{n=1}^\infty$.
Moreover,~Proposition~1.10(iii) is exactly (iii) in the case under consideration.
Hence the $(C,F)$-action $T$ of $G$ associated with 
$(C_n,F_{n-1},\kappa_n,\nu_{n-1})_{n=1}^\infty$ is well defined.
Let
$(X,\mu)$ stand for the space of this action. 
It follows from the $(C,F)$-construction that 
$$
X=D_1\times D_2\times\cdots=D,
$$
 $\mu$ is non-atomic and $T$ is free (mod $\mu$).
 In view of \thetag{3-3}, we can  identify $X$ with $Y$.
 Hence we consider $\mu$ as a probability on $Y$.
 Moreover, \thetag{3-2} yields that $T$ is conjugate to $O$.
 Thus, $(Y,\mu,O)$ is a nonsingular odometer.
It remains to apply~Theorem~1.13.
\qed
\enddemo

\comment

We leave the proof of the following simple statement to the reader.

\proclaim{Proposition 3.3}
If there is a F{\o}lner sequence $(F_n)_{n=1}^\infty$ in $G$ such that $F_n$ is a $\Gamma_n$-cross-section in $G$ then the odometer $(Y,O)$ equipped with the Haar measure is of rank one.
\endproclaim

\endcomment

We now show that the classical nonsingular $\Bbb Z$-odometers of product type are covered by Definition~3.1.

\example{Example 3.3} Let $G=\Bbb Z$ and let $(a_n)_{n=1}^\infty$ be a sequence of   integers such that $a_n>1$ for each $n\in\Bbb N$.
We set $\Gamma_n:=a_1\cdots a_n\Bbb Z$.
Then $\Gamma_1\supsetneq\Gamma_2\supsetneq\cdots$ and
$\bigcap_{n=1}^\infty\Gamma_n=\{0\}$.
The set $D_n:=a_1\cdots a_{n-1}\cdot\{0,1,\dots,a_n-1\}$ is a $\Gamma_n$-cross-section in $\Gamma_{n-1}$.
Hence, in view of  \thetag{3-4}, the space $Y$ of the $\Bbb Z$-odometer $O=(O_n)_{n\in\Bbb Z}$ associated with $(\Gamma_n)_{n=1}^\infty$
is homeomorphic to the infinite product $D=D_1\times D_2\times\cdots$.
We identify $D_n$ naturally with the set $\{0,1,\dots,a_n-1\}$.
Then
 $$
 D=
\{0,1,\dots,a_1-1\}\times\{0,1,\dots,a_2-1\}\times\cdots.
$$
To define $O$ explicitly on this space we take 
 $y=(y_n)_{n=1}^\infty\in D$.
 It is a routine to check that if  there is $k> 0$ such that $y_j=a_{j}-1$ for each $j<k$ and $y_k\ne a_k-1$ then
 $$
 O_1y=(0,\dots,0,y_k+1,y_{k+1}, y_{k+2},\dots).
 $$
 If such a $k$ does not exist, i.e. $y_j=a_j-1$ for each $j>0$, then $O_1y=(0,0,\dots)$.
 Let $\kappa_n$ be a non-degenerated probability measure on $\{0,1,\dots,a_n-1\}$ and let
 $\prod_{n>0}\max_{0\le d<a_n}\kappa_n(d)=0$.
 This means that (i) and (ii) of Proposition~3.2 hold.
 Of course, Proposition~3.2(iii) holds also.
 Hence, by Proposition~3.2, the nonsingular odometer
 $$
 \Bigg(\bigotimes_{n=1}^\infty \{0,\dots,a_n-1\}, \bigotimes_{n=1}^\infty\kappa_n, O\Bigg)
 $$
 is  of rank one.
Thus, in this case our definition of nonsingular odometer coincides with the classical definition of nonsingular $\Bbb Z$-odometers of product type (see \cite{Aa}, \cite{DaSi}).
Moreover, the $\Bbb Z$-odometers of product type are of rank one.
\endexample

It is routine to verify that if $G=\Bbb Z^d$ with $d\in\Bbb N$ then each probability preserving $G$-odometer is of rank one.
This follows from Proposition~3.2 if
 one  chooses the $\Gamma_n$-cross-sections $D_n$ in $\Gamma_{n-1}$ in such a way that the sum
$D_1+\cdots+D_n$ is a parallelepiped $\{0,1,\dots,a_{1,n}\}\times\cdots\times\{0,1,\dots,a_{d,n}\}$ for some $a_{1,n},\dots,a_{d,n}\in\Bbb N$ with $\lim_{n\to\infty}a_{j,n}=\infty$ for each $j$.
We leave details to the reader (see also \cite{JoMc, Theorem~2.11}).

Note,  however,  that there exist probability preserving free $G$-odometers which are not of rank one.

\example{Example 3.4}
The free group with 2 generators $\boldsymbol F_2$ is residually finite.
Hence there is a sequence $N_1\supsetneq N_2\supsetneq\cdots$ of normal subgroups in $\boldsymbol F_2$ such that $\bigcap_{n=1}^\infty N_n=\{1_{\boldsymbol F_2}\}$.
Then the topological $\boldsymbol F_2$-odometer associated with $(N_n)_{n=1}^\infty$
is a free minimal $\boldsymbol F_2$-action by translation on a compact group $Y$.
Let $\chi$ denote the Haar measure on $Y$.
Then $(K,\chi, O)$ is an ergodic probability preserving $\boldsymbol F_2$-odometer.
If $(K,\chi, O)$ were of rank one then $\boldsymbol F_2$ would be amenable 
by Corollary~1.11(ii), a contradiction.
This argument works also for each non-amenable residually finite group in place of  $\boldsymbol F_2$.
\endexample

We also provide two examples of non-rank-one odometer actions for amenable groups $G$.
In the first example $G$ is  locally finite, and in the second one $G$ is non-locally finite periodic.

\example{Example 3.5}
Let $Z=\{0,1\}^\Bbb N$.
Endow $Z$ with the infinite product $\eta$ of the equidistributions on $\{0,1\}$.
Fix a sequence $\boldsymbol{s}=(s_n)_{n=0}^\infty$ of mappings 
$$
s_n:\{0,1\}^n\to \text{Homeo}(\{0,1\}),\qquad n\ge 0.
$$ 
Consider the following  transformation $T_{\boldsymbol{s}}$ of $Z$:
$$
T_{\boldsymbol{s}}(z_1,z_2,\dots):=(s_0z_1, s_1(z_1)z_2,s_2(z_1,z_2)z_3,s_3(z_1,z_2,z_3)z_4,\dots).
$$
Of course, $T_{\boldsymbol{s}}$ preserves $\eta$.
Let  
$$
G:=\{T_{\boldsymbol{s}}\mid \boldsymbol{s}=(s_n)_{n=0}^\infty \text{ with $s_n\equiv I$ eventually}\}.
$$
Then $G$ is a locally finite (and hence amenable) countable group.
\comment

For $n\in\Bbb N$, let
$$
A_n:=\{(s_k)_{k=0}^\infty\in \Cal A\mid s_k\equiv I\text{ for  each }k=0,\dots,n\}.
$$
Then $A_n$ is a cofinite normal subgroup  of $\Cal A$.
Of course, $A_1\supsetneq A_2\supsetneq\cdots$ and $\bigcap_{n}A_n=\{I\}$.

\endcomment
We claim that the dynamical system $(Z,\eta, G)$ is a $G$-odometer.
Indeed, for each $n\in\Bbb N$, we denote by $\pi_n:Z\to\{0,1\}^{n}$ the projection to the first $n$ coordinates.
Of course, there is  a natural transitive action of $G$ on $\{0,1\}^{n}$:  
$$
T_{\boldsymbol{s}}*(z_1,\dots,z_n):=\big(s_0z_1, s_1(z_1)z_2,\dots s_{n-1}(z_1,\dots, z_{n-1})z_n\big).
$$
Then  $\pi_n$  is a $G$-equivariant mapping.
Thus, $\{0,1\}^n$ is a finite factor of  $(Z,\eta, G)$.
Let 
$$
\Gamma_n:=\{g\in G\mid g*(0,\dots,0)=(0,\dots,0)\}.
$$
In other words, $\Gamma_n$ is the stabilizer of a point $(0,\dots,0)\in\{0,1\}^n$.
It is straightforward to verify that 
$$
\Gamma_n=\{T_{\boldsymbol{s}}\mid \boldsymbol{s}=(s_k)_{k=0}^\infty\text{ with } s_0= s_1(0)=s_2(0,0)=\cdots=s_{n-1}(0,\dots,0)=I\}.
$$
Hence, $\Gamma_n$ is a cofinite subgroup in $G$ for each $n$.
Moreover, 
$\Gamma_1\supsetneq\Gamma_2\supsetneq\cdots$
Thus, we obtain, for each $n$, a $G$-equivariant bijection $\phi_n:\{0,1\}^n\to
G/\Gamma_n$ such that the following diagram commutes:
$$
\CD
G/\Gamma_1 @<<< G/\Gamma_2 @<<< G/\Gamma_3 @<<<\cdots\\
@A\phi_1 AA @A\phi_2 AA    @A\phi_3 AA \\
\{0,1\} @<\tau_1<<     \{0,1\}^2  @<\tau_2<<    \{0,1\}^3 @<\tau_3<<\cdots
\endCD,
\tag3-5
$$
where $\tau_n(z_1,\dots,z_{n+1}):=(z_1,\dots,z_{n})$ for each $(z_1,\dots,z_{n+1})\in\{0,1\}^{n+1}$
and all $n\in\Bbb N$.
It is routine to check that
$$
\widetilde\Gamma_n:=\bigcap_{g\in G}g\Gamma_n g^{-1}=
\{T_{\boldsymbol{s}}\mid \boldsymbol{s}=(s_k)_{k=0}^\infty\text{ with } s_0=I, s_1\equiv I,s_2\equiv I,\dots, s_{n-1}\equiv I\}
$$
and $\bigcap_{n=1}^\infty \widetilde\Gamma_n=\{I\}$.
Denote by $(Y,O)$ the topological $G$-odometer associated  with the sequence $(\Gamma_n)_{n=1}^\infty$.
Furnish it with the Haar measure $\nu$.
Then \thetag{3-5} yields 
 a $G$-equivariant isomorphism $\phi:(Z,\eta)\to (Y,\nu)$, as desired.

We now show that $(Z,\eta,G)$ is not free.
Take a point $z=(z_n)_{n=1}^\infty\in Z$.
Then the $G$-stabilizer $G_z$ of $z$ is the group 
$$
\bigcap_{n=1}^\infty 
\{T_{\boldsymbol{s}}\mid \boldsymbol{s}=(s_k)_{k=0}^\infty\text{ with } s_0= s_1(z_1)=s_2(z_1,z_2)=\cdots=s_{n}(z_1\dots,z_n)=I\}.
$$
Let $r_1:\{0,1\}\to \text{Homeo}(\{0,1\})$ be the only mapping such that $r_1(z_1)=I$ but $r_1\not\equiv I$. 
We define a transformation $R$ of $(Z,\eta)$ by setting
$$
R(z_1,z_2,z_3,\dots):=(z_1,r_1(z_1)z_2, z_3,z_4,\dots).
$$
Then $G_z\ni R\ne I$.
Hence, $O$ is not free.
Therefore $O$ is not of rank one.
\endexample

\example{Example 3.6}
Let $(Z,\eta)$ and $\pi_n$ be as in Example~3.5.
Denote by $\Cal R$   the tail equivalence relation on $Z$.
We let
$$
{ \Cal A}:=\{T_{\boldsymbol{s}}\mid \text{for each $z\in Z$, there is $N>0$ with  $s_n(z_1,\dots,z_n)=I$ if $n>N$}\}.
$$
Then 
 ${ \Cal A}$ is a subgroup of $[\Cal R]$.
 Of course, ${ \Cal A}$ generates $\Cal R$.
Let $\theta$ denote the non-identity bijection of $\{0,1\}$.
Define 4 transformations $a,b,c,d\in \Cal A$  by the following formulae 
$$
\align
a(z_1,z_2,\dots) &:=(\theta(z_1),z_2,\dots)\\
b(1^n,0,z_{n+2},\dots)&:=
\cases
(1^n,0,\theta(z_{n+2}),z_{n+3},\dots) &\text{if }n\not\in 3\Bbb Z_+,\\
(1^n,0,z_{n+2},z_{n+3},\dots)
  &\text{otherwise}
  \endcases,
  \\
  c(1^n,0,z_{n+2},\dots)&:=
\cases
(1^n,0,\theta(z_{n+2}),z_{n+3},\dots) &\text{if }n\not\in 1+3\Bbb Z_+,\\
(1^n,0,z_{n+2},z_{n+3},\dots)
  &\text{otherwise}
   \endcases,
   \\
  d(1^n,0,z_{n+2},\dots)&:=
\cases
(1^n,0,\theta(z_{n+2}),z_{n+3},\dots) &\text{if }n\not\in 2+3\Bbb Z_+,\\
(1^n,0,z_{n+2},z_{n+3},\dots)
  &\text{otherwise}.
  \endcases
\endalign
$$
We remind that  the group $G$ generated by  $a,b,c,d$ is called the Grigorchuk group.
It was introduced in  \cite{Gr}.
The group is residually finite, amenable, non-locally finite. 
Every proper quotient subgroup of $G$ is finite.
Of course, $G\subset {\Cal A}$.
It is routine to verify $\Cal R$ is the  $G$-orbit equivalence relation.
Hence $G$ is an ergodic transformation group of $(Z,\eta)$.
Since $\pi_n$ is a $G$-equivariant mapping of $Z$ onto $\{0,1\}^n$ and
the dynamical system $(Z,\eta, G)$ is ergodic, $G$ acts transitively on $\{0,1\}^n$.
Therefore, repeating our reasoning in Example~3.5 almost literally, we obtain that
$(Z,\eta, G)$ is isomorphic to the  probability preserving $G$-odometer associated with the following sequence
$(\Gamma_n)_{n=1}^\infty$ of cofinite subgroups $\Gamma_n\subset G$: 
$$
\Gamma_n:=\{T_{\boldsymbol{s}}\in G\mid \boldsymbol{s}=(s_k)_{k=0}^\infty\text{ with } s_0= s_1(0)=\cdots=s_{n-1}(0,\dots,0)=I\},
$$
and  $\bigcap_{g\in G}\bigcap_{n=1}^\infty g\Gamma_ng^{-1}=\{I\}$.
Furthermore, the stabilizer $G_z$ of this odometer at a point $z=(z_n)_{n=1}^\infty\in Z$ is the group
$$
\{T_{\boldsymbol{s}}\in G\mid \boldsymbol{s}=(s_k)_{k=0}^\infty\text{ with } s_0= s_1(z_1)=s_2(z_1,z_2)=\cdots=s_n(z)=I\}.
$$
Hence $G$ 
 is  not free.
 Therefore $G$ is not of rank one.

\endexample

On the other hand, we will show that each probability preserving $G$-odometer is a factor of a rank-one $\sigma$-finite measure preserving $G$-action.

\proclaim{Theorem 3.7} 
Let $(Y,O)$ be a topological $G$-odometer associated with a decreasing sequence $\Gamma_1\supsetneq \Gamma_1\supsetneq\cdots$ of cofinite subgroups in $G$ satisfying \thetag{3-1}.
Then there is a topological $(C,F)$-action $T$ of $G$ on a locally compact Cantor space $X$
and a continuous $G$-equivariant mapping $\tau:X\to Y$.
Moreover, $\tau$ maps  the  Haar measure\footnote{See Remark~1.5.} on $X$ to a 
(non-$\sigma$-finite, in general) measure which is equivalent\footnote{This means that the two measures have the same class of subsets of zero measure.}
to the Haar measure on $Y$.
\endproclaim

\demo{Proof}
 Construct inductively sequences $(C_n)_{n=1}^\infty$, and
$(F_{n})_{n=0}^\infty$ of finite subsets in $G$ such that \thetag{1-1} and \thetag{1-4} holds,
 $C_n\subset \Gamma_n$ and the projection 
 $$
 C_n\ni c\mapsto c\Gamma_{n+1}\in \Gamma_n/\Gamma_{n+1}\qquad\text{is one-to-one and onto}\tag3-6
 $$
  for each $n\in\Bbb N$.
 Let $T$ be the topological $(C,F)$-action of $G$ associated with $(C_n,F_{n-1})_{n=1}^\infty$.
By Proposition~1.9, $T$ is defined on the entire locally compact space $X=\bigcup_{n=0}^\infty X_n$, where $X_n=F_n\times C_{n+1}\times C_{n+2}\times\cdots$.
We define $\tau:X\to Y$ by setting
$$
\tau(x)=(f_n\Gamma_1, f_n\Gamma_2, \dots, f_n\Gamma_{n+1},f_nc_{n+1}\Gamma_{n+2}, ,f_nc_{n+1}c_{n+2}\Gamma_{n+3},  \dots)\in Y
$$
if $x=(f_n,c_{n+1},c_{n+2},\dots)\in X_n$ for some $n\in X_n$.
Of course, $\tau$ is well defined,  continuous and $G$-equivariant. 
Thus, the first claim of the proposition is proved.

Let $\mu$ denote the Haar measure on $X$ and let $\chi$ be the  Haar measure on $Y$.
Then $\mu$ is the unique $T$-invariant $(C,F)$-measure such that $\mu(X_0)=1$.
It is determined by the sequence $(\kappa_n,\nu_{n-1})_{n=1}^\infty$, where
$\kappa_n$ is the equidistribution on $C_n$ and 
$\nu_n(f)=\prod_{k=1}^n\kappa_k(1_G)$ for each $f\in F_n$ and $n\in\Bbb N$.
We now show that 
$
(\mu\restriction X_0)\circ\tau^{-1}=\chi.
$
Let $\tau_n$ stand for the mapping
$$
C_1\times\cdots\times C_n\ni (c_1,\dots,c_n)\mapsto c_1\cdots c_n\Gamma_{n+1}\in G/\Gamma_{n+1}.
$$
We note that if 
$$
\tau_n(c_1,\dots, c_n)=\tau_n(c_1',\dots, c_n')
$$ for some
$c_1,c_1'\in C_1,\dots,c_n,c_n'\in C_n$ then $c_j=c_j'$ for each $j=1,\dots,n$.
Indeed,
$$
c_1\Gamma_2=c_1\cdots c_n\Gamma_{n+1}\Gamma_2=c_1'\cdots c_n'\Gamma_{n+1}\Gamma_2=c_1'\Gamma_2.
$$
Therefore \thetag{3-6} yields that $c_1=c_1'$ and hence $c_2\cdots c_n\Gamma_{n+1}=
c_2'\cdots c_n'\Gamma_{n+1}$.
Arguing in a similar way, we obtain that $c_2=c_2', \dots, c_n=c_n'$, as claimed.
It follows that $\tau_n$ is one-to-one.
Moreover, $\tau_n$ is onto in view of \thetag{3-5}.
Then it is straightforward to verify that the  diagram
$$
\CD
G/\Gamma_2 @<<<G/\Gamma_3 @<<< G/\Gamma_4 @<<<\cdots\\
@A{\tau_1}AA @A{\tau_2}AA @A{\tau_3}AA\\
C_1@<<< C_1\times C_2 @<<<  C_1\times C_2\times C_3@<<< \cdots\\
\endCD
$$
commutes.
Passing to the projective limit, we obtain that $\tau$ is a homeomorphism of $X_0$ onto $Y$.
Since $\tau_n$ maps the equidistribution on $C_1\times\cdots\times C_n$ to the 
equidistribution on  $G/\Gamma_{n+1}$ for each $n$, it follows that
$\tau$ maps $\mu\restriction X_0$ to $\chi$.
Take  a  probability measure $\mu'$ on $X$ which is is equivalent to $\mu$.
Then  
$$
\mu'\circ\tau^{-1}\gg (\mu'\restriction X_0)\circ\tau^{-1}\sim (\mu\restriction X_0)\circ\tau^{-1}=\chi.\tag3-7
$$
Since $\mu'$ is quasiinvariant and ergodic under $T$, it follows that $\mu'\circ\tau^{-1}$
is quasiinvariant and ergodic under $O$.
As the two probability Borel measures $\mu'\circ\tau^{-1}$ and $\chi$ on $Y$ are quasiinvariant and ergodic under $O$, they are either equivalent or mutually singular.
Therefore \thetag{3-7} yields that $\mu'\circ\tau^{-1}\sim\chi$.
Hence $\mu\circ\tau^{-1}$ is equivalent to $\chi$, as desired.
\qed
\enddemo

\remark{Remark 3.8} If we change the construction of $T$ in the proof of Theorem~3.7 in such a way that the mapping \thetag{3-6} is one-to-one but $\# (C_n)/\#(\Gamma_n/\Gamma_{n+1})\le  0.5 $ for each $n\in\Bbb N$ then the first claim of Theorem~3.7 still holds: there is a $G$-equivariant continuous mapping $\tau:X\to Y$.
However, the second claim fails: the $O$-quasi-invariant measure $\mu\circ\tau^{-1}$ on $Y$ will be  singular with $\chi$.
\endremark

\vskip 5pt

We note that if $G$ is Abelian then  each  ergodic nonsingular $G$-action $T$  possesses the following property. 
Let $\Gamma_1\subsetneq\Gamma_2\subsetneq\cdots$ be a sequence of cofinite subgroups in $G$ with $\bigcap_{n\in\Bbb N}\Gamma_n=\{0_G\}$.
If, for each $n\in\Bbb N$,  $T$ has a finite factor $\goth F_n$ isomorphic to the homogeneous $G$-space $G/\Gamma_n$ then $\goth F_1\subsetneq\goth F_2\subsetneq\cdots$ and  $T$ has an odometer factor $\bigvee_{n>0}\goth F_n$.
This is no longer true if $G$ is non-Abelian (see Example~2.5).
However, the following  version of the aforementioned property holds for an arbitrary $G$.

\proclaim{Theorem~3.9} 
Let $T=(T_g)_{g\in G}$ be an ergodic  nonsingular $G$-action on a standard nonatomic probability space $(X,\goth B,\mu)$.
Let  $\Gamma_1\supsetneq\Gamma_2\supsetneq\cdots$ be a sequence of cofinite subgroups in $G$ such that $\thetag{3-1}$ holds.
Denote by $(Y,O)$ the topological $G$-odometer associated with this sequence.
Suppose that for  each $n\in\Bbb N$, there exists a $T$-factor which is isomorphic to the homogeneous $G$-space $G/\Gamma_n$.
 Then
there is an  $O$-quasi-invariant measure $\nu$ on $Y$ such that the nonsingular odometer $(Y,\nu,O)$ is a factor of $(X,\mu,T)$.
\endproclaim

\demo{Proof}
We first prove an auxiliary claim.

{\it Claim A.}
Let $H$ be a cofinite subgroup in $G$.
There exist no more than $\# (G/H)$ different factors of $T$ that are
 isomorphic to the homogeneous $G$-space $G/H$.

 Let $J:=\# (G/H)+1$.
Suppose that there are $J$ pairwise different $T$-invariant $\sigma$-algebras $\goth F_j\subset\goth B$
such that $T\restriction\goth F_j$ is isomorphic to $G/H$ for each $j\in J$.
Denote by $\tau_j:X\to G/H$ the corresponding $G$-equivariant mapping.
Then the mapping 
$$
\tau:X\ni x\mapsto (\tau_j(x))_{j\in J}\in (G/H)^J
$$
 is also $G$-equivariant.
Denote by  $E$  the support of the measure $\mu\circ\tau^{-1}$.
Then
\roster
\item"---"  $E$ is a single $G$-orbit;
\item"---" the projection of $E$ onto each of the $J$ coordinates is onto.
\endroster
Take a point $(g_jH)_{j\in J}\in E$.
Since $J>\#(G/H)$, there are $i_0, j_0\in J$ such that $i_0\ne j_0$ but $g_{i_0}H=g_{j_0}H$.
Hence the projection of $E$ onto the ``plane'' generated by the $i_0$ and $j_0$ coordinates is
the diagonal $\{(gH,gH)\mid g\in G\}$ in $(G/H)^2$.
Hence, $\goth F_{i_0}=\goth F_{j_0}$.
This contradiction proves Claim~A.

To prove the theorem we
 define a graded graph $\Cal G$.
The set $V$ of vertices of $\Cal G$ is the union $\bigsqcup_{n\ge 0}V_n$, where 
$V_n$ is the set of all Borel $G$-equivariant maps from $X$ to $G/\Gamma_n$.
For the consistency of  notation, we  let $\Gamma_0:=G$.
Given $n\ge 0$, we denote by $\theta_n$ the projection 
$$
G/\Gamma_{n+1}\ni g\Gamma_{n+1}\mapsto g\Gamma_n\in G/\Gamma_n.
$$
The set $E$ of edges of $\Cal G$ is the union $\bigsqcup_{n\ge 0}E_n$, where an edge $e\in E_n$ joins 
a vertex $\pi\in V_n$ with a vertex $\tau\in V_{n+1}$ if 
$\pi=\theta_n\circ\tau$, 
It follows from Claim~A that $V_n$ is finite for each $n$.
Of course, every vertex from $V_n$ is adjacent (i.e. connected by an edge) with a vertex in $V_{n-1}$ for each $n\in\Bbb N$.
Hence, for each vertex of $\Cal G$, there is a path connecting this vertex with the only vertex from $V_0$.
Thus, $\Cal G$ is connected.
Of course, $\Cal G$ is locally finite and infinite.
Hence, by the K\H{o}nig's infinity lemma, $\Cal G$ contains a ray.
It follows that there exists  a Borel $G$-equivariant mapping $\iota:X\to Y$.
We set $\nu:=\mu\circ\iota^{-1}$.
Then $(Y,\nu,O)$ is a factor of $(X,\mu,T)$, as desired.
\qed
\enddemo

\subhead 3.2. Normal covers for nonsingular odometers
\endsubhead
Let $(Y,O)$ be a topological  $G$-odometer associated with
 a decreasing sequence $(\Gamma_n)_{n=1}^\infty$ of  cofinite subgroups in $G$ such that \thetag{3-1} holds.
  If each $\Gamma_n$ is normal in $G$ then $(Y,O)$ is 
 called {\it normal}.
In this case we have that $G/\Gamma_n$ is a finite group and hence   $Y$ is a compact totally disconnected metric group.
Moreover, there is a  one-to-one group homomorphism $\phi:G\to Y$ such that $O_gy=\phi(g)y$ for all $g\in G$ and $y\in Y$.
Of course, $\phi(g)=(g\Gamma_1,g\Gamma_2,\dots)\in Y$ for each $g\in G$.
This homomorphism embeds $G$ densely into $Y$.
Every normal odometer is free.

Given a cofinite subgroup $\Gamma$ in $G$,
 the subgroup $\widetilde \Gamma:=\bigcap_{g\in G}g\Gamma g^{-1}$ is the maximal normal (in $G$) subgroup of $\Gamma$.
Of course, $\widetilde \Gamma$ is of finite index in $G$.
The natural projection 
$G/\widetilde\Gamma\ni g\widetilde\Gamma\mapsto g\Gamma\in G/\Gamma$ is  $G$-equivariant.
Hence, for  a decreasing sequence $\Gamma_1\supset\Gamma_2\supset\cdots$ of cofinite subgroups in $G$ satisfying \thetag{3-1}, we obtain a  decreasing sequence
$\widetilde\Gamma_1\supset\widetilde\Gamma_2\supset\cdots$ 
of normal cofinite subgroups in $G$ with $\bigcap_{n=1}^\infty\widetilde\Gamma_n=\{1_G\}$.
Let 
$(\widetilde Y, \widetilde O)$ denote the normal topological  $G$-odometer 
associated with $(\widetilde\Gamma_n)_{n=1}^\infty$.
It is called {\it the topological normal cover} of $(Y,O)$.
The natural projections 
$$
G/\widetilde\Gamma_n\ni g\widetilde\Gamma_n\mapsto g\Gamma_n\in G/\widetilde\Gamma_n,\quad  n\in\Bbb N,
$$
 generate a  continuous projection $w:\widetilde Y\to Y$ that intertwines $\widetilde O$ with $O$.
Let 
$$
H:=\{(\widetilde y_n)_{n=1}^\infty\in \widetilde Y\mid \widetilde y_n\in\Gamma_n/\widetilde\Gamma_n\text{ for all }n\in\Bbb N\}.
$$
Then $H$ is a closed subgroup of $\widetilde Y$.
We claim that $\omega$ is the quotient mapping 
$$
\widetilde Y\ni\widetilde y\mapsto\widetilde yH\in\widetilde Y/H.
$$
Indeed, we first observe that $\omega(\widetilde y)=\omega(\widetilde yh)$ for all $\widetilde y\in \widetilde Y$ and $h\in H$.
Secondly, if $\omega(\widetilde y)=\omega(\widetilde z)$ for some $\widetilde y,\widetilde z\in\widetilde Y$ then $\widetilde y\widetilde z^{-1}\in H$.
Finally, the subset $\omega(\widetilde Y)$ is $G$-invariant and closed in $Y$.
Hence $\omega(\widetilde Y)=Y$.

It may seem that  the coordinate projection $H\ni (\widetilde y_n)_{n=1}^\infty\mapsto \widetilde y_n\in \Gamma_n/\widetilde\Gamma_n$ is onto for each $n\in\Bbb N$.
That is not true. 
A counterexample (in which $H$ is trivial  but $\#(\Gamma_n/\widetilde\Gamma_n)=2$ for each $n$) is constructed in Example~5.7 below.

\definition{Definition 3.10} Let  $(Y,\nu, O)$ be a nonsingular $G$-odometer,    $(\widetilde Y,\widetilde O)$ the topological normal cover  of $(Y,O)$, and $\widetilde\nu$ an $\widetilde O$-quasi-invariant probability on $\widetilde Y$.
We call the nonsingular normal odometer $(\widetilde Y,\widetilde\nu,\widetilde O)$  {\it the normal cover} of  $(Y,\nu,O)$ if 
\roster
\item"(i)" $\widetilde\nu\circ\omega^{-1}=\nu$ and 
\item"(ii)" $\frac{d\widetilde\nu\circ\widetilde O_g}{d\widetilde\nu}=\frac{d\nu\circ O_g}{d\nu}\circ\omega$  \ for each $g\in G$ and
\endroster
\enddefinition

We note that (ii) means that  $\widetilde O$ is $\omega$-relatively finite measure preserving.

   \proclaim{Proposition 3.11} Given a nonsingular   $G$-odometer  $(Y,\nu,O)$, there is a   $G$-quasi-invariant probability $\widetilde\nu$ on
$\widetilde Y$ such that $(\widetilde Y,\widetilde\nu,\widetilde O)$ is a normal cover
of  $(Y,\nu,O)$.
\endproclaim
\demo{Proof}  Without loss of generality we may assume that $\widetilde Y=Y\times H$ (as a set, not as a group) and there is a Borel map (1-cocycle) $s:G\times Y\to H$ such that
$$
\widetilde O_g(y,h)=(O_gy, s(g,y)h)\text{  \ and \ $\omega(y,h)=y$}\quad\text{for each $(y,h)\in Y\times H$}.
$$
Denote by $\lambda_H$ the Haar measure on $H$.
Then the direct product $\widetilde\nu:=\nu\otimes\lambda_H$ satisfies (i) and (ii) from Definition~3.10.
\qed

\comment

We remind that two nonsingular $G$-actions are {\it weakly equivalent} if each one of them is a factor of the other.

\proclaim{Proposition 3.9} Let two $G$-odometers furnished with Haar measures are weakly equivalent. 
Then they are isomorphic.
\endproclaim

\endcomment

\comment

For brevity, we will write below $gy$ and $g\widetilde y$ for $O_gy$ and $\widetilde O_g\widetilde y$ respectively.
Let $s:Y\to \widetilde Y$ be a Borel cross-section of the natural projection $\omega$.
Then for each $g\in G$ and $y\in Y$,
 $$
 s(gy)h(g,y)=gs(y)\quad\text{for some  \ $h(g,y)\in H$.}
 $$
Denote by $\lambda_H$ the Haar measure on $H$.
For a Borel subset $A\subset \widetilde Y$ and each $y\in Y$, we let $A_y:=A\cap\omega^{-1}(y)$
and set 
$$
\widetilde\nu(A):=\int_{Y}\lambda_H\big(s(y)^{-1} A_y\big)\,  d\nu(y)
$$
Since $\omega^{-1}(y)=s(y)H$, we obtain that $s(y)^{-1} A_y\subset H$.
Hence $\widetilde\nu(A)$ is well defined.
Of course, $\widetilde\nu$ is a well defined probability on $\widetilde Y$.
Moreover, 
 (i) holds.
 Given $g\in G$ and $y\in Y$, we obtain that
 $$
 (gA)_y=gA\cap s(y)H=g(A\cap g^{-1}s(y)H)=g(A\cap s(g^{-1}y)H)=gA_{g^{-1}y}.
 $$
 Therefore
 $$
 \align
 \widetilde\nu(gA)&=\int_{Y}\lambda_H\big(s(y)^{-1} gA_{g^{-1}y}\big)\,  d\nu(y)\\
 &=
 \int_{Y}\lambda_H\big(s(g^{-1}y)^{-1} A_{g^{-1}y}\big)\,  d\nu(y)\\
 &=
 \int_{Y}\lambda_H\big(s(y)^{-1} A_{y}\big)  \frac{d\nu\circ O_g}{d\nu}\,d\nu(y).
 \endalign
 $$
Thus yields (ii).

\endcomment
\enddemo

Let $\Gamma_1\supsetneq\Gamma_2\supsetneq\cdots$ be as above.
Suppose that there is a sequence $(b_n)_{n=1}^\infty$ of $G$-elements  such that
$$
b_1\Gamma_1b_1^{-1}\supsetneq b_2\Gamma_2b_2^{-1}\supsetneq\cdots.
$$
Denote by $(Y', O')$ the topological odometer associated with this sequence.
Let $\nu$ and $\nu'$ stand for the Haar measures on $Y$ and $Y'$ respectively.

\proclaim{Proposition 3.12} The odometers $(Y',O',\nu')$ and  $( Y, O,\nu)$ are isomorphic.
\endproclaim
\demo{Proof} 
It is easy to see 
that the normal covers  of $(Y',O',\nu')$ and  $( Y, O,\nu)$ are same.
Denote this common normal cover by $(\widetilde Y,\widetilde \nu,\widetilde O)$.
Then there are closed subgroups $H$ and $H'$ of $\widetilde Y$ such that
$( Y, O,\nu)$ is the right $H$-quotient of $(\widetilde Y,\widetilde \nu,\widetilde O)$ and 
$( Y', O',\nu')$ is the right $H'$-quotient of $(\widetilde Y,\widetilde \nu,\widetilde O)$.
It follows from Theorem~3.9 that $(Y',O',\nu')$ and  $( Y, O,\nu)$ are weakly equivalent, i.e.
$(Y',O',\nu')$ is a factor of  $( Y, O,\nu)$  and  $( Y, O,\nu)$ is a factor of $(Y',O',\nu')$.
Hence, there are compact subgroups $K$ and $K'$ of $\widetilde Y$  and
elements $a,b\in \widetilde Y$ such that
$K\supset H$, $K'\supset H'$, $K=aH'a^{-1}$ and $K'=bHb^{-1}$.
Hence, $H\subset abHb^{-1}a^{-1}$.
We claim  that this implies that  $H=abHb^{-1}a^{-1}$.
Indeed, let $V:=\{y\in\widetilde Y\mid H\subset yHy^{-1}\}$.
Then $V$ is a closed subset of $\widetilde Y$.
Of course, $V\ni (ab)^n$ for each $n\in\Bbb N$.
Hence $V$ includes the closure of the semigroup $\{(ab)^n\mid n\in\Bbb N\}$.
It follows from \cite{HeRo, Theorem~9.1} that 
the closure of $\{(ab)^n\mid n\in\Bbb N\}$ equals the closure of the group $\{(ab)^n\mid n\in\Bbb Z\}$.
Hence $(ab)^{-1}\in V$, i.e. $H\supset abHb^{-1}a^{-1}$.
Therefore, $H=abHb^{-1}a^{-1}$, as claimed.
This yields that $K= H$ and  $K'=H'$.
Thus, we obtain that $H$ and $H'$ are conjugate.
Hence, $(Y',O',\nu')$ and  $( Y, O,\nu)$ are isomorphic.
\qed
\enddemo

\comment

We now let
$$
\align
H&:=\{(\widetilde y_n)_{n=1}^\infty\in \widetilde Y\mid \widetilde y_n\in\Gamma_n/\widetilde\Gamma_n\text{ for all }n\in\Bbb N\}\quad\text{and}\\
H'&:=\{(\widetilde y_n)_{n=1}^\infty\in \widetilde Y\mid \widetilde y_n\in\Gamma_n'/\widetilde\Gamma_n=(b_n\Gamma_nb_n^{-1})/\widetilde\Gamma_n\text{ for all }n\in\Bbb N\}
\endalign
$$
Suppose that the  $G$-odometers $(Y,O,\nu)$ and $(Y',O',\nu')$ are isomorphic.
Then the corresponding isomorphism extends by continuity to an isomorphism
of  the homogeneous $\widetilde Y$-spaces $\widetilde Y/H$ and $\widetilde Y/H'$.
This, in turn, yields that $H$ and $H'$ are conjugate in $\widetilde Y$.
Thus, 
  there is  an element 
$y=(y_n)_{n=1}^\infty\in \widetilde Y$ such that  $yHy^{-1}=H'$.
Let $H_n:=\{z_n\in \Gamma_n/\widetilde\Gamma_n\mid (z_m)_{m=1}^\infty\in H\}$
and let
$\widetilde b_n:=b_n\widetilde \Gamma_n\in\Gamma_n/\widetilde\Gamma_n$,
for each $n\in\Bbb N$.
Then
$$
y_n H_ny_n^{-1}=
\widetilde b_nH_n\widetilde b_n^{-1},
$$
for each $n\in\Bbb N$.
Thus, $(Y,O,\nu)$ and $(Y',O',\nu')$ are isomorphic if and only if
there is $(y_n)_{n=1}^\infty\in \widetilde Y$ such that
 $$
 y_n{\widetilde b_n}^{-1}\in N_{G/\widetilde{\Gamma}_n}(H_n)\quad\text{ for each $n\in\Bbb N$.}
 $$
 We note that $H_n$ can be  a proper subgroup of $\Gamma_n/\widetilde\Gamma_n$.
 See Example~4.6 below, where  $H$ is trivial but $\#(\Gamma_n/\widetilde \Gamma_n)=2$ for each $n>0$.
 
On the other hand, if $(Y,O,\nu)$ and $(Y',O',\nu')$ are weakly equivalent, i.e.
$(Y,O,\nu)$ is a factor  of  $(Y',O',\nu')$ and  $(Y',O',\nu')$ is a factor of 
$(Y,O,\nu)$ then $H$ and $H'$ are conjugate.
Hence, $(Y,O,\nu)$ and $(Y',O',\nu')$ are  isomorphic.

\endcomment

\head 4. Odometer factors of nonsingular $(C,F)$-actions
\endhead 
The following concept is an ``infinite''  analogue of Definition~2.1.

\definition{Definition 4.1}
Let a sequence $\Cal T=(C_n,F_{n-1},\kappa_n,\nu_{n-1})_{n=1}^\infty$ satisfy  \thetag{1-1}--\thetag{1-3} and Proposition~1.10(ii) and let a sequence $(\Gamma_n)_{n=1}^\infty$ satisfy \thetag{3-1.}
Denote by $(Y,O)$ the topological odometer associated with $(\Gamma_n)_{n=1}^\infty$.
Given $y=(g_n\Gamma_n)_{n=1}^\infty\in Y$,
we say that 
$\Cal T$  is {\it compatible with $y$} if
$$
 \sum_{n=1}^\infty\kappa_n\big(\{c\in C_n\mid c\not\in g_n\Gamma_n g_n^{-1}\}\big)<\infty.
 $$ 
 \enddefinition

 Denote by $T$ the $(C,F)$-action of $G$  associated with $\Cal T$.
 Let $X$ be the space of $T$ and let $\mu$ stand for the 
 the nonsingular $(C,F)$-measure on $X$ determined by $(\kappa_n)_{n=1}^\infty$ and $(\nu_n)_{n=0}^\infty$. 
As 
$$
g_n\Gamma_n g_n^{-1}=g_{n+1}\Gamma_n g_{n+1}^{-1}\supset g_{n+1}\Gamma_{n+1} g_{n+1}^{-1}\qquad\text{for each $n\in\Bbb N$,}
$$ 
 it follows that  if $\Cal T$  is  compatible with $y$ then
 $\Cal T$  is  compatible with the coset $g_n\Gamma_n\in G/\Gamma_n$ in the sense of Definition~2.1 for each $n\in\Bbb N$.
 Hence the $( {\Cal T},g_n\Gamma_n)$-factor mapping  $\pi_{({\Cal T},g_n\Gamma_n)}:X\to G/\Gamma_n$ for $T$ is well defined (mod 0) for each $n\in\Bbb N$.
 Moreover, a  measurable mapping 
 $$
 \pi_{({\Cal T},y)}:=X\ni x\mapsto \Big(\pi_{({\Cal T},g_n\Gamma_n)}(x)\Big)_{n=1}^\infty\in Y
 $$ 
 is well defined (mod 0) too.
 Of course, $\pi_{({\Cal T},y)}\circ T_g=O_g\circ \pi_{({\Cal T},y)}$ for each $g\in G$.
 Hence the nonsingular odometer $\Big(Y,\mu\circ\pi_{({\Cal T},y)}^{-1}, O\Big)$ is a factor
 of $(X,\mu, T)$.

\definition{Definition 4.2} We call $\pi_{({\Cal T},y)}$ {\it the $({\Cal T},y)$-factor mapping}
for $T$.
\enddefinition

In the  proposition below we find necessary and sufficient conditions (in terms of the parameters $\Cal T$) under which $\pi_{({\Cal T},y)}$ is one-to-one, i.e. the dynamical systems $(X,\mu,T)$ and  $(Y,\mu\circ \pi_{({\Cal T},y)}^{-1},O)$ are isomorphic via  $\pi_{({\Cal T},y)}$.

\proclaim{Proposition 4.3} Let $\Cal T$ be compatible with $y$.
The following are equivalent:
\roster
\item"(i)"
$\pi_{({\Cal T},y)}$ is one-to-one (mod 0),
\item"(ii)"
for each $n>0$ and $\epsilon>0$ there are $l>0$ and a subset $D_l\subset G/\Gamma_l$ such that $\mu\Big([1_G]_n\triangle \pi_{({\Cal T},g_l\Gamma_l)}^{-1}(D_l)\Big)<\epsilon$ and 
\item"(iii)"
for each $n>0$ and $\epsilon>0$ there are $l>0$, a subset $D_l\subset G/\Gamma_l$ and $M>0$ such that
$\nu_m\Big(C_{n+1}\cdots C_m\triangle\{f\in F_m\mid fg_l\Gamma_l\in D_l\} \Big)<\epsilon$ for each $m>M$.
\endroster
\endproclaim

\demo{Proof} (i)$\Leftrightarrow$(ii) Denote by $\goth B$ and $\goth Y$ the Borel $\sigma$-algebra
 on $X$ and $Y$ respectively.
 Of course, $\pi_{({\Cal T},y)}$ is one-to-one (mod 0) if and only if $\goth B=\pi_{({\Cal T},y)}^{-1}(\goth Y)$ (mod 0).
 Since 
 \roster
 \item"---" $\goth B$ is generated by the family of all cylinders in $X$ and
 \item"---" $\goth B$ is invariant under $T$,
 \endroster
 it follows that $\goth B=\pi_{({\Cal T},y)}^{-1}(\goth Y)$ if and only if $[1_G]_n\in\pi_{({\Cal T},y)}^{-1}(\goth Y)$
 for each $n>0$.
 Let $\goth Y_l\subset\goth Y$ denote the finite sub-$\sigma$-algebra of subsets that are measurable with respect to the canonical projection $G\to G/\Gamma_l$.
 Then $\goth Y_1\subset\goth Y_2\subset\cdots$  and the union $\bigcup_{l=1}^\infty\goth Y_l$ is dense in $\goth Y$.
It follows that $[1_G]_n\in\pi_{({\Cal T},y)}^{-1}(\goth Y)$ if and only if
$$
\min_{A\in\goth Y_l}\mu\Big([1_G]_n\triangle \pi_{({\Cal T},y)}^{-1}(A)\Big)\to 0\quad\text{ as $l\to\infty$.}
$$
It remains to note that $\{ \pi_{({\Cal T},y)}^{-1}(A)\mid A\in\goth Y_l\}=\{\pi_{({\Cal T},g_l\Gamma_l)}^{-1}(D)\mid D\subset G/\Gamma_l\}$.

(ii)$\Leftrightarrow$(iii) 
We note that for all $n>0$, $l>0$ and a subset $D\subset G/\Gamma_l$,
$$
\gather
\mu\Big([1_G]_n\triangle \pi_{({\Cal T},g_l\Gamma_l)}^{-1}(D)\Big)=
\lim_{m\to\infty}\mu\Big([1_G]_n\triangle \pi_{({\Cal T},g_l\Gamma_l)}^{-1}(D)\cap[F_m]_m\Big)\quad\text{and}\\
\lim_{m\to\infty}\bigotimes_{j\ge m}\kappa_j\Big(\big\{x=(c_j)_{j\ge m}\in C_m\times C_{m+1}\times\cdots\mid c_j\in g_l\Gamma_lg_l^{-1}\  \forall j\ge m\big\}\Big)=1.
\endgather
$$
The latter follows from the fact that $\Cal T$ is compatible with $y$.
It  implies that
$$
\lim_{m\to\infty}\mu\Big(\{x=(f_m,c_{n+1},\dots)\in [F_m]_m\mid \pi_{({\Cal T},g_l\Gamma_l)}(x)=f_mg_l\Gamma_l \}\Big)=\mu([F_m]_m).
$$
Hence, for each $D\subset \Gamma_l$,
$$
\align
\mu\Big([1_G]_n\triangle \pi_{({\Cal T},g_l\Gamma_l)}^{-1}(D)\Big)&=
\lim_{m\to\infty}\mu\Bigg([1_G]_n\triangle\bigsqcup_{f\in F_m,fg_l\Gamma_l\in D}[f]_m\Bigg)\\
&=\lim_{m\to\infty}\nu_m\Big(C_{n+1}\cdots C_m\triangle\big\{f\in F_m\mid fg_l\Gamma_l\in D\big\} \Big).
\endalign
$$
This equality implies the equivalence of (ii) and (iii).
 \qed

 \enddemo

The following theorem is the main result of this section.

\proclaim{Theorem 4.4}
Let a sequence $\Cal T=(C_n,F_{n-1},\kappa_n,\nu_{n-1})_{n=1}^\infty$ satisfy  \thetag{1-1}--\thetag{1-3} and Proposition~1.10(ii).
Let $T$  be  the nonsingular $(C,F)$-action of $G$  associated with $\Cal T$
and let  $(Y,O)$ be the topological $G$-odometer associated with a sequence $(\Gamma_n)_{n=1}^\infty$ satisfying \thetag{3-1}. 
Then for each  $G$-equivariant measurable mapping $\tau:X\to Y$
 there exist an increasing sequence $\boldsymbol q$ of nonnegative integers and an element
 $y\in Y$ 
 such that the $\boldsymbol q$-telescoping $\widetilde{\Cal T}=
 (\widetilde C_n,\widetilde F_{n-1},\widetilde\kappa_n,\widetilde\nu_{n-1})_{n=1}^\infty$ of
 $\Cal T$ is compatible with $y$
 and  
 $$
  \pi_{(\widetilde{\Cal T},y)}\circ \iota_{\boldsymbol q}=\tau.\tag4-1
  $$
 Moreover, $\tau$ is one-to-one (mod 0) if and only if 
 for each $n>0$ and $\epsilon>0$ there are $l>0$, a subset $D_l\subset G/\Gamma_l$ and $M>n$ such that  for each $m>M$,
 $$
 \widetilde\nu_m\Big(\widetilde C_{n+1}\cdots \widetilde C_m\triangle\{f\in \widetilde F_m\mid fg_l\Gamma_l\in D_l\} \Big)<\epsilon.
 $$
\comment
 Under these conditions, the factor map $\widetilde\pi:\widetilde X\to Y
$ of $\widetilde T$ is as follows:
 $$
 \gather
\widetilde\pi(x)_n:=\lim_{m\to\infty}f_{n}c_{n+1}\cdots c_{q_m} g_0^{-1}g_1^{-1}\cdots
 g_{n-1}^{-1}\Gamma_n\in G/\Gamma_n,\tag3-4\\
 \widetilde\pi(x):=(\widetilde\pi(x)_n)_{n=1}^\infty\in Y.
 \endgather
 $$
 at $\widetilde\mu$-a.e. $x=(f_n, c_{n+1},c_{n+2},\dots)\in\widetilde X_n\subset \widetilde X$ for each $n>0$.
 \endcomment
\endproclaim

We preface the proof of Theorem~4.4 with  auxiliary simple but useful facts about factor mappings.

\proclaim{Lemma 4.5} Let $\Gamma,\Gamma_1$ be  two cofinite subgroups in $G$ and $\Gamma_1\subset \Gamma$.
Then
\roster
\item"(i)" if  $\Cal T$ is compatible with a coset $g\Gamma\in G/\Gamma$  then for each increasing sequence  $\boldsymbol a=(a_n)_{n=0}^\infty$ of nonnegative integers with $a_0=0$, the $\boldsymbol a$-telescoping $\widetilde{\Cal T}$ of $\Cal T$ is also compatible with $g\Gamma$ and
 $$\pi_{(\Cal  T,g\Gamma)}=\pi_{(\widetilde {\Cal T},g\Gamma)}\circ\iota_{\boldsymbol a},
 \tag4-2
 $$
 \item"(ii)"
 if  $\Cal T$ is compatible with two cosets $g\Gamma\in G/\Gamma$ and
 $g_1\Gamma_1\in G/\Gamma_1$ 
 then 
 $$
r\circ \pi_{(\Cal T,g_1\Gamma_1)}=\pi_{(\Cal T,g\Gamma)}\quad\text{if and only if}\quad g_1g^{-1}\in\Gamma,
 $$
where $r: G/\Gamma_1\to
 G/\Gamma$ denotes the natural projection.
 \endroster
 \endproclaim

\demo{Proof}
(i) 
For each $n> 0$,  we let $\widetilde C_{n}:=C_{a_{n-1}+1}\cdots C_{a_n}$ and
$\widetilde\kappa_n:=\kappa_{a_{n-1}}*\cdots *\kappa_{a_{n}}$.
Then
$$
\widetilde \kappa_n
\big(\{c\in \widetilde C_n\mid c\not\in g\Gamma g^{-1}\}\big)\le \sum_{j=a_{n-1}+1}^{a_{n}}\kappa_j\big(\{c\in C_j\mid c\not\in g\Gamma g^{-1}\}\big).
$$
Hence,
$$
\sum_{n=1}^\infty\widetilde \kappa_n
\big(\{c\in \widetilde C_n\mid c\not\in g\Gamma g^{-1}\}\big)\le \sum_{j=1}^\infty\kappa_j\big(\{c\in C_j\mid c\not\in g\Gamma g^{-1}\}\big)<\infty.
$$
Thus, $\widetilde{\Cal T}$ is compatible with  $g\Gamma$.

The formula \thetag{4-2} and the claim (ii) are verified straightforwardly. 
\qed
\enddemo

\demo{Proof of Theorem~4.4} 
We note that $\tau(x)=(\tau_n(x))_{n=1}^\infty$ for each $x\in X$, where
 $\tau_n:X\to G/\Gamma_n$ is  a $G$-equivariant mapping for  every $n$.

From now on, we will argue inductively.
At the first step  we apply~Remark~2.4 to $T$ and $\tau_1$: there exist a coset $g_1\Gamma_1\in G/\Gamma_1$ and an
increasing sequence $\boldsymbol q^1$ of nonnegative integers such that the $\boldsymbol q^1$-telescoping
 $\Cal T_1$ of $\Cal T$ is
 $g_1\Gamma_1$-compatible  and
$$
\pi_{(\Cal T_1,g_1\Gamma_1)}\circ \iota_{\boldsymbol q^1}=\tau_1.\tag4-3
$$
At the second step we we apply~Remark~2.4 to the $(C,F)$-action $ \widetilde T$ of $G$ associated with $\Cal T_1$ and the factor mapping $\tau_2\circ{\iota_{\boldsymbol q^1}}^{-1}$ of $ \widetilde T$: there exist
 a coset $g_2\Gamma_2\in G/\Gamma_2$ and
and an
increasing sequence $\boldsymbol q^2$ of nonnegative integers such that the $\boldsymbol q^2$-telescoping
 $\Cal T_2$ of $\Cal T_1$ is
 $g_2\Gamma_2$-compatible  and
$$
\pi_{(\Cal T_2,g_2\Gamma_2)}\circ \iota_{\boldsymbol q^2}=\tau_2\circ {\iota_{\boldsymbol q^1}}^{-1}.
\tag4-4
$$
Consider the natural projection $\omega_{2,1}:G/\Gamma_2\to G/\Gamma_1$.
Since $\omega_{2,1}\circ\tau_2=\tau_1$, it follows from~\thetag{4-3} and~\thetag{4-4} that 
$$
\omega_{2,1}\circ\pi_{(\Cal T_2,g_2\Gamma_2)}\circ \iota_{\boldsymbol q^2}=\pi_{(\Cal T_1,g_1\Gamma_1)}
$$
Then~Lemma~4.5(i),(ii) imply that $g_2g_1^{-1}\in\Gamma_1$.
Continuing inductively, we obtain a sequence $(g_n)_{n=1}^\infty$ of elements in $G$ and a sequence $(\boldsymbol q^n)_{n=1}^\infty$ of increasing sequences of nonnegative integers such that for each $n>0$,
\roster
\item"($\alpha_1$)"
$g_ng_{n-1}^{-1}\in \Gamma_{n-1}$,  
\item"($\alpha_2$)"
the $\boldsymbol q^n$-telescoping $\Cal T_n=(C_k^{(n)},F_{k-1}^{(n)},\kappa_k^{(n)},\nu_{k-1}^{(n)})_{k=1}^\infty$  of $\Cal T_{n-1}$ is compatible with the coset $g_{n}\Gamma_n$
 and
\item"($\alpha_3$)" $\pi_{(\Cal T_n,g_{n}\Gamma_n)}\circ \iota_{\boldsymbol q^1\circ\cdots\circ \boldsymbol q^n}=\tau_n$.
\endroster 
We now choose an integer $a_n>0$ large so that
\roster
\item"($\alpha_4$)"
$
\sum_{k\ge a_n}\kappa_k^{(n)}(\{c\in C_k^{(n)}\mid c\not\in g_n\Gamma_n g_n^{-1} \})<\frac 1{n^2}
$
\endroster
for each $n>0$.
It follows from ($\alpha_1$) that the sequence $y:=(g_{n}\Gamma_n)_{n=1}^\infty$ is a well-defined element of $Y$.
Of course, 
\roster
\item""$\boldsymbol q^1\circ\boldsymbol q^2$ is a subsequence of $\boldsymbol q^1$,
\item""
$\boldsymbol q^1\circ\boldsymbol q^2\circ\boldsymbol q^3$ is a subsequence of $\boldsymbol q^1\circ\boldsymbol q^2$
\endroster
 and so on.
Hence, utilizing  the diagonalization method, we can construct a sequence $\boldsymbol q=(q_n)_{n=1}^\infty$ of integers such that
\roster
\item"---" $0=q_0<q_1<\cdots <q_n$,
\item"---"
 $\boldsymbol q$  is a subsequence of $\boldsymbol q^1\circ\cdots\circ\boldsymbol q^n$  and
 \item"---"
 $q_n\ge a_n$
 \endroster
 for every $n>0$.
 Denote by  $\widetilde{\Cal T}$
 the $\boldsymbol q$-telescoping
 $\widetilde{\Cal T}=(\widetilde C_n,\widetilde F_{n-1},\widetilde \kappa_n,\widetilde\nu_{n-1})_{n=1}^\infty$ of $\Cal T$.
 We are going to show that  $\widetilde{\Cal T}$ is compatible with $y$.
 By the construction of  $\boldsymbol q$, for each $n>0$, there are integers
 $
 d_{n,2}\ge d_{n,1}\ge a_n
 $
such that  
$$
\widetilde C_n=C^{(n)}_{d_{n,1}}\cdots C^{(n)}_{d_{n,2}}\qquad\text{and}\qquad
\widetilde\kappa_n=\kappa_{d_{n,1}}*\cdots *\kappa_{d_{n,2}}.
$$
 Therefore we deduce from  ($\alpha_4$) that
 $$
\widetilde\kappa_n(\{c\in \widetilde C_n\mid c\not\in g_n\Gamma_n g_n^{-1} \})
\le\sum_{k=d_{n,1}}^{d_{n,2}}
\kappa_k^{(n)}(\{c\in C_k^{(n)}\mid c\not\in g_n\Gamma_n g_n^{-1} \})
<\frac 1{n^2}.
 $$
Hence, $\sum_{n=1}^\infty\widetilde\kappa_n(\{c\in \widetilde C_n\mid c\not\in g_n\Gamma_n g_n^{-1} \})<\infty$, i.e. $\widetilde{\Cal T}$ is compatible with $y$, as desired.

 We now prove \thetag{4-1}.
Of course, 
 $\widetilde{\Cal T}$
  is  a telescoping of the $\boldsymbol q^1\circ\cdots\circ\boldsymbol q^n$-telescoping of $\Cal T$ for each $n$.
  In view of   ($\alpha_2$),   the $\boldsymbol q^1\circ\cdots\circ\boldsymbol q^n$-telescoping of $\Cal T$ equals $\Cal T_n$.
  Thus, $\widetilde{\Cal T}$ is a telescoping of $\Cal T_n$.
  Denote by $\theta_n$ the canonical isomorphism corresponding to this telescoping.
 Then  $\iota_{\boldsymbol q}=\theta_n\circ \iota_{\boldsymbol q^1\circ\cdots\circ\boldsymbol q^n}$.
It follows from this and ($\alpha_3$) that
 $$
 \pi_{(\widetilde{\Cal T},g_{n}\Gamma_n)}\circ\iota_{\boldsymbol q}=
  \pi_{(\widetilde{\Cal T},g_{n}\Gamma_n)}\circ\theta_n\circ\iota_{\boldsymbol q^1\circ\cdots\circ\boldsymbol q^n}
  =\pi_{({\Cal T_n},g_{n}\Gamma_n)}\circ \iota_{\boldsymbol q^1\circ\cdots\circ\boldsymbol q^n}
  =\tau_n
 $$
 for each $n\in\Bbb N$.
 Hence $\pi_{(\widetilde{\Cal T},y)}=\tau$, as desired.

 The second (the last) claim of the theorem follows from the first one and Proposition~4.3.
 \qed

\enddemo

As a corollary, we obtain a criterion for the  existence (or non-existence) of odometer factors for rank-one nonsingular actions.

\proclaim{Corollary 4.6}
Let $T$  be  the nonsingular $(C,F)$-action of $G$  associated with  a sequence  $\Cal T=(C_n,F_{n-1},\kappa_n,\nu_{n-1})_{n=1}^\infty$ satisfying  \thetag{1-1}--\thetag{1-3} and Proposition~1.10(ii).
Then $T$ has no nonsingular odometer factors if and only if
for each decreasing sequence $\Gamma_1\supsetneq\Gamma_2\supsetneq\cdots$ of cofinite subgroups in $G$ satisfying \thetag{3-1}, no  telescoping of 
$ \Cal T$ is  compatible with the sequence $(\Gamma_n)_{n=1}^\infty$, i.e.
 for each sequence $0=q_1<q_2<\cdots$, 
$$
\sum_{n=0}^\infty
\kappa_{q_n+1}*\cdots*\kappa_{q_{n+1}}(\{c\in C_{q_n+1}\cdots C_{q_{n+1}}\mid c\not\in\Gamma_n  \})=\infty.
$$
\endproclaim

\demo{Proof} It is sufficient to utilize Theorem~4.4 and the following remark.

Let $Y$ be a  $G$-odometer associated  with a decreasing sequence $(\Gamma_n)_{n=1}^\infty$ of cofinite subgroups in $G$ such that \thetag{3-1} holds.
Let $y\in Y$.
Then $y=(g_n\Gamma_n)_{n=1}^\infty$ with $g_ng_{n+1}^{-1}\in \Gamma_n$ for each $n\in\Bbb N$.
Of course, $g_1\Gamma_1g_1^{-1}\supset g_2\Gamma_2g_2^{-1}\supset \cdots$
and the sequence
$(g_n\Gamma_ng_n^{-1})_{n=1}^\infty$ satisfies \thetag{3-1}.
Denote by $Y_y$ the space of the $G$-odometer associated with 
$(g_n\Gamma_ng_n^{-1})_{n=1}^\infty$.
Then there is a canonical $G$-equivariant homeomorphism $\varphi_y:Y\to Y_y$.
It is well defined by the formula
$$
\varphi_y\big((z_n\Gamma_n)_{n=1}^\infty\big):=(z_n\Gamma_ng_n^{-1})_{n=1}^\infty=
\big(z_ng_n^{-1}(g_n\Gamma_ng_n^{-1})\big)_{n=1}^\infty.
$$

It follows that  there is a $G$-equivariant map from $X$ to $Y$ if and only if there is  
 $G$-equivariant map from $X$ to $Y_y$.
\qed
\enddemo

In a similar way, we obtain a criterion when a nonsingular $(C,F)$-action is not isomorphic to any nonsingular odometer.

\proclaim{Corollary 4.7}
Let $T$  be  the nonsingular $(C,F)$-action of $G$  associated with  a sequence  $\Cal T$ satisfying  \thetag{1-1}--\thetag{1-3} and Proposition~1.10(ii).
Then $T$ is not isomorphic to any nonsingular odometer if and only if 
for each decreasing  sequence $\Gamma_1\supsetneq\Gamma_2\supsetneq\cdots$ of cofinite subgroups in $G$ satisfying \thetag{3-1} 
and each increasing  sequence  $\boldsymbol q$ of nonnegative integers such that
the $\boldsymbol q$-telescoping $\widetilde {\Cal T}$ of $\Cal T$ is compatible with the point
$(\Gamma_n)_{n=1}^\infty\in\projlim_{n\to\infty}G/\Gamma_n$, 
there exist $n>0$ and $\epsilon_0>0$
such that for each $l>0$,  $D_l\subset G/\Gamma_l$ and $M>n$, there is $m>M$ with 
$$
 \widetilde\nu_m\Big(\widetilde C_{n+1}\cdots \widetilde C_m\triangle\{f\in \widetilde F_m\mid fg_l\Gamma_l\in D_l\} \Big)>\epsilon_0.
 $$
\endproclaim

We state one more corollary from Theorems~4.4 and Theorem~1.19 on the existence of minimal Radon uniquely ergodic topological models for  rank-one nonsingular extensions of nonsingular odometers.

\proclaim{Corollary 4.8} Let $(X,\mu,T)$ be a rank-one nonsingular action of $G$.
Let $T$ have a nonsingular odometer factor $(Y,\nu, O)$ and let $\pi:X\to Y$ stand for the corresponding $G$-equivariant factor mapping with $\nu=\mu\circ\pi^{-1}$. 
Then there exist a locally compact Cantor space $\widetilde X$, a minimal Radon uniquely ergodic 
free continuous action $\widetilde T$ of $G$ on $\widetilde X$, a continuous $G$-equivariant mapping $\widetilde\pi:\widetilde X\to Y$ and a Borel isomorphism  $R:X\to 
\widetilde X$ such that 
\roster
\item"---" $\widetilde\mu:=\mu\circ R^{-1}$ is a Radon measure on $\widetilde X$, 
\item"---" $RT_g=\widetilde T_gR$ for each $g\in G$,
\item"---" the function $\frac{d\widetilde\mu\circ\widetilde T_g}{d\widetilde\mu}:\widetilde X\to\Bbb R^*$ is continuous for each $g\in G$,
 \item"---" $\widetilde T$ is Radon $\big(\frac{d\widetilde\mu\circ\widetilde T_g}{d\widetilde\mu}\big)_{g\in G}$-uniquely ergodic
and
\item"---" $\widetilde \pi R=\pi$.
\endroster
\endproclaim

We can also characterize the class of quasi-invariant measures for odometers that appear as factors of rank-one actions.
Let $(Y,O)$ be the topological $G$-odometer associated with a decreasing sequence $(\Gamma_n)_{n=1}^\infty$ of cofinite subgroups $\Gamma_n$ of $G$ satisfying \thetag{3-1}.
Let  a sequence $\Cal T=(C_n,F_{n-1},\kappa_n,\nu_{n-1})_{n=1}^\infty$ satisfy \thetag{1-1}--\thetag{1-4} and $C_n\subset\Gamma_n$ for each $n>0$.
Denote by $\mu_{\Cal F}$  the $(C,F)$-measure determined by the sequence $(\kappa_n,\nu_{n-1})_{n=1}^\infty$.
Let $X$ stand for the space of $\mu_\Cal F$.
We define a mapping $\pi_{\Cal F}:X\to Y$ by setting
$$
\pi_{\Cal F}(x)=(f_n\Gamma_1,\dots, f_n\Gamma_{n+1},f_nc_{n+1}\Gamma_{n+2},
f_nc_{n+1}c_{n+2}\Gamma_{n+3},\dots)
$$
if $x=(f_n,c_{n+1},c_{n+2},\dots)\in X_n\subset X$ for some $n\ge 0$.
Then $\pi_{\Cal F}$ is well defined  and continuous. 
Let
$$
\Cal M_{Y}:=\{\mu_\Cal F\circ\pi^{-1}_{\Cal F}\mid \text{$\Cal F$ satisfies  \thetag{1-1}--\thetag{1-4} and $C_n\subset\Gamma_n$ $\forall n>0$}\}.
$$
We deduce the following claim from Corollary~4.8.

\proclaim{Corollary 4.9}
Each measure $\nu\in \Cal M_{Y}$ is quasininvariant under $O$.
A Borel measure $\nu$ on $Y$ is equivalent to a measure belonging to $\Cal M_{Y}$ if and only if there is a rank-one nonsingular $G$-action $T$ such that $(Y,O,\nu)$ is a measurable factor of $T$.
\endproclaim

\head 5.  Examples\endhead

\subhead 5.1. Non-odometer  rank-one $\Bbb Z$-action with
odometer factor
\endsubhead
In \cite{Fo--We}, an example of classical rank-one finite measure preserving $\Bbb Z$-action $T$ is constructed such that 
\roster
\item"---" $T$ has the 2-adic odometer as a factor but 
\item"---" $T$ is not isomorphic to any odometer.
\endroster
We remind that given a prime $p$, the {\it p-adic odometer} is associated with the sequence 
$$
p\Bbb Z\supset p^2\Bbb Z\supset p^3\Bbb Z\supset\cdots
$$
 of cofinite subgroups in $\Bbb Z$.
The argument  in \cite{Fo--We} is based on their description of the odometer factors of rank-one transformations (that result is generalized in our Theorem~4.4).
We now consider their example from another point of view, bypassing the use of any version of Theorem~4.4.
Our approach is more direct and  leads to  stronger results.

\example{Example 5.1} Let $G=\Bbb Z$.
We construct a measure preserving (classical) rank-one $\Bbb Z$-action $T$ on a probability space $(X,\mu)$
 such that 
 \roster
 \item"---" 
 $(X,\mu,T)$ has  a proper 2-adic odometer  factor $(Y,\nu,O)$,
 \item"---"
 $(Y,\nu,O)$ is the {\it Kronecker factor} of $(X,\mu,T)$, i.e. $O$ is the maximal  factor of $T$ with a pure discrete spectrum,
 \item"---" 
 the projection $(X,\mu)\to (Y,\nu)$ is  uncountable-to-one (mod 0), i.e. the corresponding conditional measures on fibers are non-atomic.
\endroster
We set $h_0:=0$ and $h_{n+1}:=4h_n+2^{n+1}$ for each $n\in\Bbb N$.
It follows that  $h_n=2^n(2^{n+1}-1)$ for each $n\ge 0$. 
We let 
$$
\gather
F_n:=\{0,\dots,h_n-1\},  \quad C_{n+1}:=\{0,h_n,2h_n+2^{n+1},3h_n+2^{n+1}\},\\
\nu_n(f)=\frac 1{4^n}\quad
\text{for each $f\in F_n$ and $\kappa_n(c)=\frac 14$ for each $c\in C_{n+1}$}
\endgather
$$
 for every $n\ge0$.
 Then the sequence $\Cal T:=(C_n,F_{n-1},\kappa_n,\nu_{n-1})_{n=1}^\infty$
 satisfies \thetag{1-1}--\thetag{1-3} and   Proposition~1.10(ii).
 Denote by $(X,\mu, T)$ the $(C,F)$-action of $\Bbb Z$ associated with $\Cal T$.
Then $T$ is of classical rank one along $(F_n)_{n=0}^\infty$.
We note that $T$ is the transformation that was studied in \cite{Fo--We}.
Of course,  $T$ preserves $\mu$
 and $\mu(X)<\infty$.
 We have that
$$
\mu(X)=\mu(X_0)+\sum_{n=1}^\infty 2^n\mu([0]_n)=1+\sum_{n=1}^\infty \frac{2^n}{4^n}=2.
$$
Denote by $(Y, O)$ the 2-adic $\Bbb Z$-odometer.
Then the transformation $O_1$ acts on the compact metric group $Y:=\projlim_{n\to\infty} \Bbb Z/2^n\Bbb Z$ by translation with the element $(1+2\Bbb Z, 1+2^2\Bbb Z, 1+2^3\Bbb Z,\dots)\in Y$.
Let $\nu$ stand for the Haar  measure on $Y$.
Since each element of $C_n$ is divisible by $2^{n-1}$ for every $n>0$, it follows 
that $O$ is a factor of $T$.
The corresponding $(\Cal T, 0)$-factor mapping $\pi:X\to Y$ is well defined by the formula:
$$
\pi(x)=(f_n+ 2^{n}\Bbb Z, f_n+c_{n+1}+2^{n+1}\Bbb Z,f_n+c_{n+1}+c_{n+2}+2^{n+2}\Bbb Z,\cdots)\in Y,
$$
if  $x=(f_n,c_{n+1},c_{n+1},\dots)\in X_n=F_n\times C_{n+1}\times C_{n+2}\times\cdots\subset X$ for some  $n\ge 0$ (see Definition~2.2). 
Since the measure  $\mu\circ\pi^{-1}$ is invariant under $O$, it follows that $\mu\circ\pi^{-1}$ is proportional to $\nu$.
More precisely,  $\mu\circ\pi^{-1}=\mu(X)\cdot \nu=2\nu$.

We now show that if $\lambda$ is an eigenvalue of $T$ then there is $n>0$ such that $\lambda^{2^n}=1$.
Since $\# C_m=4$ for each $m$, it follows from  \cite{DaVi, Corollary~3.8} that
$$
\lim_{m\to\infty}\max_{c\in C_m}|1-\lambda^c|=0.
$$
As
$
4^{m}=2h_{m-1}+2^m\in C_m,
$
we obtain that $\lim_{m\to\infty}\lambda^{4^m}=1$.
This is only possible if $\lambda$ is a dyadic root of $1$, as desired.
\footnote{Indeed, observe that $\lambda^{4^{m+1}}=\big(\lambda^{4^m}\big)^4$ and iterate.}
On the other hand, denote by $\goth Y$ the $\sigma$-algebra of all  measurable subsets in $Y$.
Let $\goth F:=\{\pi^{-1}(B)\mid B\in\goth Y\}$.
Then each eigenfunction of $T$ whose eigenvalue is a 2-adic root of $1$ is $\goth F$-measurable.
 It follows that the  $\goth F$
is the Kronecker factor of $T$.

We now show that the Kronecker factor is proper, i.e. that the spectrum of $T$ has a continuous
component.
Moreover, we prove that the extension $T\to O$ is uncoutable-to-one.
For each $n>0$, we let 
$C_n^{(1)}:=\{0,h_{n-1}\}$ and $C_n^{(2)}:=\{0,2h_{n-1}+2^n\}=\{0,4^n\}.$
Then $C_n=C_n^{(1)}+C_n^{(2)}$.
For $j=1,2$, let $X_0^{(j)}:=C_1^{(j)}\times C_2^{(j)}\times\cdots$.
Then $X_0^{(j)}$ is a compact subset of $X_0$.
Given $x=(c_1,c_2,\dots)\in X_0^{(1)}$ and $z=(d_1,d_2,\dots)\in X_0^{(2)}$, the sum
$$
x+z:=(c_1+d_1,c_2+d_2, \dots)\in X_0
$$
 is well defined.\footnote{Note that $X_0^{(1)}$, $X_0^{(2)}$ and $X_0$ are compact subsets of the Polish Abelian group $\Bbb Z^\Bbb N$.} 
Moreover, the mapping
$ (x,z)\mapsto x+z
$
is a homeomorphism of the Cartesian product $X_0^{(1)}\times X_0^{(2)}$ onto $X_0$ and 
 $$
 \pi(x+z)=\pi(x)+\pi(z)\qquad\text{
for all  $x,z\in X_0$.}\tag5-1
 $$
Endow $X_0^{(1)}$ and $X_0^{(2)}$ with the infinite products $\mu^{(1)}$ and $\mu^{(2)}$ of the equidistributions on  $C_n^{(1)}$ and $C_n^{(2)}$ respectively, $n\in\Bbb N$.
We claim that the restriction of $\pi$ to $X_0^{(j)}$ is one-to-one for $j=1,2$.
It is straightforward to verify that for each $x=(c_m)_{m=1}^\infty\in X_0$,
$$
\pi(x)=\bigg(\bigg(\sum_{m=1}^nc_m^1+\sum_{1\le m< n/2}c_m^2\bigg)+2^n\Bbb Z\bigg)_{n=1}^\infty\in Y,
\tag5-2
$$
where $c_m^j\in C_m^{(j)}$ and $c_m=c_m^1+c_m^2$ for each $m$.
Take  two points $x=(c_1^1,c_2^1,\dots)\in X_0^{(1)}$ and $y=(d_1^1,d_2^1,\dots)\in X_0^{(1)}$
 such that
$\pi(x)=\pi(y)$. 
It follows from \thetag{5-2} that $c_1^1=d_1^1\pmod 2$, $c_1^1+c_2^1=d_1^1+d_2^1\pmod {2^2}$, \dots.
 The first equality
implies that $d_1^1=c_1^1$.
 Therefore, the second equality 
is equivalent to $c_2^1=d_2^1\pmod {2^2}$,   which, in turn, yields that $d_2^1=c_2^1$.
By the induction, $d^1_n=c^1_n$ for each $n>0$, i.e. $x=y$, as desired.
It now follows that the mapping 
$$
C^{(1)}_1\times\cdots\times C^{(1)}_n\ni(c_1^1,\dots,c_n^1)\mapsto\bigg(\sum_{m=1}^n c_m^1\bigg)+2^n\Bbb Z\in \Bbb Z/2^n\Bbb Z
$$
 is bijective for each $n>0$.
Hence $\pi$ maps $X_0^{(1)}$ bijectively (and homeomorphically)  onto $Y$.
Moreover, $\pi$ maps $\mu^{(1)}$ to $\nu$.
In a similar way, one can check that  $\pi$ maps $X_0^{(2)}$ bijectively (and homeomorphically) onto the closed subset $\pi(X_0^{(2)})\subset Y$.
Given $y\in Y$ and $z\in X_0^{(2)}$, we define an element $Q(y,z)\in X_0^{(1)}$ by the formula
$$
\pi(Q(y,z))=y-\pi(z).
$$
Of  course, $Q(y,z)$ is well defined.
Moreover, the mapping 
$$
Q:Y\times X_0^{(2)}\ni (y,z)\mapsto Q(y,z)\in X_0^{(1)}
$$
 is continuous.
 If we fix $y\in Y$ then the mapping $ X_0^{(2)}\ni z\mapsto Q(y,z)\in X_0^{(1)}$ is one-to-one.
 Of course, for each $y\in Y$,
 $$
X_0\cap  \pi^{-1}(y)=\{Q(y,z)+z\mid z\in X_0^{(2)}\}.
 $$
Hence, $\pi^{-1}(y)$ is uncountable.
Moreover, the corresponding conditional measure on the fiber $X_0\cap  \pi^{-1}(y)$ is the 
image of $\mu^{(2)}$ under the mapping 
$$
X_0^{(2)}\ni z\mapsto Q(y,z)+z\in X_0.
$$
Hence, the conditional  measure on $X_0\cap  \pi^{-1}(y)$ is non-atomic for each $y\in Y$.
In particular, $\pi$ is uncoutable-to-one.
\comment

The coordinates $d_n\in C_n$ will be determined inductively in $n$.
For each $n>0$, we let $C_n^{(1)}:=\{0,h_{n-1}\}$ and $C_n^{(2)}:=\{0,2h_{n-1}+2^n\}=\{0,4^n\}$.
Then $C_n=C_n^{(1)}+C_n^{(2)}$.
Hence each element $c\in C_n$ can be written as a sum $c=c^1+c^2$ for some
uniquely determined $c^1\in C_n^{(1)}$ and $c^2\in C_n^{(2)}$.
Therefore there are (uniquely determined)  4 sequences $(c_n^j)_{n=1}^\infty$ and $(d_n^j)_{n=1}^\infty$, $j=1,2$, such that
$c^j_n,d^j_n\in C^{(j)}_n$, $c_n=c^1_n+c^2_n$ and  $d_n=d_n^1+d_n^2$
for all $n\in\Bbb N$ and $j=1,2$.
Given $c\in C^{(j)}_n$, we define the {\it complimentary} element $\widehat c\in C^{(j)}_n$ by the formula
$C^{(j)}_n=\{c,\widehat c\}$.
It is straightforward to verify that
$$
\align
\pi(x)&=(c_1^1+2\Bbb Z, c_1^1+c_2^1+2^2\Bbb Z,c_1+c_2^1+c_3^1+2^3\Bbb Z, c_1+c_2^1+c_3^1+c^1_4+2^4\Bbb Z,\dots)\in Y,\\
\pi(y)&=(d_1^1+2\Bbb Z, 
d_1^1+d_2^1+2^2\Bbb Z,
d_1+d_2^1+2^3\Bbb Z,d_1+d_2^1+d_3^1+2^4\Bbb Z,\dots)\in Y.
\endalign
$$
As $\pi(x)=\pi(y)$, we obtain that 
 $$
 \align
 c_1^1&=d_1^1\pmod 2, 
 \\ c_1^1+c_2^1&=d_1^1+d_2^1\pmod {2^2}, \\ 
 c_1+c_2^1+c_3^1&=d_1+d_2^1+d_3^1\pmod {2^3}, \\
c_1+c_2^1+c_3^1+c_4^1&=d_1+d_2^1+d_3^1+d_4^1\pmod {2^4},\\
c_1+c_2+c_3^1+c_4^1+c_5^1&=d_1+d_2+d_3^1+d_4^1+d_5^1\pmod {2^5},\dots.
 \endalign 
 $$
The first equality  $c_1^1=d_1^1\pmod 2$ implies that $d_1^1=c_1^1$.
 Hence, the second equality $c_1^1+c_2^1=d_1^1+d_2^1\pmod {2^2}$
is equivalent to $c_2^1=d_2^1\pmod {2^2}$,   which, in turn, yields that $d_2^1=c_2^1$.
 Hence, the third equality $c_1+c_2^1+c_3^1=d_1+d_2^1+d_3^1\pmod {2^3}$ is equivalent to
 $c_1^2+c_3^1=d_1^2+d_3^1\pmod {2^3}$.
This implies that either  $(d_1^2,d_3^1)=(c_1^2,c_3^1)$ or 
$(d_1^2,d_3^1)=(\widehat{c_1^2},\widehat{c_3^1})$.
We now pass to an arbitrary $n$.
Consider separately two cases: $n$ is odd, and $n$ is even.

Let first  
$$
\sum_{j=1}^nc_j+\sum_{j=n+1}^{2n+1}c_{j}^1=
\sum_{j=1}^nd_j+\sum_{j=n+1}^{2n+1}d_{j}^1
\pmod{ 2^{2n+1}}
$$
for some $n\ge1$.
Our purpose on this step is to determine 
$d_{2n+2}^1$.
We have that  either
$$
\gather
\sum_{j=1}^nc_j+\sum_{j=n+1}^{2n+1}c_{j}^1=
\sum_{j=1}^nd_j+\sum_{j=n+1}^{2n+1}d_{j}^1
\pmod{ 2^{2n+2}}\quad\text{or}\tag4-1\\
\sum_{j=1}^nc_j+\sum_{j=n+1}^{2n+1}c_{j}^1=\bigg(
\sum_{j=1}^nd_j+\sum_{j=n+1}^{2n+1}d_{j}^1\bigg)+2^{2n+1}
\pmod{ 2^{2n+2}}.
\tag4-2
\endgather
$$
If \thetag{4-1} holds  then we let
 $d_{2n+2}^1:=c_{2n+2}^1$.
 If \thetag{4-2}  holds then we let
 $d_{2n+2}^1:=\widehat{c_{2n+2}^1}$.
 It is straightforward to verify that in the two cases,
 $$
\sum_{j=1}^nc_j+\sum_{j=n+1}^{2n+2}c_{j}^1=
\sum_{j=1}^nd_j+\sum_{j=n+1}^{2n+2}d_{j}^1
\pmod{ 2^{2n+2}},
\tag4-3
$$
 as desired.

Consider now the second case: \thetag{4-3} holds for some $n\ge1$.
Our purpose on this step is to determine $d_{n+1}^2$ and
$d_{2n+3}^1$ simultaneously.
We have that  either
$$
\gather
\sum_{j=1}^nc_j+\sum_{j=n+1}^{2n+2}c_{j}^1=
\sum_{j=1}^nd_j+\sum_{j=n+1}^{2n+2}d_{j}^1
\pmod{ 2^{2n+3}}\quad\text{or}\tag4-4\\
\sum_{j=1}^nc_j+\sum_{j=n+1}^{2n+2}c_{j}^1=\bigg(
\sum_{j=1}^nd_j+\sum_{j=n+1}^{2n+2}d_{j}^1\bigg)+2^{2n+2}
\pmod{ 2^{2n+3}}.
\tag4-5
\endgather
$$
If \thetag{4-4} holds, we have two (exactly two) possible choices for the pair $(d_{n+1}^2,
d_{2n+3}^1 )$:
$$
\text{either  $(d_{n+1}^2,d_{2n+3}^1 ):=(c_{n+1}^2,
c_{2n+3}^1 )$ or $(d_{n+1}^2,d_{2n+3}^1 ):=(\widehat{c_{n+1}^2},
\widehat{c_{2n+3}^1} )$}.
$$
If \thetag{4-5} holds, we also have two (exactly two) possible choices for  $(d_{n+1}^2,
d_{2n+3}^1 )$:
$$
\text{either  $(d_{n+1}^2,d_{2n+3}^1 ):=(\widehat{c_{n+1}^2},
c_{2n+3}^1 )$ or $(d_{n+1}^2,d_{2n+3}^1 ):=(c_{n+1}^2,
\widehat{c_{2n+3}^1} )$}.
$$
It is straightforward to verify that in every of the four possible (and only in these four) choices,  
 $$
 \align
\bigg(\sum_{j=1}^nc_j\bigg)+c_{n+1}^2+\sum_{j=n+1}^{2n+3}c_{j}^1&=
\bigg(\sum_{j=1}^nd_j\bigg)+d_{n+1}^2+\sum_{j=n+1}^{2n+3}d_{j}^1
\pmod{ 2^{2n+3}},\quad\text{i.e.}\\
\sum_{j=1}^{n+1}c_j+\sum_{j=n+2}^{2n+3}c_{j}^1&=
\sum_{j=1}^{n+1}d_j+\sum_{j=n+2}^{2n+3}d_{j}^1
\pmod{ 2^{2n+3}},
\endalign
$$
 as desired.



\endcomment

\comment

We now claim that $T$ is not isomorphic to any transformation with pure discrete spectrum.
Indeed, suppose that  $T$ has pure discrete spectrum.
 Then $T$ is isomorphic to a rotation on a monothetic compact group $K$ equipped with the Haar measure.
Moreover, there exists a closed subgroup $K_0\subset K$ such that $\pi$ is the projection
$$
\pi:K\ni k\mapsto kK_0\in K/K_0
$$ 
and $K/K_0=Y$.
As $\pi$ is 2-to-1, $\# K_0=2$.
Hence $K_0=\Bbb Z/2\Bbb Z$ and we obtain a short exact sequence of compact metric Abelian groups
$$
\{0\}\to \Bbb Z/2\Bbb Z\to K\to Y\to\{0\},
$$
where $Y$ is the group of 2-adic integers.
The only possible $K$ with this property is the direct group product $(\Bbb Z/2\Bbb Z)\times Y$
with the natural embedding of $\Bbb Z/2\Bbb Z$ into $(\Bbb Z/2\Bbb Z)\times Y$ and the natural projection of $(\Bbb Z/2\Bbb Z)\times Y$ onto $Y$.
However, $(\Bbb Z/2\Bbb Z)\times Y$ is not monothetic, a contradiction.

Since each 2-to-1 extension of a dynamical system is a group extension (more precisely, $\Bbb Z/2\Bbb Z$-extension) of this system, it follows from \cite{Ki} that $T$ is rigid.

\endcomment
\endexample

\remark{Remark 5.2} \roster
\item"(i)"
The existence of  the 2-adic odometer factor $O$ of  $T$  in Example~5.1 was proved in \cite{Fo--We}.
It was shown there that  $O$ is maximal in the family of  odometer factors of $T$.
We refine this result by showing  in Example~5.1  that $O$ is  the Kronecker  factor of $T$.
The claims that  $T$ is an uncountable-to-one extension of $O$, and  $T$ has a continuous part in the spectrum are new.
\item"(ii)" Perhaps, the simplest example of non-odometer rank-one $\Bbb Z$-action with the maximal odometer factor $O=(O_n)_{n\in\Bbb Z}$ is the following one.
Let $R$ be an irrational rotation on the circle. 
Then the transformation $S := O_1 \times R$ is an ergodic transformation with pure discrete spectrum. 
Hence it is of rank one \cite{dJ1}. 
Of course, $O_1$ is the maximal odometer factor of $S$. 
Of course, $S$ is not isomorphic to any odometer (as the discrete spectrum of $S$ has elements  of infinite order).
Obviously, $S$ is rigid.
However, in contrast with $T$ from Example~5.1, the spectrum of $S$ is purely discrete.
\endroster
\endremark

\subhead 5.2. Nonsingular counterparts of Example 5.1
\endsubhead
We first remind briefly the concepts of the associated flow of an ergodic equivalence relation, Krieger's type and AT-flow. 
For details we refer to the survey \cite{DaSi} and references therein.
Let $\Cal R$ be a countable Borel equivalence relation on a standard $\sigma$-finite measure space $(Z,\gamma)$.
Assume that $\Cal R$ is $\gamma$-nonsingular and ergodic.
This means that if a Borel subset $A\subset Z$ is $\gamma$-null then the $\Cal R$-saturation of $A$ is also $\gamma$-null, and each $\Cal R$-invariant (i.e. $\Cal R$-saturated) Borel subset of $Z$ is either $\gamma$-null or $\gamma$-conull. 
Denote by $\rho_\gamma:\Cal R\to\Bbb R^*_+$ the Radon-Nikodym cocycle of $\Cal R$.
Endow the product space $Z\times \Bbb R$ with the product measure $\gamma\otimes\text{Leb}$.
We define an equivalence relation $\Cal R(\log\rho_\gamma)$ on $Z\times \Bbb R$ by setting
$$
   (z,t)\sim_{\Cal R(\log\rho_\gamma) } (z',t')\qquad\text{if $(z,z')\in\Cal R$ and $t'=t-
   \log\rho_\gamma(z,z')$.}
$$
Then $\Cal R(\log\rho_\gamma)$ is countable, $(\mu\otimes\text{Leb})$-nonsingular but not necessarily ergodic. 
Denote by $\goth I$ the $\sigma$-algebra of $\Cal R(\log\rho_\gamma)$-invariant subsets.
Let $V=(V_s)_{s\in \Bbb R}$ denote the action of $\Bbb R$ on $Z\times \Bbb R$ by translations along the second coordinate, i.e. $V_s(z,t)=(z,t+s)$.
Of course,  $\goth I$ is invariant under $V$.
The dynamical system $(Z\times\Bbb R, \goth I, \mu\otimes\text{Leb},V)$ is called {\it the flow associated with $\Cal R$}.
The associated flow is nonsingular and ergodic.
We denote it by $W^\Cal R$.
We will need the following two well known properties of the associated flows:
\roster
\item"$(*)$"
If $A\subset Z$ is of positive measure $\gamma$ then the associated flow of $\Cal R\restriction A$
is isomorphic to the associated flow of $\Cal R$.
\item"$(**)$"
If $(Z,\gamma)=(Z_1,\gamma_1)\otimes(Z_2,\gamma_2)$, $\Cal R=\Cal R_1\otimes\Cal R_2$
and $W^{\Cal R_1}$ is free and transitive then $W^{\Cal R}$ is isomorphic to $W^{\Cal R_2}$.
\endroster
An ergodic nonsingular flow $V$ is called an {\it AT-flow} if there is a sequence $(A_n,\alpha_n)_{n=1}^\infty$ of finite subsets $A_n$ and non-degenerated probability measures $\alpha_n$ on $A_n$ such that the infinite product measure $\bigotimes_{n=1}^\infty\alpha_n$ is non-atomic and $V$ is isomorphic to the associated flow of the tail equivalence relation on the probability space
$\bigotimes_{n=1}^\infty(A_n,\alpha_n)$.
If, moreover, $\# A_n=2$ for each $n$ then we call the corresponding AT-flow {\it finitary}\footnote{The class of finitary AT-flows coincides with the class of flows of weights for the ITPFI$_2$ factors (in the sense of the theory of von Neumann algebras).}.
For instance, every  ergodic flow with pure point spectrum is a finitary AT-flow \cite{BeVa}.
If $S$ is an ergodic nonsingular $G$-action on a standard $\sigma$-finite measure space then the 
associated flow of the $S$-orbit equivalence relation is called {\it the associated flow of $S$}.
If the associated flow of $S$ is transitive and aperiodic then $S$ is called {\it of Krieger type $II$.}
If the associated flow of $S$ is transitive and periodic with the least positive period $-\log \lambda$ for some $\lambda\in (0,1)$ then $S$ is called {\it of Krieger type $III_\lambda$.}
If the associated flow of $S$ is the trivial flow on a singleton then $S$ is called {\it of Krieger type $III_1$.}
If $S$ is neither of type $II$ nor of type $III_\lambda$ for any $\lambda\in (0,1]$ then $S$
is called {\it of Krieger type $III_0$.}

Let $(X,T)$, $(Y,\nu,O)$ and $\pi:X\to Y$ be the same as in Example~5.1.

\proclaim{Proposition 5.3}
For each finitary AT-flow $V$, 
there exists a $(C,F)$-measure $\mu$ on $X$ such that 
\roster
\item"(i)" $\mu$ is quasiinvariant under $T$ and hence the nonsingular system $(X,\mu, T)$ is of rank one,
\item"(ii)"  the associated flow of the system $(X,\mu,T)$ is isomorphic to $V$,
\item"(iii)"  $\mu\circ\pi^{-1}\sim\nu$, i.e. 
the 2-adic probability preserving odometer $(Y,\nu,O)$ is a factor of $(X,\mu,T)$,
\item"(iv)"
$(Y,\nu,O)$ is the maximal odometer factor of $(X,\mu,T)$  and 
\item"(v)" the factor mapping $\pi$ is uncoutable-to-one (mod 0).
Hence, $(X,\mu,T)$ is not an odometer.
\endroster
In particular, for each $\lambda\in[0,1]$, there is a  $(C,F)$-measure $\mu$ on $X$ such that $(X,\mu,T)$ is  of Krieger type $III_\lambda$  satisfying \rom{(i), (iii)--(v)}.
\endproclaim

\demo{Proof} 
We will use below the notation from \S5.1.

Let $\kappa_n^1$ stand for the equidistribution on $C_n^{(1)}$ and let
$\kappa^1:=\bigotimes_{n=1}^\infty\kappa_n^1$.
As $V$ is finitary AT, there is a sequence $(\kappa_n^{2})_{n=1}^\infty$
of non-degenerated probability measures $\kappa_n^{2}$ on $C_n^{(2)}$
such that the infinite product measure $\kappa^2:=\bigotimes_{n=1}^\infty\kappa_n^{2}$ is nonatomic and
the associated flow of the tail equivalence relation on $(X_0^{(2)},\kappa^{2})$ is isomorphic to $V$.
Denote by $\kappa_n$ the convolution $\kappa_n^1*\kappa_n^2$.
Then $\kappa_n$ is a non-degenerated probability measure on $C_n$ for each $n\in\Bbb N$.
We now select inductively a sequence $(\nu_n)_{n=1}^\infty$ of measures on $G$ such that $\nu_n$ is supported on $F_n$ for each $n\in\Bbb N$ and
\thetag{1-3} holds.
Consider now the sequence  $\Cal T:=(C_n,F_{n-1},\kappa_n,\nu_{n-1})_{n=1}^\infty$.
Of course, \thetag{1-2} holds.
It is straightforward  to verify that Proposition~1.10(iii) holds for  $\Cal T$.
Let $\mu$ denote the  $(C,F)$-measure  determined by $(\kappa_n)_{n=1}^\infty$ and
$(\nu_{n})_{n=0}^\infty$.
It follows from Proposition~1.10 that 
$\mu$ is quasi-invariant under $T$.
Thus, (i) holds.

According to $(*)$, the associated flow of $T$ is isomorphic to the associated flow of the tail equivalence relation on $(X_0,\mu\restriction X_0)$.
Then $(**)$ yields that the later flow is isomorphic to the associated flow of the tail equivalence relation on $(X_0^{(2)},\kappa^{2})$.
Hence,  the associated flow of $T$ is isomorphic to $V$, i.e. (ii) is proven.

We deduce from~\thetag{5-1} that 
$$
(\mu\restriction X_0)\circ \pi^{-1}=(\kappa^1\circ\pi^{-1})*(\kappa^2\circ\pi^{-1}).
$$
It was shown in \S5.1 that $\kappa^1\circ\pi^{-1}=\nu$.
As $\nu$ is the Haar measure on the compact Abelian group $Y$, we obtain that 
$\nu*(\kappa^2\circ\pi^{-1})=\nu$.
Thus, $(\mu\restriction X_0)\circ \pi^{-1}=\nu$.
Hence, $\mu\circ\pi^{-1}\gg \nu$.
As the two measures, $\mu\circ\pi^{-1}$ and $\nu$ on $Y$ are quasi-invariant and ergodic under $T$, it follows that  $\mu\circ\pi^{-1}\sim \nu$.
Thus, (iii) is proven.

If $O$ is not the maximal odometer factor of $T$ then there is a prime number $p>2$ such that the homogeneous $\Bbb Z$-space 
$\Bbb Z/p\Bbb Z$ is a factor of $T$.
Since $\# C_n=4$ for each $n>0$, it can be deduced from Theorem~2.3 that each element of $C_n$ is divisible by $p$ eventually in $n$.
However, $4^n\in C_n$ for each $n$, a contradiction.
Thus, (iv) is proven.

As for (v), it is proved almost literally in the same way as in Example~5.1.

The second claim of the proposition follows directly from the first one.
\qed

\enddemo

\subhead 5.3. Rank-one $\Bbb Z^2$-action without odometer factors but whose  generators
have $\Bbb Z$-odometer factors
\endsubhead
In this section $G=\Bbb Z^2$.
Only finite measure preserving actions of $G$ are considered in this section.
In \cite{JoMc, \S6},  a rank-one $\Bbb Z^2$-action $T$ is constructed such that 
\roster
\item"(a)" each of the generators $T_{(0,1)}$ and
$T_{(0,1)}$ of $T$ has an odometer factor (as a $\Bbb Z$-action) but 
\item"(b)" $T$  has no $\Bbb Z^2$-odometer factors.
\endroster
The corresponding construction is rather involved (see \cite{JoMc, Theorem~6.1}).
We provide a different, elementary  example of $T$ possessing (a) and (b).
To prove that, we do not utilize any machinery developed in the previous sections.

\example{Example 5.4} Let $R$ be an irrational rotation on the circle $(\Bbb T,\lambda_\Bbb T)$
and let $S$ be an ergodic rotation on a compact totally disconnected Abelian infinite group  $Y$ endowed with the Haar measure $\lambda_Y$.
Then the dynamical system $O=(S^n)_{n\in\Bbb Z}$ is a classical odometer on $(Y,\lambda_Y)$.
We now define a $\Bbb Z^2$-action $T=(T_g)_{g\in\Bbb Z^2}$ on the product space $(\Bbb T\times 
 Y,\lambda_\Bbb T\otimes\lambda_Y)$ by setting
$
T_{(n,m)}:=R^n\times S^m.
$
Of course, $T$ is of rank-one along a F{\o}lner sequence of rectangles $(F_n)_{n=1}^\infty$ in $\Bbb Z^2$.
Let $\theta:=\pmatrix 1&0\\1&1\endpmatrix\in SL_2(\Bbb Z)$.
Then the action $\widetilde T:=(T_{\theta(g)})_{g\in\Bbb Z^2}$ is of rank-one along the sequence
$(\theta(F_n))_{n=1}^\infty$.
Of course, $(\theta(F_n))_{n=1}^\infty$ is also F{\o}lner.
The generators $\widetilde T_{(1,0)}=R\times S$ and $\widetilde T_{(0,1)}=I\times S$ of $\widetilde T$
 have $S$ as a factor.
 Thus, (a) holds for  $\widetilde T$.

It remains to show that 
$\widetilde T$ satisfies~(b), i.e.  that $\widetilde T$ has no   $\Bbb Z^2$-odometer factors.
Suppose that  $\widetilde T$ has an odometer factor  $\widetilde O$.
Let $K$ be the space of $\widetilde O$ and
let  $\pi:\Bbb T\times Y\to K$ denote the corresponding $\Bbb Z^2$-equivariant factor mapping.
Then $K$ is a totally disconnected compact Abelian group.
Since $\widetilde T$ has pure discrete spectrum, there is a compact subgroup $H\subset
\Bbb T\times Y$ such that $K=(\Bbb T\times Y)/H$ and  $\pi(g)=gH$ for each $g\in \Bbb T\times Y$.
Since $\Bbb T$ is connected and $K$ is totally disconnected, the closed subgroup $\pi(\Bbb T\times\{0\})\subset K$ is  trivial.
This means that $H$ contains the subgroup $\Bbb T\times\{0\}$.
Hence $K$ is  a quotient group of $Y=(\Bbb T\times Y)/(\Bbb T\times\{0\})$ indeed and $\pi$ is the corresponding projection map. 
It follows that $\widetilde O$ is not faithful:
 $$
 \widetilde O_{(n,-n)}=\pi\circ \widetilde T_{(n,-n)}=\pi\circ T_{(n,0)}=\pi\circ (S^n\times I)=I
 $$
  for each $n\in\Bbb Z$.
Since each odometer  action is faithful,   $\widetilde O$ is not  odometer. 
Thus,~(b) is proven.
\endexample

\comment

\example{Example 4.3} Let $R$ and $S$ be two irrational rotations on the circle $(\Bbb T,\lambda_\Bbb T)$ such that
the direct product $R\times S$ is ergodic.
Let $O$ stand for the 2-adic  odometer.
We consider $O$ as a rotation on 
the (compact totally disconnected) group of 2-adic integers $Y$
 furnished with the Haar measure $\nu$.
We now define a $\Bbb Z^2$-action $T=(T_g)_{g\in\Bbb Z^2}$ on the product space $(\Bbb T\times 
\Bbb T\times Y,\lambda_\Bbb T\otimes\lambda_\Bbb T\otimes\nu)$ by setting
$$
T_{(n,m)}(x_1,x_2,x_3):=(R^nx_1,S^mx_2,O^{n+m}x_3).
$$
Of course, $T$ is free and ergodic.
We claim that $T$ is of rank one.
Since $O$ is a 2-adic odometer, there is sequence Rokhlin towers $(A_k, F_k)_{k=1}^\infty$ for $O$
that generates the entire Borel  $\sigma$-algebra on $Y$,  where $F_k=\{0,1,\dots,2^k-1\}$ for each $k\in \Bbb N$.
By \cite{dJ1},
 the transformation $R\times S^{-1}$ is of rank one.
 Hence, there is a sequence of Rokhlin towers $(B_k, F_k')_{k=1}^\infty$ for $R\times S^{-1}$ such that  $F_k'=\{0,\dots,m_k-1\}$ for some $m_k>0$  
 and
 for each Borel subset $D\subset\Bbb T^2$,
$$
\max_{j\in F_k}\min_{I\subset F_k'}\lambda_\Bbb T\otimes\lambda_{\Bbb T}\Bigg(D\triangle \bigsqcup_{i\in I}(R^i\times S^{-i})(I\times S^j)B_k\Bigg)\to 0
$$
as $k\to\infty$.
It follows 
 that 
 $$
 \big(B_k\times A_k, \{(i,j-i)\mid j\in F_k, i\in F_k'\}\big)_{k=1}^\infty
 $$
  is a sequence of Rokhlin towers for $T$ that  generates the entire Borel  $\sigma$-algebra on $\Bbb T^2\times Y$.
Of course, $\big(\{(i,j-i)\mid j\in F_k, i\in F_k'\}\big)_{k=1}^\infty$ is a F{\o}lner sequence in $\Bbb Z^2$.
Hence $T$ is of rank one. 
It is obvious that $O$ is an odometer factor of transformations $T_{(1,0)}$ and $T_{(0,1)}$. 

It remains to show that $T$ has no $\Bbb Z^2$-odometer factors.
Suppose that  $T$ an odometer factor  $\widetilde O$.
Let $K$ be the space of $\widetilde O$ and
let  $\pi:\Bbb T^2\times Y\to K$ denote the corresponding $\Bbb Z^2$-equivariant factor mapping.
Then $K$ is a totally disconnected compact Abelain group.
Since $T$ is a minimal rotation on  $\Bbb T^2
\times Y$, there is a compact subgroup $H\subset
\Bbb T^2\times Y$ such that $K=(\Bbb T^2\times Y)/H$ and  $\pi(g)=gH$ for each $g\in \Bbb T^2\times Y$.
Since $\Bbb T^2$ is connected and $K$ is totally disconnected, the closed subgroup $\pi(\Bbb T^2\times\{0\})\subset K$ is  trivial.
This means that $H$ contains $\Bbb T^2\times\{0\}$.
Hence $K$ is  a quotient group of $Y=(\Bbb T^2\times Y)/(\Bbb T^2\times\{0\})$ indeed and $\pi$ is the corresponding projection map. 
It follows that $\widetilde O$ is not free:
 $\widetilde O_{(n,-n)}=\pi\circ T_{(n,-n)}=I$ for each $n\in\Bbb Z$.
 However, each  odometer $\Bbb Z^2$-action  is free because every odometer of an Abelian group is normal.
\endexample

\endcomment

According to the terminology of \cite{JoMc} (which is different from ours),  $\Bbb Z^2$-odometers can be non-free.
Hence, (b)
can be interpreted as  ``every rank-one factor of $T$ is non-free''.

\remark{Remark 5.5} We note that  in  the example \cite{JoMc, \S6}, the generators $T_{(0,1)}$ and $T_{(1,0)}$ are non-ergodic. 
In Example~5.4,  $\widetilde T_{(1,0)}$ is ergodic but $\widetilde T_{(0,1)}$ is not.
However, if we change $\theta$ with the matrix $\theta':=\pmatrix 2&1\\1&1\endpmatrix$ then we obtain a new example of rank-one $\Bbb Z^2$-action  $\widetilde T$ possessing (a), (b) and
\roster
\item"(c)" each of the generators $\widetilde T_{(0,1)}$ and
$\widetilde T_{(0,1)}$ of $\widetilde T$ is ergodic.
\endroster
\endremark

\subhead 5.4. Heisenberg group actions of rank one
\endsubhead
In this section we consider  the 3-dimensional discrete Heisenberg group $H_3(\Bbb Z)$.
We remind that 
$$
H_3(\Bbb Z)=
\left\{
\pmatrix 
1 & x & z\\
0&1&y\\
0&0&1
\endpmatrix\Bigg|\, x,y,z\in\Bbb Z\right\}.
$$
This group is non-Abelian, nilpotent and residually finite.
For brevity, we will write $(x,y,z)$ for the matrix $\pmatrix 
1 & x & z\\
0&1&y\\
0&0&1
\endpmatrix$.
It is straightforward to verify that 
$$
\gather
(x,y,z)\cdot(x_1,y_1,z_1)=(x+x_1,y+y_1, z+z_1+xy_1).\\
(x,y,z)\cdot(x_1,y_1,z_1)\cdot(x,y,z)^{-1}=(x_1,y_1, z_1+xy_1-yx_1).
\endgather
$$
The center of $H_3(\Bbb Z)$ is $\{(0,0,z)\mid z\in\Bbb Z\}$.
Given $a,b,c\ge 0$, we let
$$
\Pi(a,b,c):=\{(x,y,z)\in H_3(\Bbb Z)\mid 0\le x<a,0\le y<b,0\le z<c\}.
$$
It is straightforward to verify that if $a_n\to+\infty$, $b_n\to+\infty$, $c_n\to+\infty$
and $\frac{b_n}{c_n}\to 0$ as $n\to\infty$ then
$(\Pi(a_n,b_n,c_n))_{n=1}^\infty$ is a left F{\o}lner sequence in $H_3(\Bbb Z)$.
If, in addition, $\frac{a_n}{c_n}\to 0$ then
$(\Pi(a_n,b_n,c_n))_{n=1}^\infty$ is a 2-sided F{\o}lner sequence in $H_3(\Bbb Z)$.

\example{Example 5.6} We  construct a probability preserving $(C,F)$-action $T$ of $H_3(\Bbb Z)$
that has an odometer factor but $T$ itself is not isomorphic to any odometer.
For that, we first define recurrently a sequence  $(h_n)_{n=0}^\infty$ of positive integers 
 by setting 
$$
h_0:=1\text{ and $h_{n+1}:=16h_n+9\cdot 4^{n+1}$.}
$$
It is easy to check that  
$h_{n}=4^{2n+1}- 3\cdot 4^n$ for each $n\ge 0$.
We set
$$
\align
C_{n+1}^{(1)}&:=\big\{(a,b,0)\mid a,b\in\{0, 2^n\} \big\}\text{ and}\\
C_{n+1}^{(2)}&:=\{(0,0,jh_n)\mid j=0,\dots,7\}\cdot\{(0,0,0), (0,0, 8h_n+2\cdot 4^{n+1})\}.
\endalign
$$
We now define a sequence $(C_n,F_{n-1})_{n=1}^\infty$ by setting
$$
F_{n}:=\Pi(2^n,2^n,h_n)\text{ \ and \ }C_{n+1}:=C_{n+1}^{(1)}C_{n+1}^{(2)}\qquad\text{for each $n\ge 0$}.
$$
It is straightforward to check that \thetag{1-1} is satisfied.
We define measures $\kappa_n$ on $C_n$ and $\nu_n$ on $F_n$ by setting $\nu_0(1)=1$ and 
$$
\kappa_n(c):=\frac1{\# C_n}=\frac 1{64},\quad\nu_n(f):=\frac1{64^n}\qquad\text{for each $c\in C_n$, $f\in F_n$, $n>0$.}
$$
Let $\Cal T:=(C_n,F_{n-1},\kappa_n,\nu_{n-1})$.
Then \thetag{1-2} and \thetag{1-3} hold for $\Cal T$.
Since $\# F_n=4^n(4^{2n+1}-3\cdot 4^n)$ and $\# C_{n+1}=64$, we obtain that
$$
\prod_{n>0}\frac{\# F_{n+1}}{\# F_n\#C_{n+1}}=\prod_{n>0}\frac{4^{n+2}-3}{4(4^{n+1}-3)}
=\prod_{n>0}\Bigg(1+\frac9{4(4^{n+1}-3)}\Bigg)<\infty.\tag5-3
$$
Of course, $(F_n)_{n=1}^\infty$ is a 2-sided F{\o}lner sequence in $H_3(\Bbb Z)$.
Hence,
Proposition~1.10(ii) holds.
Then the $(C,F)$-action $T=(T_g)_{g\in H_3(\Bbb Z)}$ 
associated with $\Cal T$
is well defined on a measure space $(X,\mu)$, where $\mu$ is the $(C,F)$-measure determined by
$(\kappa_n,\nu_{n-1})_{n=1}^\infty$.
Moreover, $T$ preserves $\mu$, i.e. $\mu$ is the Haar measure for the $(C,F)$-equivalence relation on $X$, and $\mu(X)<\infty$ in view of~\thetag{5-3} (see Remark~1.5).

Next, we define a measure preserving odometer action of $H_3(\Bbb Z)$.
Given $n>0$, we let 
$$
\Gamma_n:=\{(i\cdot 2^n,j\cdot 2^n,k\cdot 2^n)\in H_3(\Bbb Z)\mid i,j,k\in\Bbb Z\}.
$$
Then $
\Gamma_n$ is a normal cofinite subgroup of $H_3(\Bbb Z)$,
$\Gamma_1\supsetneq\Gamma_2\supsetneq\cdots$ and $\bigcap_{n=1}^\infty\Gamma_n=\{1\}$.
Denote by $O$  the $H_3(\Bbb Z)$-odometer
associated with the sequence $(\Gamma_n)_{n=1}^\infty$.
We call it {\it the 2-adic odometer} action of $H_3(\Bbb Z)$.
This odometer is normal.
It is defined on the compact metric group
 $$
          Y:=\projlim_{n\to\infty} H_3(\Bbb Z)/\Gamma_n.
 $$ 
Denote by $\nu$ the Haar measure on $Y$.
Since $C_{n+1}\subset\Gamma_n$ for each $n\in\Bbb N$, it follows that $\Cal T$ is compatible with $(\Gamma_n)_{n=1}^\infty\in Y$.
Hence the $(\Cal T, (\Gamma_n)_{n=1}^\infty)$-factor mapping $\pi_{(\Cal T,(\Gamma_n)_{n=1}^\infty)}:X\to Y$ intertwines $T$ with $O$.
For brevity, instead of $\pi_{(\Cal T,(\Gamma_n)_{n=1}^\infty)}$ below we will write $\pi$.
The measure $\mu\circ \pi^{-1}$ on $Y$
is finite and invariant under $O$.
Since $O$ is uniquely ergodic, it follows that $\mu\circ \pi^{-1}=\mu(X)\cdot\nu$.
Thus, the 2-adic $H_3(\Bbb Z)$-odometer is a finite measure preserving factor of  $(X,\mu, T)$.

We claim that  $(Y,\nu, O)$ is the {\it maximal odometer factor} of $T$, i.e. every odometer factor of $T$ is a factor of $O$.
Let $\Lambda_1\supset\Lambda_2\supset\cdots$ be a sequence of cofinite subgroups in $H_3(\Bbb Z)$ with $\bigcap_{n=1}^\infty\bigcap_{g\in G}g\Lambda_ng^{-1}=\{1\}$ and let $Q=(Q_g)_{g\in H_3(\Bbb Z)}$ stand for the associated $H_3(\Bbb Z)$-odometer which is
 defined on the space $Z:=\projlim_{n\to\infty}H_3(\Bbb Z)/\Lambda_n$ equipped with Haar measure.
If there is an $H_3(\Bbb Z)$-equivariant mapping $\tau:X\to Z$ then,
by Theorem~4.4, there exist a sequence $\boldsymbol q=(q_n)_{n=0}^\infty$ and an element
$z=(g_n\Lambda_n)_{n=1}^\infty\in Z$ 
such that $0=q_0<q_1<q_2<\cdots$,  the $\boldsymbol q$-telescoping $\widetilde {\Cal T}$ of $\Cal T$ is compatible with $z$
and $\tau=\pi_{(\widetilde{\Cal T},z)}\circ\iota_{\boldsymbol q}$.
Replacing $\Lambda_n$ with $g_n^{-1}\Lambda_ng_n$ for each $n\in\Bbb N$, we pass to 
an isomorphic (to $Q$) odometer  as a factor of $X$ (see the proof of Corollary~4.6).
We denote it by the same symbol $Z$.
Therefore, without loss of generality, we may assume that 
$\widetilde {\Cal T}$ is compatible with 
$z=(\Lambda_n)_{n=1}^\infty$.
It follows that there is $N>0$ such that for each $n>N$
$$
\#\{c\in C_{q_n+1}\cdots C_{q_{n+1}}\mid c\not\in \Lambda_n\}<64^{-1}.
$$ 
As $\# C_{q_n+1}=64$ and 
$\#(C_{q_n+1}\cdots C_{q_{n+1}})=\#(C_{q_n+1})\#(C_{q_n+2}\cdots C_{q_{n+1}})$,  a version of Fubuni theorem yields that
there exists at least  one
element $d\in C_{q_n+2}\cdots C_{q_{n+1}}$ such that  $C_{q_n+1}d\subset \Lambda_n$.
Hence, 
$$
\{(2^{q_n},0,0), (0,2^{q_n},0)\}\subset\{\widetilde cc^{-1}\mid \widetilde c,c\in C_{q_n+1} \}\subset \Lambda_n.
$$
Therefore, the subgroup 
$$
\Sigma_n:=\{(i\cdot 2^{q_n},j\cdot 2^{q_n},k\cdot 4^{q_n})\in H_3(\Bbb Z)\mid i,j,k\in\Bbb Z\}
$$ 
is contained in $\Lambda_n$.
This implies that $\Gamma_{2q_n}\subset\Sigma_n\subset \Lambda_n$ for each $n>0$.
The  natural projection 
$$
H_3(\Bbb Z)/\Gamma_{2q_n}\to H_3(\Bbb Z)/\Lambda_n
$$
 is $H_3(\Bbb Z)$-equivariant for each $n$.
Passing to the projective limit as $n\to\infty$, we obtain a $H_3(\Bbb Z)$-equivariant projection $\eta:Y\to Z$.
It is straightforward to verify that  $\tau=\eta\circ\pi$, as desired.

It remains to show that $\pi$ is not one-to-one (mod\,0).
For that, it suffices to show that the restriction of $\pi$ to $X_0$ is not one-to-one.
Let  $\kappa_n^1$ and $\kappa_n^2$ be the equidistributions on  $C^{(1)}_n$ and $C^{(2)}_n$
respectively, $n\in\Bbb N$.
We let
$$
\align
(X_0^{(1)}, \kappa^1)&:=\bigg(C_1^{(1)}\times C_2^{(1)}\times\cdots,\bigotimes_{n=1}^\infty\kappa^1_n\bigg)\quad\text{and}\\
(X_0^{(2)}, \kappa^2)&:=\bigg(C_1^{(2)}\times C_2^{(2)}\times\cdots,\bigotimes_{n=1}^\infty\kappa^2_n\bigg).
\endalign
$$
Then $X_0, X_0^{(1)}$ and $X_0^{(2)}$ are compact subsets of the Polish group $H_3(\Bbb Z)^\Bbb N$.
The mapping 
$$
\alpha:X_0^{(1)}\times X_0^{(2)}\ni (x^1,x^2)\mapsto x^1x^2\in X_0
$$
 is a homeomorphism that maps the product measure $\kappa^1\otimes\kappa^2$ to $\mu\restriction X_0$.
  Moreover,  for each $(x^1,x^2)\in X_0^{(1)}\times X_0^{(2)}$,
$$
\pi(x^1x^2)=\pi(x^1)\pi(x^2).\tag5-4
$$
Let $\Cal Z$ stand for the center of $H_3(\Bbb Z)$.
Denote by $Y_0$ the center of $Y$.
It is routine to verify that
$$
Y_0=\projlim_{n\to\infty}\Cal Z\Gamma_n/\Gamma_n =\projlim_{n\to\infty}
\Cal Z/(\Cal Z\cap \Gamma_n).
$$
Denote by $Y_2$ the quotient group $Y/Y_0$ and by $\omega$ the quotient homomorphism $Y\to Y_2$.
 Then we obtain a short exact sequence of compact totally disconnected groups
$$
1\longrightarrow Y_0\longrightarrow Y\overset{\omega}\to\longrightarrow Y_2\longrightarrow 1.
$$
Denote by $\lambda_0$ and $\lambda_2$ the Haar measures on $Y_0$ and $Y_2$ respectively.
It is straightforward to verify that 
$$
\gather
Y_2=\projlim_{n\to\infty}H_3(\Bbb Z)/(\Cal Z\Gamma_n)=\projlim_{n\to\infty}\Bbb Z^2/2^n\Bbb Z^2\quad\text{and}\\
\omega\circ \pi(x^1)=\big((a_1,b_1)+2\Bbb Z^2,(a_2,b_2)+2^2\Bbb Z^2,\dots \big)
\endgather
$$
for each $x^1=\big((a_1,b_1,0),(a_2,b_2,0),\dots\big)\in X^{(1)}_0$.
Hence $\omega\circ \pi$ is a measure preserving homeomorphism of $(X_0^{(1)},\kappa^{1})$ onto $(Y_2,\lambda_2)$.
We now define a continuous  mapping $s:Y_2\to Y$ by setting
$$
s(\omega(\pi(x^1))):=\pi(x^1).\tag5-5
$$
Then $s$ is a cross-section of $\omega$.
Hence, the mapping
$$
\beta:Y\ni y\mapsto (ys(\omega(y))^{-1},\omega(y))\in Y_0\times Y_2
$$
is a well defined measure preserving homeomorphism of  $(Y,\nu)$ onto the product measure space $\big(Y_0\times Y_2,\lambda_0\otimes\lambda_2\big)$.
It follows from \thetag{5-4} and \thetag{5-5} that
$$
\beta\circ \pi\circ \alpha(x^1,x^2)=
\omega(\pi(x^1)\pi(x^2))=
\big(\pi(x^2),\omega\circ \pi(x^1)\big)
$$
for each $(x^1,x^2)\in X_0^{(1)}\times X_0^{(2)}$.
Moreover, 
$$
(\kappa^1\otimes\kappa^2)\circ(\beta\circ \pi\circ \alpha)^{-1}=\lambda_0\otimes\lambda_2.
$$
Therefore, $\pi$ is not one-to-one ($\mu$-mod\,0) if and only if
the mapping $\pi\restriction X^{(2)}_0\to Y_0$ is 
 not one-to-one ($\mu^{(2)}$-mod\,0).
 
Our purpose now is to show that $\pi\restriction X^{(2)}_0\to Y_0$ is 
 not one-to-one.
Let
 $$
 \align
 C^{(3)}_{n+1}&:=\{0, h_n,2h_n,3h_n\},\\
  C^{(4)}_{n+1}&:=\{0,4h_n\}+\{0,8h_n + 2 á\cdot4
^{n+1}\},
\endalign
$$
$\kappa_n^3$ and $\kappa_n^4$ are the equidistributions  on $C^{(3)}_{n}$ and $C^{(4)}_{n}$
respectively
 for each $n\ge 0$ and let
$$
\align
(X_0^{(3)}, \kappa^3)&:=\bigg(C_1^{(3)}\times C_2^{(3)}\times\cdots,\bigotimes_{n=1}^\infty\kappa_n^3\bigg)\quad\text{and}\\
(X_0^{(4)}, \kappa^4)&:=\bigg(C_1^{(4)}\times C_2^{(4)}\times\cdots,\bigotimes_{n=1}^\infty\kappa_n^4\bigg).
\endalign
$$
Then $X_0^{(3)}$ and $X_0^{(4)}$ are compact subsets of  $H_3(\Bbb Z)^\Bbb N$.
The mapping 
$$
X_0^{(3)}\times X_0^{(4)}\ni (x^3,x^4)\mapsto x^3x^4\in X_0^{(2)}
$$
 is a well defined homeomorphism that maps the product measure $\kappa^3\otimes\kappa^4$ to $\kappa^2$.
  Moreover,  for each $(x^3,x^4)\in X_0^{(3)}\times X_0^{(4)}$,
$$
\pi(x^3x^4)=\pi(x^3)\pi(x^4).
$$
In view of that, it suffices to show that the mapping $\pi\restriction X^{(3)}\to Y_0$ is a
a bijection and $\kappa^3\circ(\pi\restriction X^{(3)})^{-1}=\lambda_0$.
We leave a routine verification of these facts to the reader.

\comment

For each $n>0$, the set $\Cal Z\Gamma_n/\Gamma_n$ is a subgroup of $H_3(\Bbb Z)/\Gamma_n$.
We obtain a projective system $\Cal Z\Gamma_1/\Gamma_1\leftarrow\Cal Z\Gamma_2/\Gamma_2 \leftarrow\cdots$, which is a subsystem of the system
$H_3(\Bbb Z)/\Gamma_1\leftarrow H_3(\Bbb Z)/\Gamma_2\leftarrow\cdots$.
Since $\Cal Z\Gamma_n/\Gamma_n$ is identified naturally with $\Cal Z/(\Cal Z\cap \Gamma_n)$ for each $n$,
we obtain a commutative diagram
$$
\CD
H_3(\Bbb Z)/\Gamma_1 @<<< H_3(\Bbb Z)/\Gamma_2 @<<<\cdots\\
@AAA @AAA\cdots\\
\Cal Z/(\Cal Z\cap \Gamma_1)@<<<
\Cal Z/(\Cal Z\cap \Gamma_2)@<<<
\cdots
\endCD
\tag4-4
$$
We note that all the arrows in the diagram are $\Cal Z$-equivariant mappings.
Denote by $(Y_0,S)$ the topological $\Cal Z$-odometer associated with the nested sequence
$(\Cal Z\cap \Gamma_n)_{n=1}^\infty$ of subgroups in $\Cal Z$.
It follows from $\thetag{4-4}$ that $Y_0$ is a closed normal subgroup of $Y$ and  
$O_g\restriction Y_0=S_g$ for each $g\in\Cal Z$.
Denote by $Y_2$ the quotient group $Y/Y_0$ and by $\omega$ the quotient homomorphism $Y\to Y_2$.
It is straightforward to verify that 
$$
Y_2=\projlim_{n\to\infty}H_3(\Bbb Z)/(\Cal Z\Gamma_n)=\projlim_{n\to\infty}\Bbb Z^2/(2^n\Bbb Z)^2.
$$
Denote by $\nu_0$ and $\nu_2$ the Haar measures on $Y_0$ and $Y_2$ respectively.
Then there is a Borel isomorphism of $(Y,\nu)$ onto $(Y_0\times Y_2,\nu_0\otimes\nu_2)$.
It is straightforward to verify that  $\omega\circ\pi_{(\Cal T,(\Gamma_n)_{n=1}^\infty)}$ maps  $X^{(1)}$ bijectively
onto $Y_2$.
Of course, $\pi_{(\Cal T,(\Gamma_n)_{n=1}^\infty)}(x)\in Y_0$ for each $x\in X_0^{(2)}$.

\endcomment

\endexample

With the example below we illustrate that the common concept of normality for  odometers 
is not invariant under isomorphism.
The normality depends on the choice of the underlying sequence of $(\Gamma_n)_{n=1}^\infty$
of cofinite subgroups in $G$.

\example{Example  5.7}
Let $\Gamma_n:=\{(i2^{n-1},j2^n,k2^n)\in H_3(\Bbb Z)\mid i,j,k\in\Bbb Z\}$.
Of course, $\Gamma_n$ is a non-normal cofinite subgroup of $H_3(\Bbb Z)$ and $\Gamma_1\supsetneq\Gamma_2\supsetneq\cdots$ with $\bigcap_{n=1}^\infty\Gamma_n=\{1\}$.
It is easy to see that the largest normal subgroup $\widetilde \Gamma_n$ of $\Gamma_n$ is
$$
\widetilde\Gamma_n:=\{(i2^{n},j2^n,k2^n)\in H_3(\Bbb Z)\mid i,j,k\in\Bbb Z\}.
$$
Thus, $\widetilde\Gamma_n$ is of index 2 in $\Gamma_n$ for each $n>0$.
Denote by $(Y,O)$ and $(\widetilde Y,\widetilde O)$ the topological odometers associated with
 $(\Gamma_n)_{n=1}^\infty$ and $(\widetilde\Gamma_n)_{n=1}^\infty$ respectively.
 Then $(\widetilde Y,\widetilde O)$ is a normal odometer while $(Y,O)$ is not.
Moreover, $(\widetilde Y,\widetilde O)$ is the topological normal cover of $(Y,O)$.
 As was shown in \S3.2, $(Y,O)$ is a factor of $(\widetilde Y,\widetilde O)$ under the natural projection 
 $\omega:\widetilde Y\ni\widetilde y\mapsto\widetilde y H\in \widetilde Y/H=Y$, where 
 $$
 H:=\{(\widetilde y_n)_{n=1}^\infty\in\widetilde Y\mid \widetilde y_n\in\Gamma_n/\widetilde \Gamma_n\text{ for all }n\in\Bbb N\}
 $$
 is a closed subgroup of $\widetilde Y$.
Since $\Gamma_{n+1}\subset\widetilde\Gamma_n\subset\Gamma_n$ for each $n$, it follows that
$\omega$ is one-to-one.
Hence $\omega$ is an isomorphism of $\widetilde O$ with $O$, and
 $H=\{1\}$.
\endexample

\comment

\remark{Remark 5.8}
We encounter with a similar phenomenon in Examples~3.5 and 3.6.
For instance, let $O$ stand for the odometer defined in Example~3.5.
Then $O$ is not normal but the canonical factor mapping $\omega$ from the normal cover $\widetilde O$ to $O$ is one-to-one. 
However, in this case we do not have the ``alternating'' property $\Gamma_{n+1}\subset\widetilde \Gamma_n\subset\Gamma_n$, which was crucial for the invertibility of $\omega$ in Example~5.7.
In Example~3.5 we have that $\Gamma_{m}\not\subset\widetilde \Gamma_n$ for any $m\ge n$.
\endremark

\endcomment

\comment

Let $k_1:=1$ and $k_n:=k_{n-1}(k_{n-1}+1)$ for $n>1$.
We now set
$$
\Gamma_n:=\{c(k_n^3i_3)b(k_n^2i_2)a(k_ni_1)\mid i_1,i_2,i_3\in\Bbb Z\}.
$$
Then $\Gamma_n$ is a cofinite subgroup in $H_3(\Bbb Z)$,   $\Gamma_1\supsetneq\Gamma_2\supsetneq\cdots$ and $\bigcap_{n=1}^\infty\Gamma_n=\{1\}$.

\endcomment
\comment

Let $A$ stand for the multiplicative group $\{-1,1\}$.
We now consider $\Bbb Z^2$ as an $A$-module by setting $(-1)\cdot (l,m):=(m,l)$ for each $(l,m)\in\Bbb Z^2$.
Then the semidirect product  $G:=\Bbb Z^2\rtimes A$ is well defined.
For each $n>0$, let $\Gamma_n:=2^n\Bbb Z$.
Then $\Gamma_n$ is a cofinite subgroup of $\Bbb Z$ and $\Gamma_1\supsetneq\Gamma_2\supsetneq\cdots$ with $\bigcap_{n>0}\Gamma_n=\{0\}$.
We now let 
$$
\Sigma_n:=(\Gamma_n\oplus\Gamma_{n+1})\times \{1\}\subset G\quad\text{and}\quad
\Sigma_n':=(\Gamma_{n+1}\oplus\Gamma_{n})\times \{1\}\subset G.
$$
Then $\Sigma_n$ and $\Sigma_n'$ are cofinite subgroups of $G$,
$$
\Sigma_1\supsetneq\Sigma_2\supsetneq\cdots,\quad 
\Sigma_1'\supsetneq\Sigma_2'\supsetneq\cdots\quad\text{and}\quad
\bigcap_{n>0}\Sigma_n=\bigcap_{n>0}\Sigma_n'=\{1_G\}.
$$
Let $\Cal G:=\bigoplus_{j\in\Bbb N} G$,
$$
\Cal G_n:=\underbrace{\Sigma_n\oplus\cdots\oplus\Sigma_n}_{n\text{ times}}\oplus G\oplus G\oplus\cdots\quad\text{and}\quad
\Cal G_n':=\underbrace{\Sigma_n'\oplus\cdots\oplus\Sigma_n'}_{n\text{ times}}\oplus G\oplus G\oplus\cdots.
$$
Then $\Cal G_n$ and $\Cal G_n'$ are  cofinite subgroup of $\Cal G$,
$$
\Cal G_1\supsetneq\Cal G_2\supsetneq\cdots,\quad \Cal G_1'\supsetneq\Cal G_2'\supsetneq\cdots\quad\text{and}\quad\bigcap_{n>0}\Cal G_n=\bigcap_{n>0}\Cal G_n'=\{1_{\Cal G}\}.
$$
We now let
$$
g_n:=\big(\underbrace{(0,0,-1),(0,0,-1),\dots,(0,0,-1)}_{n\text{ times}},1_G,1_G,\dots\big)\in\Cal G.
$$
It is straightforward to verify that  $g_n \Cal G_ng_n^{-1}=\Cal G_n'$
for each $n\in\Bbb N$.

Denote by $(Y,O)$  and $(Y',O')$ the topological $\Cal G$-odometers associated with the sequences $(\Cal G_n)_{n=1}^\infty$ and $(\Cal G_n')_{n=1}^\infty$ respectively.
Furnish $Y$ and $Y'$ with the Haar measures $\nu$ and $\nu'$ respectively.
We now show that the probability preserving $\Cal G$-odometers 
$(Y,O,\nu)$  and $(Y',O',\nu')$ are of rank one.
For each $n\in\Bbb N$, we let
$$
\align
E_n&:=\{0,2^{n-1}\}\times\{0,2^n\}\times\{1\}\subset \Sigma_{n-1},\\
R_n&:=\{0,1,\dots,2^n-1\}
\times\{0,1,\dots,2^{n+1}-1\}\times A\subset G\quad\text{and}\\
D_n&:=\underbrace{E_n\times\cdots\times E_n}_{n-1\text{ times}}\times \,R_n\times\{1_G\}\times\{1_G\}\times\cdots\subset\Cal G_{n-1}.
\endalign
$$
It is straightforward to verify that  $E_n$ is a $\Sigma_n$-cross-section in $\Sigma_{n-1}$,
$R_n$ is a $\Sigma_n$-cross-section in $G$,
and
$D_n$ is a $\Cal G_{n}$-cross-section in $\Cal G_{n-1}$.
Let $\kappa_n$ stand for the equidistribution on $D_n$.
It is straightforward to check that (i)--(iii) of Proposition~3.2 hold.
Hence, by Proposition~3.2, $(Y,O,\nu)$ is of rank one.
In a similar way, $(Y',O',\nu')$ is of rank one.
Moreover, $(Y,O,\nu)$ is isomorphic to a $(C,F)$-action for which the corresponding sequence $(C_n)_{n=1}^\infty$ in the underlying parameters equals  $(D_n)_{n=1}^\infty$.

\endcomment
\comment
We now show that the probability preserving $H_3(\Bbb Z)$-odometer $(Y,O,\nu)$ is of rank one.
For each $n\in\Bbb N$, we let
$$
D_n:=\{c(k_{n-1}^3j_3)b(k_{n-1}^2)j_2a(k_{n-1}j_1)\mid 0\le j_3<k_n^3,0\le j_2<k_n^2,0\le j_1<k_n\}.
$$
It is straightforward to verify that $D_n$ is a $\Gamma_{n}$-cross-section in $\Gamma_{n-1}$.
Let $\kappa_n$ stand for the equidistribution on $D_n$.
It is straightforward to check that (i)--(iii) of Proposition~3.2 hold.
Hence, by Proposition~3.2, $(Y,O,\nu)$ is of rank one.

We now set
$$
\Gamma_n':=b(k_n)\Gamma_nb(k_n)^{-1}=\{c\big(k_n^3i_3-k_n^2i_1\big)b(k_n^2i_2)a(k_ni_1)\mid i_1,i_2,i_3\in\Bbb Z\}.
$$
As $b(k_n)b(k_{n-1})^{-1}=b(k_{n-1}^2)\in\Gamma_{n-1}$, we obtain that
$\Gamma_1'\supsetneq\Gamma_2'\supsetneq\cdots$.
Denote by $(Y',O')$ the topological odometer associated with $(\Gamma_n')_{n=1}^\infty$.
Furnish $Y'$ with the Haar measure $\nu'$.

\endcomment
\comment

\proclaim{Claim E}  $(Y',\nu',O')$ is not a factor of $(Y,\nu,O)$.
\endproclaim
\demo{Proof}
Indeed, suppose that the contrary holds.
Since $(Y,\nu,O)$ is of rank one, it follows from Theorem~3.13 that there is an element
$(h_n\Cal G_n')_{n=1}^\infty\in Y'$ and an increasing sequence $(q_n)_{n=1}^\infty$ of positive integers such that $q_n> n$ and
$$
\frac{\#\big(\{d\in D_{q_n}\cdots D_{q_{n+1}-1}\mid d\in h_n\Cal G_n'h_n^{-1}\}\big)}
{\#(D_{q_n}\cdots D_{q_{n+1}-1})}\to 1\qquad\text{as $n\to\infty$.}
$$
This implies that $\Cal G_{q_n-1}\subset h_n\Cal G_n'h_n^{-1}$ eventually in $n$.
The inequality yields, in turn, that $\Cal G_{n}\subset h_n\Cal G_n'h_n^{-1}$ and there is $d_n\in \bigoplus_{j\in\Bbb Z}(\Bbb Z^2\times\{1\})\subset\Cal G$
such that $h_n=d_ng_n$ eventually.
It follows that $h_{n}h_{n+1}^{-1}\not\in\Cal G_n'$.
Hence, $(h_n\Cal G_n')_{n=1}^\infty\not\in Y'$, a contradiction.
\enddemo

\endcomment
\comment

\proclaim{Claim E} Let $\Gamma_{q_n}\subset a(\alpha)b(\beta)\Gamma_nb(\beta)^{-1}a(\alpha)^{-1}$ for some $n>0$ and $\alpha,\beta\in\Bbb Z$.
Then 
$
\Gamma_n\subset a(\alpha)b(\beta)\Gamma_nb(\beta)^{-1}a(\alpha)^{-1}
$
and $k_n^2|\beta$ and $k_n|\alpha$.
\endproclaim
\demo{Proof} Since 
$$
a(\alpha)b(\beta)\Gamma_nb(\beta)^{-1}a(\alpha)^{-1}=\{c(k_n^3i_3+\alpha k_n^2i_2-\beta k_ni_1)b(k_n^2i_2)a(k_ni_1)\mid i_1,i_2,i_3\in\Bbb Z\}
$$
and $\Gamma_{q_n}$ is generated by  the two elements $a(k_{q_n})$ and $b(k_{q_n}^2)$,
it follows that $\Gamma_{q_n}$ is a subgroup of $a(\alpha)b(\beta)\Gamma_nb(\beta)^{-1}a(\alpha)^{-1}$ if and only if there exist integers $i_3$ and $i_3'$ such that
$k_n^3i_3-\beta k_{q_n}=0$ and $k_n^3i_3'-\alpha k_{q_n}^2=0$, i.e. $k_n^3|\beta k_{q_n}$ and 
$k_n^3|\alpha k_{q_n}^2$.
If $l>0$ then $k_{n+l}$ depends on $k_n$ in a polynomial way such that $k_{n+l}=k_n(k_n^{l}+\cdots+1)$.
It follows that $k_n|k_{n+l}$ and $k_n$ is relatively prime with the ratio $k_{n+l}/k_n$.
Therefore, $k_n^3|\beta k_{q_n}$ and 
$k_n^3|\alpha k_{q_n}^2$ are equivalent to $k_n^2|\beta $ and 
$k_n|\alpha$ respectively.
Thus, we obtain that
$$
\Gamma_{q_n}\subset a(\alpha)b(\beta)\Gamma_nb(\beta)^{-1}a(\alpha)^{-1}\iff k_n^2|\beta \text{ and }k_n|\alpha.
$$
Since the righthand side  does not depend on $q_n$ (but the inequality $q_n\ge n$ is our standing assumption), it follows that 
$
\Gamma_n\subset a(\alpha)b(\beta)\Gamma_nb(\beta)^{-1}a(\alpha)^{-1}
$, as desired.
\enddemo

For each $n\in\Bbb N$, select $\alpha_n,\beta_n,\gamma_n\in\Bbb Z$ such that $g_n=a(\alpha_n)b(\beta_n)c(\gamma_n)$.
It follows from Claim~E that  
$$
\Gamma_n\subset a(\alpha_n)b(\beta_n)\Gamma_nb(\beta_n)^{-1}a(\alpha_n)^{-1}\quad\text{and $k_n|\alpha_n$ and $k_n^2|(\beta_n+k_n)$.}
$$
Since $g_n$ is determined up to the right multiplication with $\Gamma_n'$, we may assume without loss of generality that $g_n$ belongs to the  ``natural'' $\Gamma_n$-fundamental subset in $H_3(\Bbb Z)$, i.e.
$$
g_n\in b(k_n)\big\{c(p)b(q)a(s)\mid 0\le p<k_n,0\le q<k_n^2,0\le s<k_n^3\big\}b(k_n)^{-1}.
$$
Then $k_n|\alpha_n$ and $k_n^2|(\beta_n+k_n)$ are only possible if $\alpha_n=0$ and $\beta_n=-k_n$.
This yields that  $g_n=b(-k_n)c(\gamma_n)$.
Hence, $g_{n+1}g_n^{-1}=b(k_n-k_{n+1})c(\gamma_{n+1}-\gamma_n)\not\in\Gamma_n'$ because
$k_{n+1}-k_{n}=k_n^2$

\endcomment
\comment

Since $(g_n\Gamma_n')_{n=1}^\infty\in Y'$, there are $m_n,m_n',m_n''\in\Bbb Z$
such that  $g_n=c(k_{n-1}^3m_n)b(k_{n-1}^2m_n')a(k_{n-1}m_n'')$.
Thus,
$$
\Gamma_{q_n}\subset
\{c\big(k_n^3i_3-k_n^2i_1+k_n^2i_2k_{n-1}m_n''- k_ni_1k_{n-1}^2m_n'\big)b(k_n^2i_2)a(k_ni_1)\mid i_1,i_2,i_3\in\Bbb Z\}.
$$
\endcomment

\head 6. Comments on the paper \cite{JoMc}
\endhead

The article \cite{JoMc} by A.S.A.~Johnson and D.M.~McClendon is also devoted to generalization of  
\cite{Fo--We}.
However,  \cite{JoMc} deals only with   probability preserving $\Bbb Z^d$-actions.  In this section we discuss the results from \cite{JoMc} and compare them with ours.
 
 \roster
 \item {\bf On the freeness of rank-one actions from \cite{JoMc}.}
According to  \cite{JoMc},  a probability preserving  $\Bbb Z^d$-action $T$ is called  of rank one     if 
Definition~1.1(i) holds.
It is not even assumed there that $T$ is free.
However, omitting the freeness, Johnson and McClendon  do not actually increase the class of actions.
Indeed, Definition~1.1(i) implies that $T$ is ergodic, but every ergodic action of an Abelian group $G$ is, in fact, a free action of the  quotient group $G/H$, where the subgroup $H\subset G$ is the stabilizer of the action, i.e.
$$
H=\{g\in G\mid T_g=I\}.
$$
If $G=\Bbb Z^d$ then $G/H$ is isomorphic to the direct product of $\Bbb Z^{d_1}$, for some $d_1\le d$,
with a finite Abelian group $A$.
Thus, $T$ is, upon the factorization $G\to G/H$, a free probability preserving action of $\Bbb Z^{d_1}\times A$.
\item{\bf ``F{\o}lner rank one'' is rank one according to Definition~1.1.}
By \cite{JoMc}, a $G$-action $T$ is called {\it of F{\o}lner rank one} if Definition~1.1(i) holds and $(F_n)_{n=1}^\infty$
is a  F{\o}lner sequence in $\Bbb Z^d$.
We first claim that then $T$ is free.
Indeed, denote by $H\subset G$ the stabilizer of $T$.
It follows from the definition of Rokhlin tower that $(F_n^{-1}F_n)\cap H=\{0\}$ for every $n>0$.
Hence 
$$
\Bigg(\bigcup_{n=1}^\infty F_n^{-1}F_n\Bigg)\cap H=\{0\}.
$$
However, $\bigcup_{n=1}^\infty F_n^{-1}F_n=\Bbb Z^d$ because $(F_n)_{n=1}^\infty$
is   F{\o}lner.
It follows that $H=\{0\}$, i.e. $T$ is free.
Since $(F_n)_{n=1}^\infty$
is   F{\o}lner, Definition~1.1(ii) is satisfied. 
Thus, $T$ is of rank one according to Definition~1.1.
The converse follows from
Corollary~1.11(ii) and Theorem~1.13.
Thus,  we obtain
that a probability preserving $\Bbb Z^d$-action  $T$ is of  F{\o}lner rank-one if and only if $T$ is of  rank
one according to Definition~1.1.
Therefore the main results of \cite{JoMc}: their Theorems~3.1, 4.7 and 5.1  (which are stated for the F{\o}lner rank-one $\Bbb Z^d$-actions) follow from  our Theorems~2.3 and 4.4.
\item 
In  \cite{JoMc, Theorem~3.2}, a sort of  conterexample to \cite{JoMc,  Theorem~3.1} is given.
It is shown there that the condition ``F{\o}lner rank one'' for $T$ can not be replaced with their ``rank-one'' for $T$ when describing the odometer factors of $T$.
We do not think that this  example is important, because  the whole trick is based entirely on the non-freeness of 
the
corresponding $\Bbb Z^2$-action $T$.
However, as we showed in (2), each  ``F{\o}lner rank one'' action is free.
\item
An ergodic  $G$-action is totally ergodic if and only if it has no non-trivial finite factors.
However, it is claimed in \cite{JoMc, Theorem~3.3} that there exist non-totally ergodic F{\o}lner rank-one $\Bbb Z^2$-actions without non-trivial finite factors.
This seeming contradiction is caused by the non-standard definition of total ergodicity in \cite{JoMc}.
As we understood from the proof of \cite{JoMc, Theorem~3.3}, by the total ergodicity of a $\Bbb Z^d$-action $T$ they mean
the {\it individual ergodicity} of $T$, i.e. that   the transformation $T_g$ is ergodic for each non-zero $g\in \Bbb Z^d$.
The proof of \cite{JoMc, Theorem~3.3} given there
is based on their analysis of finite factors for rank-one systems.
However, 
an easy alternative proof
follows from the joining theory (see \cite{Ru}, \cite{dJRu}) and has no direct relation to the rank one.
Indeed, let $S$ be a  transformation with MSJ.
Let $(X,\goth B,\mu)$ be the space of $S$.
Denote by $T$ the following $\Bbb Z^2$-action:   $T_{(n,m)}=S^n\times S^m$, $n,m\in\Bbb Z$,
on $(X\times X,\goth B\otimes\goth B,\mu\otimes\mu)$.
Of course, each factor of $T$ is a factor of  the transformation $T_{(1,1)}=S\times S$.
However, $S\times S$ has only 3 nontrivial proper factors: $\goth B\otimes\{\emptyset, X\}$, $\{\emptyset, X\}\otimes \goth B$ and $\goth B^{\odot 2}$  \cite{dJRu}.
The first two $\sigma$-algebras are also factors of $T$, while the latter one is not.
Thus, $T$ has only two proper factors which are weakly mixing. 
Hence $T$ is totally ergodic.
Since $T_{(1,0)}$ is not ergodic,  $T$ is not individually ergodic.
It remains to note that the 3-cut Chacon transformation, utilized  in the proof of \cite{JoMc, Theorem~3.3},
has MSJ \cite{dJRaSw}.
\item 
\cite{JoMc, Theorem~4.1} is a particular case of our Theorem~3.9: $G=\Bbb Z^d$ and the actions are probability preserving.
However, as we noted just above Theorem~3.9, if $G$ is Abelian then no proof is needed at all.
This follows immediately from the spectral theory.
Indeed, there is an alternative spectral definition for the probability preserving odometers: an ergodic
$G$-action $T$ is isomorphic to the odometer associated with $(\Gamma_n)_{n=1}^\infty$ if and only if $T$
has a pure discrete spectrum and the discrete spectrum $\Lambda(T)$ of $T$ equals
$\bigcup_{n=1}^\infty\widehat{G/\Gamma_n}$, which is a subgroup  of the dual group $\widehat G$ of characters on $G$.
\comment

 Halmos-von Neumann theory of systems with purely discrete spectrum.
Indeed, let $(X,\mu,T)$ is an ergodic $G$-action and let $(\Gamma_n)_{n=1}^\infty$ be a decreasing sequence of cofinite subgroups in $G$. 
Denote by $\Lambda(T)$ the discrete spectrum of $T$.
Then $\Lambda(T)$ is countable subgroup of the dual $\widehat G$.

\endcomment
If $T$ has a factor $G/\Gamma_n$ for each $n$ then $\widehat{G/\Gamma_n}\subset\Lambda(T)$
for each $n$.
Hence $L:=\bigcup_{n=1}^\infty\widehat{G/\Gamma_n}\subset\Lambda(T)$.
Then the factor of $T$  generated by the eigenfunctions with ``eigenvalues'' from $L$
 is isomorphic to  the odometer  associated with $(\Gamma_n)_{n=1}^\infty$, as desired.
Also, \cite{JoMc, Corollary~4.4} is a criterion when two $\Bbb Z^d$-odometers are conjugate.
It follows immediately from the Halmos-von Neumann theorem: two $\Bbb Z^d$-odometers $T$ and $T'$ associated with 
$(\Gamma_n)_{n=1}^\infty$ and $(\Gamma_n')_{n=1}^\infty$ respectively are conjugate if and only if $\Lambda(T)=\Lambda(T')$, i.e. if $\bigcup_{n=1}^\infty\widehat{G/\Gamma_n}=\bigcup_{n=1}^\infty\widehat{G/\Gamma_n'}$.
We also note that these spectral facts (or, rather, their non-Abelian counterparts) are no longer true for general non-Abelian odometers (see \cite{DaLe}, for example).
\item 
\S6 from \cite{JoMc} is devoted entirely to construction of a $\Bbb Z^2$-action $T$ that has no $\Bbb Z^2$-odometer factors but whose generators $T_{(0,1)}$ and $T_{(1,0)}$ have $\Bbb Z$-odometer factors.
We provide a simpler construction in Example~5.3.
\endroster

\comment

\head Appendix A: a nonsingular version of the Veech theorem
\endhead

\subhead A1. 
Simplex of relatively f.m.p. quasi-invariant measures
\endsubhead
Let $S=(S_g)_{g\in G}$ be an ergodic action of $G$ on a standard non-atomic probability space
$(Z,\goth Z,\xi)$.
Let $L$ be a locally compact Polish group.
A Borel map $\beta: G\times Z\to L$ is  called a {\it  cocycle} of $S$ if
$$
\beta(g_1g_2,z)=\beta(g_1, S_{g_2}z)\beta({g_2},z)\quad\text{at a.e. $z\in Z$ for every pair $g_1,g_2\in G.$}
$$
Let $\rho_\xi:G\times Z\to\Bbb R^*$ stand for the Radon-Nikodym cocycle of $S$, i.e.
$$
\rho_\xi(g,z):=\frac{d\xi\circ S_g}{d\xi}(z)\in\Bbb R_+^*.
$$
Denote by $\Delta_{\xi}$ the set of $S$-quasi-invariant probability measures $\omega$ on $Z$
such that $\rho_\omega(z)=\rho_\xi(z)$  at $\omega$-a.e. $z\in Z$.
Let $\Delta_\xi^0$ stand for the subset of ergodic measures in  $\Delta_\xi$.
Endow $\Delta_\xi$ with the natural   $\sigma$-algebra of Borel subsets.
Then $\Delta_\xi$ is a standard Borel space.
It  has the following simplicial properties \cite{GrSc}:
\roster
\item"---" $\Delta_\xi^0$ is a Borel subset of $\Delta_\xi$.
\item"---" For each $\gamma\in \Delta_\xi$, there is  a unique probability measure $m_\gamma$ on $\Delta_\xi^0$
such that  $\gamma=\int_{\Delta_\xi^0}\omega\,dm_\gamma(\omega)$.
This integral is called {\it the ergodic decomposition of $\gamma$.}
\item"---" In particular, if $\gamma\in \Delta_\xi^0$ then $m_\gamma=\delta_\gamma$, i.e. $m_\gamma$ is supported at $\gamma$.
\endroster

Let $\goth F\subset \goth B$ be  a factor of $S$.
Suppose that $S$ is  {\it $\goth F$-relatively f.m.p.}, i.e. 
the function $Z\ni z\mapsto\rho_\xi(g,z)\in\Bbb R^*_+$ is $\goth F$-measurable for each $g\in G$.
Denote by $\xi\otimes_{\goth F}\xi$ the $\goth F$-relatively independent product of $\xi$ with itself.
Let $\tau_j:X\times X\to X$ stand for the projection onto the $j$-th coordinate, $j=1,2$.
Of course,  
$$
(\xi\otimes_{\goth F}\xi)\circ\tau_j^{-1}=\xi\quad\text{for each $j=1,2$.}
$$
The dynamical system
 $(X\times X, \xi\otimes_{\goth F}\xi, (S_g\times S_g)_{g\in G})$ 
is nonsingular  and  $\rho_{\xi\otimes_{\goth F}\xi}=\rho_\xi$.
It follows that  $(X\times X, \xi\otimes_{\goth F}\xi, (S_g\times S_g)_{g\in G})$
is 
$\goth F$-relatively f.m.p. 
We note that $\xi\otimes_{\goth F}\xi$ is not necessarily ergodic.
Let $\xi\otimes_{\goth F}\xi=\int_{\Delta_{\xi}}\omega\,dm_{\xi\otimes_{\goth F}\xi}(\omega)$
be the ergodic decomposition of $\xi\otimes_{\goth F}\xi$.
Then
$$
\xi=(\xi\otimes_{\goth F}\xi)\circ \tau^{-1}_j=\int_{\Delta_{\xi}}\omega\circ\tau^{-1}_j\,dm_{\xi\otimes_{\goth F}\xi}(\omega).
$$
We deduce from the  simplicial  properties of $\Delta_{\xi}$ that  
$$
\omega\circ\tau^{-1}_j=\xi\quad\text{ at $m_{\xi\otimes_{\goth F}\xi}$-a.e. $\omega$ for $j=1,2$.}
\tag A-1
$$

\subhead A2. Compact skew product extensions and  the Veech theorem
\endsubhead
Let $K$ be a  compact Polish group.
Fix a cocycle  $\beta: G\times Z\to K$  of $S$.
The {\it $\beta$-skew product (compact) extension $S^\beta=(S_g^\beta)_{g\in G}$} of $S$ is the following $G$-action on the product space $(Z\times K,\xi\otimes\lambda_K)$:
$$
S^\beta_g(z,k):=(S_gz,\beta(g,z)k)\qquad\text{for all  $g\in G$}. 
$$
Of course, $S^\beta$ is $(\xi\otimes\lambda_K)$-nonsingular and 
$\rho_{\xi\otimes\lambda_K}(g,z,k)=\rho_{\xi}(g,z)$ at $(\xi\otimes\lambda_K)$-a.e. $(z,k)$ for each $g\in G$, 
i.e. $S^\beta$ is  a relatively f.m.p. extension of $S$.
The cocycle 
$$
G\times Z\ni(g,z)\mapsto(\beta(g,x),\beta(g,x))\in K\times K
$$
of $S$ is denoted by $\beta\otimes\beta$.
Given $h\in K$, we define a probability measure $\lambda_K^h$ on $ K\times K$ by setting
$$
\lambda_K^h(B\times C):=\lambda_K(B\cap Ch)\qquad\text{for all Borel subsets $B,C\subset K$.}
$$
It is straightforward to verify that if $S^\beta$ is ergodic then the dynamical system
$(Z\times K\times K,\xi_h, S^{\beta\otimes\beta})$ 
is ergodic and 
 $\rho_{\xi\otimes \xi_h}=\rho_\xi$ for each
$k\in K$.
Thus, we obtain a Borel embedding $K\ni h\mapsto\xi\otimes \lambda_K^h\in\Delta^0_\xi$ of $K$ into $\Delta^0_\xi$.

\proclaim{Claim A.1} Let $\zeta$ be a probability measure on $Z\times K\times K$.
If $\zeta$ is $S^{\beta\otimes\beta}$-quasi-invariant and ergodic, the push-forward of $\zeta$ to the space $Z\times K$ under the mapping  
$$
\tau_j:(z,k_1,k_2)\mapsto (z, k_j)
$$
 is $\xi\otimes \lambda_K$ 
 either for $j=1$ or for $j=2$
then there exists an element $h\in H$ such that $\zeta=\xi\otimes\lambda_K^h$.
\endproclaim
\demo{Proof}
It follows from the condition of the proposition that the push-forward of $\zeta$ to $Z$ is $\xi$.
Let 
$$
\zeta=\int_Z\delta_z\otimes\zeta_z\,d\xi(z)
$$ 
stand for the disintegration of $\zeta$ with respect to $\xi$.
We note that $\zeta_z$ is a probability on $K\times K$.
It is straightforward to verify that the mapping 
$$
\alpha:Z\times K\times K\ni(z,k_1,k_2)\mapsto k_1^{-1}k_2\in K
$$
is $S^{\beta\otimes\beta}$-invariant.
As $\zeta$ is ergodic, there exists $h\in K$ such that $\alpha=h$ almost everywhere in $\zeta$.
Therefore $\zeta$ is supported on the subset 
$$
\{(z, k, kh)\mid z\in Z,k\in K\}.
$$
This yields, in turn, that $\zeta_z$ is supported on the graph of the right rotation on $K$ by $h$, i.e.
there is a measurable field $Z\ni z\mapsto\eta_z$ of probabilities on $K$ such that
$\zeta_z(B\times C)=\eta_z(B\cap Ch)$ for all Borel subset $B,C\subset K$.
Therefore, the push-forward of $\zeta$ to $X\times K$ is $\int_Z\delta_z\otimes\eta_z\,d\xi(z)$.
On the other hand, by the condition of the proposition, this push-forward equals $\xi\otimes\lambda_K$.
Hence, by the uniqueness of the disintegration, $\eta_z=\lambda_K$ at a.e. $z\in Z$.
Thus,
$$
\zeta(A\times B\times C)=\int_Z\delta_z(A)\lambda_K(B\cap Ch)\,d\xi(z)=\xi(A)\lambda_K(B\cap Ch),
$$
as desired.
\qed
\enddemo

We now prove a nonsingular version of Veech theorem.
Denote by 
$$
R_h:Z\times K\to Z\times K
$$ 
the right translation by $h$ along the second coordinate.
Then the transformation $R_h$ preserves $\xi\otimes\lambda_K$ and commutes with $S^\beta$.
Let $\goth N_K$ stand for the trivial $\sigma$-algebra on $K$.

\proclaim{Theorem A.2} Let $\goth F$ be a factor of $S^\beta$ and let $\goth F\supset\goth Z\otimes\goth N_K$.
Then there is a compact subgroup $H\subset K$ such that $\goth F=\{A\subset Z\times K\mid R_hA=A\text{ for each }h\in H\}$.
\endproclaim

\demo{Proof}
Consider a relatively independent product of $S^\beta$ with itself over $\goth F$.
Formally, it is a  dynamical system $\big((Z\times K)^2,\zeta, (S^\beta_g\times S^\beta_g)_{g\in G}\big)$, where $\zeta=(\xi\otimes\lambda_K)\otimes_\goth F(\xi\otimes\lambda_K)$.
However, since $\goth F\supset\goth Z\otimes \goth N_K$, we can  assume (upon the natural identification) that $\zeta$ is a measure on $Z\times K\times K$ and
the corresponding dynamical system is $(Z\times K\times K, \zeta,S^{\beta\otimes\beta})$.
Let 
$$
\zeta=\int_{\Delta_\xi}\omega\, d m(\omega)\tag A-2
$$ 
stand for the ergodic decomposition of $\zeta$.
Of course,  $\zeta\circ \tau_1^{-1}=\xi\otimes\lambda_K$.
We remind that  $\xi\otimes\lambda_K$ is $S^\beta$-ergodic.
Hence  $\omega\circ\tau_1^{-1}=\xi\otimes\lambda_K$ at $m$-a.e. $\omega$ by \thetag{A-1}.
Then by~Claim~A.1, there is $h\in K$ such that $\omega=\xi\otimes\lambda_K^h$.
Therefore, we can rewrite~\thetag{A-2} as
$\zeta=\xi\otimes \int_{K}\lambda_K^h\, d m(h)$ for some probability measure $m$ on $K$.
Let 
$$
H:=\{k\in K\mid R_hA=A\text{ for each Borel }A\subset Z\times K\}.
$$
Then $H$ is a closed subgroup of $K$.
Arguing as in the proof of \cite{dJRu, Theorem~1.8.1}, we deduce that
$m(H)=1$ and that  $m$ is invariant under right rotations by the elements of $H$.
Hence $m=\lambda_H$.
We note that a Borel subset $A\subset Z\times K$ belongs to $\goth F$ if and only if
$\zeta(A\times A^c)=0$, where $A^c$ is the complement of $A$ in $Z\times K$.
The latter is equivalent to 
$$
(\xi\otimes\lambda_K)(A\triangle R_hA)=0\quad\text{for $\lambda_H$-a.e. $h\in H$},
$$
 which is, in turn, is equivalent to  
 $(\xi\otimes\lambda_K)(A\triangle R_hA)=0$ for each $h\in H$, as desired.
\qed
\enddemo

\endcomment

\comment
\proclaim{Lemma A2} Let $\zeta$ be a probability measure on $Z\times K\times K$.
If $\zeta$ is $S^{\beta\otimes\beta}$-quasiinvariant, the push-forward of $\zeta$ to the space $Z\times K$ under the mapping  
\thetag{A1}
 is $\xi\otimes \lambda_K$ 
  for $j=1$ (or for $j=2$)
  and 
  $$
  \rho_\zeta(g,z,k_1,k_2)=\rho_\xi(g,z)\quad\text{ at $\zeta$-a.e. $(z,k_1,k_2)\in Z\times K\times K$ }
  $$
for each $g\in G$ then a.e. ergodic component of 
\endproclaim

\endcomment

\comment
*****************************************************************

\head Appendix
\endhead

\definition{Definition} Given two finite measures $\lambda$ and $\rho$  on a standard Borel space $(Z,\goth Z)$, we write
$\lambda\prec\rho$ if $\lambda\ll\rho$ and $\frac{d\lambda}{d\rho}(z)\in\{0,1\}$ at $\rho$-a.e. $z\in Z$.
Equivalently, $\lambda\prec\rho$ iff there is a subset $A\subset Z$ such that $\lambda(B)=\rho(A\cap B)$ for each $B\in\goth Z$.

\enddefinition

Given a sequence $(\lambda_n)_{n=1}^\infty$ of finite measures on $(Z,\goth Z)$ such that $\lambda_1\prec\lambda_2\prec\cdots$, a $\sigma$-finite measure
$\lambda=\bigvee_{n=1}^\infty\lambda_n$ on $(Z,\goth Z)$ is well defined.
We note that there is an increasing sequence of subsets $A_1\subset A_2\subset\cdots$ such that $\lambda_n$ is supported on $A_n$ and $\lambda\restriction A_n=\lambda_n$ for each $n>0$.

Let $\kappa_n$ be a probability measure on $\Gamma_n$ and let $\nu_n$ be a finite measure on $G$ such that $\nu_n\succ\nu_{n-1}*\kappa_n$ for each $n>0$.
We assume that $\nu_0=\delta_{1_G}$. 
Suppose that given $g\in G$ and $n\in\Bbb N$, there is $N>0$ such that for each $m>N$,
$$
\delta_g*\nu_n*\kappa_{n+1}*\cdots*\kappa_m\ll\nu_m.
$$

For each $n\ge 0$, we let  $Z_n:=G\times\Gamma_{n+1}\times\Gamma_{n+2}\times\cdots$.
Then $Z_n$ is a totally disconnected Polish space.    
We endow $Z_n$ with a finite measure $\zeta_n:=\nu_n\otimes \bigotimes_{m>n}\kappa_m$.
Define a mapping
$\omega_n:Z_n\to Y$ by setting
$$
\omega_n(g,\gamma_{n+1},\gamma_{n+2},\dots)=
(g\Gamma_{n+1},g\gamma_{n+1}\Gamma_{n+2}, g\gamma_{n+1}\gamma_{n+2}\Gamma_{n+3},\dots).
$$
Let $\mu_n:=(\nu_{n}\otimes\kappa_{n+1}\otimes\kappa_{n+2}\otimes\cdots)\circ\omega_n^{-1}$.
Then $\mu_n$ is a finite Borel measure on $Y$ and $\mu_n\prec\mu_{n+1}$ for each $n\in\Bbb N$.
Let $\mu:=\bigvee_{n\in\Bbb N}\mu_n$.
Then $\mu$ is a $\sigma$-finite measure on $Y$.

**************************************************************

\endcomment

\comment

Let $\Cal R$ be a countable equivalence relation on a standard Borel space $X$.
Let $B$ be a Borel subset of $X$ and let $\mu$ be an ergodic $\Cal R\restriction B$-quasiinvariant finite measure on
$B$.

\proclaim{Claim 6.1} There exists a $\sigma$-finite $\Cal R$-ergodic Borel measure $\widetilde \mu$ on $X$ such that 
$\widetilde \mu\restriction B=\mu$.
\endproclaim

\demo{Proof} 
Let $\widetilde B:=\Cal R(B)$.
Select a Borel map $p:\widetilde B\to B$ such that $(x,p(x))\in\Cal R$
and $p(x)=x$ for each $x\in B$.
We now let 
$$
\widetilde\mu=\int_B\sum_{y\in p^{-1}(x)}\delta_y\,d\mu(x).
$$
\qed

\enddemo

\endcomment

\enddocument


\subhead Odometers
\endsubhead
Let $T$ be an ergodic transformation of a standard probability space $(X,\goth B,\mu)$.
Denote by $\Lambda(T)$ the set of all eigenvalues of $T$.
Then $\Lambda(T)$ is a countable subgroup of $\Bbb T$. 
If $T$ has pure point spectrum and $\Lambda(T)$ is a torsion group then $T$ is called {\it an odometer}.
By the Halmos-von Neumann theorem, two odometers $T$ and $S$ are isomorphic if and only if $\Lambda(T)=\Lambda(S)$.
Fix a dense countable torsion subgroup $K$ of $\Bbb T$.
Then there is a sequence $k_1,k_2,\dots$ of positive integers such that $k_1|k_2$, $k_2|k_3$, \dots and 
$K$ is  the union of an increasing  sequence of cyclic subgroups of period $k_n$, $n=1,2,\dots$.
The Pontryagin dual group $\widehat K$ 
to $K$  is the projective limit  
$\Bbb Z/k_1\Bbb Z\leftarrow \Bbb Z/k_2\Bbb Z\leftarrow\cdots$.
Endow $\widehat K$  with the Haar measure $\lambda$.
Fix an element $a\in\widehat K$   whose projection to $ \Bbb Z/ k_n\Bbb Z$ is $1+ k_n\Bbb Z$ for each $n=1,2,\dots$.
Denote by $T_K$ 
 the rotation of $(\widehat K,\lambda)$ by $a$.
 Then $T_K$ is an odometer and $\Lambda(T_K)=K$.

\subhead Rank-one transformations
\endsubhead
Let $T$ be a probability preserving rank-one transformation on $(X,\goth B,\mu)$.
Let $h_n$ be the hight of the $n$-th column and let $F_n:=\{0,1,\dots,h_n-1\}$.
For  $a\in F_n$, denote by $[a]_n$ the $a$-th level in the $n$-th column. 
Then the $n$-th column $X_n=\bigsqcup_{a\in F_n}[a]_n$.
Let $C_{n+1}$ be a subset of $F_{n+1}$ such that $[0]_n=\bigsqcup_{c\in C_{n+1}}[c]_{n+1}$.
Then for each $a\in F_n$, we have that 
$[a]_n=\bigsqcup_{c\in C_{n+1}}[a+c]_{n+1}$.
More generally, for each $m>n$ and $a\in F_n$, let $C_{n,m}:=C_n+C_{n+1}+\cdots+C_m$.
Then
$
[a]_n=\bigsqcup_{c\in C_{n+1,m}}[a+c]_{m}.
$
\proclaim{Theorem 1} The following are equivalent.
\roster\item"\rom{(i)}"
There is  $k>1$ and a factor map $\pi_k:X\to\Bbb Z/k\Bbb Z$, i.e. 
$$
\pi_k(Tx)=\pi_k(x)+1\pmod k\quad\text{ for a.e. $x\in X$.}
$$
(We identify $\Bbb Z/k\Bbb Z$ with $\{0,1,\dots,k-1\}$.)
\item"\rom{(ii)}"
There exist a sequence $(j_n)_{n=1}^\infty$ of integers  $j_n\in \{0,\dots,k-1\}$ such that
$$
\lim_{N\to\infty}
\sup_{m\ge n\ge N}
\frac{\#\{c\in
 C_{n,m}\mid c=j_{n}+\cdots+j_m\pmod k\}}
{\# C_{n,m}}=0.
$$
\item"\rom{(iii)}"
There is a sequence $q_1<g_2<\cdots$ such that 
$$
\sum_{n=1}^\infty
\frac{\#\{c\in C_{q_n+1,q_{n+1}}\mid c\ne 0\pmod k\}}
{\# C_{q_n+1,q_{n+1}}}<+\infty.
$$

\endroster
\endproclaim     
\demo{Proof} (i)$\Longrightarrow$(ii)
Consider a partition $X=\bigsqcup_{j=1}^k Y_j$, where $Y_j=\pi^{-1} (\{j\})$.
Then
$$
T^{b-a}([a]_n\cap Y_i)=[b]_n\cap Y_{i + b-a \pmod k}\quad\text{whenever $a,b\in F_n$}.\tag1
$$
For each $n>0$, $i\in\{0,\dots,k-1\}$ and $a\in F_n$, let
$d_{n,i}(a):=\mu([a]_n\cap Y_i)/\mu([a]_n)$.
Then the mapping 
$
F_n\ni a\mapsto\{d_{n,i}(a)\mid 0\le i < k\}
$
is constant, i.e. does not depend on $a$.
We consider it as {\it the  distribution of $k$ colors on $[a]_n$}.
Since the sequence $(\{[a]_n\mid a\in F_n\})_{n=0}^\infty$ of the $n$-th levels approximates  $\goth B$ as $n\to\infty$,
there is a map $i_n:F_n\to\Bbb Z/k\Bbb Z$ for each $n>0$ such that 
$$
\delta_n:=\max_{0\le i<k}d_{n,i}(a)=d_{n,i_n(a)}\to 1\qquad\text{as $n\to\infty$.}
$$
We consider  only $n$   large  so that $\delta_n>\frac12$.
Then we call $i_n(a)$ {\it the dominating color} on $[a]_n$.
We note the dominating color covers the same proportion $\delta_n$ of $[a]_n$ for every $a\in F_n$.
If follows from \thetag{1} that 
$$
\text{
$i_n(a)=i_n(b)$ \,  if and only if $b=a\pmod k$.}\tag2
$$
Given $\epsilon<\frac12$, there is
$N>0$ such that  $\delta_n>1-\epsilon^2$ for all $n\ge N$.
For each $m>n\ge N$, $[0]_n=\bigsqcup_{c\in C_{n+1,m}}[c]_m$.
By the  ``discrete Fubini theorem''
$$
\#\{c\in C_{n+1,m}\mid \mu([c]_m\cap Y_{i_n(0)})>
(1-\epsilon)\mu([c]_m)\}>(1-\epsilon)\# C_{n+1,m}.
$$
Thus, $i_n(0)$ is the dominating color for ``most'' of $c\in C_{n+1,m}$, i.e.
$$
\frac{\#\{c\in C_{n+1,m}\mid i_m(c)=i_n(0)\}}{\# C_{n+1,m}}>1-\epsilon.
$$
It follows from this and \thetag{2} that there is $j_{n+1,m}\in\{0,1,\dots,k-1\}$ such that
$$
\frac{\#\{c\in C_{n+1,m}\mid c=j_{n+1,m}\pmod k\}}
{\# C_{n+1,m}}
>1-\epsilon.
$$
We now let $j_{n+1}:=j_{n+1,n+1}$ for $n>N$.
Of course, $j_{n+1}$ does not depend on $\epsilon$.
If $N\le n\le m<r$ then $C_{n,r}=C_{n,m}+C_{m+1,r}$.
Hence $j_{n,r}=j_{n,m}+j_{m+1,r} \pmod k$.
This implies that $j_{n+1,m}=j_{n+1}+j_{n+2}+\cdots+j_m\pmod k$, as desired.

(ii)$\Longrightarrow$(iii)
Select an infinite increasing  sequence $(q_n)_{n=1}^\infty$ such that $q_1>N$ and $j_{N,q_1}=j_{N,q_2}=\cdots$.
It follows that 
 $j_{q_n+1}+\cdots+j_{q_{n+1}}=0\pmod k$ for each $n\in\Bbb N$.
 Passing to a subsequence in $(q_n)_{n=1}^\infty$ and using (ii), we may assume that
 $$
\frac{\#\{c\in C_{q_n+1,q_{n+1}}\mid c= 0\pmod k\}}
{\# C_{q_n+1,q_{n+1}}}>1-\frac1{2^n}
\tag3
$$
for all $n\in\Bbb N$.
This yields (iii).

(iii)$\Longrightarrow$(i)
We recall that for each $p>0$, there is a natural measure preserving isomorphism $\phi_p$ of the probability space $(X_p, \mu/\mu(X_p))$ onto the infinite product space $F_p\times C_{p+1}\times C_{p+2}\times\cdots$ equipped with the product of the equidistributions on $F_p$ and $C_j$, $j>p$, respectively such that if $\phi_p(x)=(f_p,c_{p+1},c_{p+2}\dots)$ for $x\in X_p$ and $f_p\ne h_p-1$ then $\phi_p(Tx)=(f_p+1, c_{p+1},\dots)$.\footnote{$\phi_p(x)=(f_p,c_{p+1},c_{p+2},\dots)$ if and only if $x\in\bigcap_{m>p}[f_p+c_{p+1}+\cdots+c_m]_m$ for a.e. $x$.}
We deduce from \thetag{3} and the Borel-Cantelli lemma that for a.e. $x\in X$, if $x=(f_{q_n},c_{q_n+1},\dots)\in X_{q_n}$ for some $n>0$ then  there is $R>n$ such that $\sum_{j=q_l+1}^{q_{l+1}}c_l=0\pmod k$ for each $l>R$.
We now can define $\pi_k(x)\in \Bbb Z/k\Bbb Z$ by setting
$$
\pi_k(x)=\Big(f_{q_n}+(c_{q_n+1}+\cdots +c_{q_{n+2}})+(c_{q_{n+2}+1}+\cdots +c_{q_{m+3}})+\cdots \Big)+ k\Bbb Z.\tag4
$$
It is well defined because the sum in the big parentheses  is finite for a.e. $x$.
Since $X=\bigcup_{n>0}X_{q_n}$, we obtain a mapping $\pi_k:X\supset X_{q_n}\ni x\mapsto \pi_k(x)\in \Bbb Z/k\Bbb Z$.  
It is well defined.
Moreover,  $\pi_k(Tx)=\pi_k(x)+1$ a.e., as desired.
Of course, $\pi_k$ is onto.
\qed
\enddemo

\comment

I prefer the following version of Theorem 1 (the above argument works for Theorem 1* as well as for Theorem 1):

\proclaim{Theorem 1*} The following are equivalent.
\roster\item"\rom{(i)}"
There is  $k>1$ and a factor map $\pi_k:X\to\Bbb Z/k\Bbb Z$
\item"\rom{(ii)}"
There is a sequence $q_1<g_2<\cdots$ such that 
$$
\sum_{n=1}^\infty
\frac{\#\{c\in C_{q_n+1,q_{n+1}}\mid c=0\pmod k\}}
{\# C_{q_n+1,q_{n+1}}}<+\infty.
$$
\endroster
\endproclaim

\endcomment

\remark{Remark 2} For each $a\in\{0,1,\dots,k-1\}$,
$$\lim_{n\to\infty}\mu\bigg(\Big(\pi_k^{-1}(a)\cap X_{q_n}\Big)\triangle\bigsqcup_{F_{q_n}\ni b=a\pmod k}[b]_{q_n}\bigg)= 0.\tag 5
$$
\endremark

Let $K:=\{k_n\in\Bbb N: k_n|k_{n+1}\text{ for each }n>0\}$.
Let $Z_K$ stand for the inverse limit $\Bbb Z/k_1\Bbb Z\leftarrow\Bbb Z/k_2\Bbb Z\leftarrow\cdots$.
Then $Z_K$ is a compact Abelian monothetic group.
Furnished with the Haar measure and an ergodic rotation, it is called an odometer.  
Of course, $Z_K$ is a factor of $X$ if and only if $\Bbb Z/k_n\Bbb Z$ is a factor of $X$ for each $n$.

Let  $K$ be a countable torsion subgroup of $\Bbb T$.
It is isomorphic to the union of an increasing sequence of cyclic subgroups of period $k_n$, $n=1,2,\dots$.

\proclaim{Corollary 3}
The following are equivalent.
\roster\item"\rom{(i)}"
There is   a factor map $\pi_K:(X,\mu)\to (\widehat K, \lambda)$, i.e. $\pi_K T=T_K\pi_K$.
\item"\rom{(ii)}"
There is a sequence $q_1<q_2<\cdots$ such that 
$$
\sum_{n=1}^\infty
\frac{\#\{c\in C_{q_n+1,q_{n+1}}\mid c\ne 0\pmod {k_n}\}}
{\# C_{q_n+1,q_{n+1}}}<+\infty.
\tag6
$$
\endroster
\endproclaim

\demo{Proof}
(ii)$\Longrightarrow$(i) It follows from \thetag{4} that the condition of Theorem~1(iii) is satisfied with $k_n$ in place of $k$ for each $n$.
Hence the factor maps $\pi_{k_n}:X\to \Bbb Z/k_n\Bbb Z$ are well defined via \thetag{4}.
Moreover, it is straightforward to verify  that  $\pi_{k_{n}}(x)=\pi_{k_{n+1}}(x)+\Bbb Z/k_n\Bbb Z$ for each $n$.
Hence a map $\pi_K:X\ni x\mapsto \pi_K(x):=(\pi_{k_1}(x),\pi_{k_2}(x),\dots)\in \widehat K$ is well defined.
It  intertwines $T$ with $T_K$.
\qed
\enddemo

The next problem is to find out when $\pi_K$ is an isomorphism, i.e. $T$ is isomorphic to 
$T_K$.
This happens if and only if for each $l>0$ and $\epsilon>0$ there is $n$ and a subset
$D\subset \Bbb Z/k_n\Bbb Z$ such that $\mu([0]_l\triangle \pi_{k_n}^{-1}(D))<\epsilon$.

\comment
Thus, it suffices to show that for for each $N>0$, there is $m>N$ such that $q_m>l$ and
$\mu(([0]_l\triangle \pi_k^{-1}(D))\cap X_{q_m})<\epsilon$.
Since
$$
\mu([0]_l\triangle (\pi_k^{-1}(D)\cap X_{q_m}))=
\frac{\#\{C_{l+1,q_m}\triangle\{f\in F_{q_m}\mid f+k\Bbb Z\in D\} \}}{\#C_{l+1,q_m}}+\overline o(1),
$$
we obtain the following theorem.
\endcomment

We deduce the following from Corollary 2* and (1).

\proclaim{Theorem 4} $(X,\mu,T)$ is isomorphic to  $Z_K$ if and only if
\roster
\item"\rom{(i)}"
There is a sequence $q_1<q_2<\cdots$ such that 
$$
\sum_{n=1}^\infty
\frac{\#\{c\in C_{q_n+1,q_{n+1}}\mid c\ne 0\pmod {k_n}\}}
{\# C_{q_n+1,q_{n+1}}}<+\infty.
$$
\item"\rom{(ii)}"
For each $l>0$ and $\epsilon>0$ there is $n,m\in\Bbb N$ and a subset
$D\subset \Bbb Z/k_n\Bbb Z$ such that  $q_m>l$ and 
$$
\frac{\#\{f\in C_{l+1,q_m}\mid f+k_n\Bbb Z\in D\} }{\#C_{l+1,q_m}}>1-\epsilon.
$$
\endroster
\endproclaim

\example{Example 5} Let $h_0:=0$ and $h_{n+1}:=4h_n+2^{n+1}$ for each $n\in\Bbb N$.
It follows that  $h_n=2^n(2^{n+1}-1)$ for each $n$. 
We let $F_n:=\{0,\dots,h_n-1\}$ and $C_{n+1}:=\{0,h_n,2h_n+2^{n+1},3h_n+2^{n+1}\}$ for each $n\ge0$.
Denote by $T$ the corresponding rank-one transformation.
We are going to apply Theorem 4. 
Let $q_n:=n$ and $k_n:=2^n$ for each $n\in\Bbb N$.
Then $C_{q_n+1,q_{n+1}}=C_{n+1}$ and every element of $C_{n+1}$ is divisible by $k_n$.
By Corollary~3, the 2-adic odometer is a factor of $T$.
\comment
For $x=(f_s,c_{s+1}, c_{s+2},\dots)$, we let
 $\pi_{2^n}(x)=f_s+2^n\Bbb Z$ for each $s\ge n$.
 Hence  
 $$
 \pi_{2^n}^{-1}(\{d\})=\bigcup_{s\ge 2^n}\bigsqcup_{ f+2^{n}\Bbb Z=d}[f]_s
 $$
 Since $[0]_r=[C_{r+1,{n-3}}+C_{{n-2}}]_{{n-2}}$  Hence 
  $$
 [0]_r\cap \pi_{2^n}^{-1}(\{d\})=\bigcup_{s\ge 2^n}\bigsqcup_{ f=d\pmod {2^n}}[f]_s
 $$

 **************************
 
 \noindent
 \endcomment
 We note that $C_{n+1}=C_{n+1}^{(1)}+C_{n+1}^{(2)}$, where $C_{n+1}^{(1)}=\{0,h_n\}$
 and $C_{n+1}^{(2)}=\{0,2h_n+2^{n+1}\}$.
 For each $m\ge n$, every element of $C_{m+1}$ is divisible by $2^{n+1}$.
Hence, for every subset
$D\subset \Bbb Z/2^{n+1}\Bbb Z$,
 $$
 \align
 \frac{\#\{f\in C_{l+1,m}\mid f+2^{n+1}\Bbb Z\in D\} }{\#C_{l+1,m}}&=
\frac{\#\{f\in C_{l+1,n}\mid f+2^{n+1}\Bbb Z\in D\} }{\#C_{l+1,n}}\\
&=
\frac{\#\{f\in C_{l+1,n-1}+C_{n+1}^{(2)}\mid f+2^{n+1}\Bbb Z\in D\} }{\#C_{l+1,n}}\\
&\le
\frac{\#( C_{l+1,n-1}+C_{n+1}^{(2)}) }{\#C_{l+1,n}}=\frac12.
\endalign
$$
 By Theorem~4, $T$ is not isomorphic to the 2-adic odometer.

 Let $\lambda\in\Lambda(T)$ and $\lambda^p=1$ for some $p$ which is not a power of $2$.
  Then the transformation $T^q$ is not ergodic.
  On the other hand, for each $n$, we have that either $h_n$ or $2h_n+2^{n+1}$ is not divisible by $q$.
  Hence, by \cite{Da, Theorem H(ii)}, $T^q$ is ergodic, a contradiction.
  Thus, 2-adic odometer is the maximal odometer factor of $T$.
\endexample


\Refs
\widestnumber\key{dJRaSw}




\ref\key{Aa}\by J. Aaronson
\book Infinite ergodic theory
\publ Amer. Math. Soc.
\yr 1997
\publaddr Providence, R. I.
\endref



\ref\key AdFrSi
\paper Rank one power weakly mixing nonsingular transformations
\by T. Adams,   N. Friedman and C. E. Silva
\jour  Ergodic Theory \& Dyn. Sys. 
\vol 21 
\yr 2001
\pages 1321--1332
\endref


\ref\key BeVa
\by  T. Berendschot, S. Vaes
\paper
Bernoulli actions of type $III_0$ with prescribed associated flow
\jour Inst. Math. Jussieu 
\toappear
\endref


\ref
\key CoWo
\by A. Connes and E. J. Woods
\paper  Approximately transitive flows and ITPFI factors
\jour  Ergod. Th. \& Dynam. Sys. 
\vol 5 \yr 1985
\pages 203--236
\endref


\comment


\ref\key CorPe
\by M. I. Cortez and S. Petite
\paper 
$G$-odometers and their almost one-to-one extensions
\jour J. Lond. Math. Soc. 
\vol 78\yr 2008
\pages 1--20
\endref 


\endcomment


\ref
\key Da1
 \by A. I. Danilenko
 \paper  Funny rank-one weak mixing for nonsingular abelian actions
 \jour  Israel J. Math. \vol 121 \yr 2001
 \pages  29--54
 \endref 
 
 
 \ref
 \key Da2
 \bysame  
 \paper Infinite rank one actions and nonsingular Chacon transformations
 \jour Illinois J. Math. \vol 48 \yr 2004
 \pages 769--786
 \endref



\ref
\key Da3
 \bysame
 \paper Actions of finite rank: weak rational ergodicity and partial rigidity
 \jour Ergodic Theory \& Dynam. Systems \vol 36 \yr 2016
 \pages 2138--2171
\endref

\ref
\key Da4
 \bysame
 \paper Rank-one actions, their (C,F)-models and constructions with bounded parameters
 \jour  J. d'Anal. Math. 
 \vol 139 \yr 2019
 \pages 697--749
 \endref
 
 
 \ref
 \key DadJ
  \paper Almost continuous orbit equivalence for non-singular homeomorphisms
\by A. I. Danilenko and A. del Junco
\jour Israel J. Math. 
\vol 183 
\yr 2011
\pages 165--188
\endref
 


\ref
\key DaLe
\by A. I.  Danilenko and  M.  Lema{\'n}czyk
\paper{ Odometer actions of the Heisenberg group}
\jour J. d'Anal. Math. 
\vol 128 
\yr 2016
\pages 107--157
\endref

\ref
\key DaSi
\by A. I.  Danilenko and C. E. Silva
\paper
Ergodic Theory: Nonsingular Transformations
\inbook Ergodic Theory (Encyclopedia of Complexity and Systems Science Series) 
\yr 2023
\pages 233--292
\publ Springer
\publaddr New York, NY
\endref


\ref
\key DaVi
\by A. I.  Danilenko and M. I. Vieprik
\paper
Explicit rank-1 constructions for irrational rotations
\jour Studia Math.
\vol 270
\yr 2023
\pages 121--144
\endref


\ref
\key DoHa
\by A. H. Dooley and T. Hamachi
\paper  Markov odometer actions not of product type
\jour  Ergodic Theory \& Dynam. Syst.
\vol  23 
\yr 2003
\pages 813--829
\endref



\ref \key dJ1
\by A. del Junco
\paper Transformations with discrete spectrum are stacking transformations
\jour Canadian J. Math.
\vol 28
\yr 1976
\pages 836--839
\endref


\ref \key dJ2
\bysame
\paper A simple map with no prime factors
\jour
Israel J. Math. 
\vol 104 \yr 1998
\pages 301--320
\endref

\ref\key dJRaSw
\paper Chacon's automorphism has minimal self-joinings
\by A. del Junco, M. Rahe, L. Swanson
\jour J. Analyse Math. 
\vol 37 
\yr 1980
\pages 276--284
\endref



\ref\key dJRu
\paper On ergodic actions whose self-joinings are graphs
\by A. del Junco and D. Rudolph
\jour Ergodic Theory Dynam. Systems 
\vol 7 \yr 1987\pages 531--557
\endref



\ref\key Ef
\by E. G. Effros
\paper Transformation groups and $C^*$-algebras
\jour Annals of Mathematics
\vol 81
\yr 1965
\pages 38--55
\endref


\ref\key HeRo
\by E. Hewitt  and K. A. Ross
\book
 Abstract harmonic analysis
 \vol I
 \yr 1963
\publ Springer
\publaddr Berlin
\endref 


\ref\key Fe
\by
S. Ferenczi
\paper Systems of finite rank
\jour Colloq. Math.
\vol 73
\yr 1997
\pages 35--65
\endref


\ref\key Fo--We
\paper Rank-one transformations, odometers, and finite factors
\by
    M. Foreman, S. Gao, A. Hill, C. E. Silva and  B. Weiss 
    \jour Isr J. Math.
 \vol 255\yr 2023
 \pages 231--249
    \endref
    
    
    
    \ref\key GrSc   
  \by G.  Greschonig and K. Schmidt
\paper Ergodic decomposition of quasi-invariant probability measures
\jour
Colloq. Math.
\vol 84/85
\yr 2000
\pages 495--514
\endref
    
    
    
    \ref\key Gri
    \by R. I. Grigorchuk
\paper On Burnside's problem on periodic groups
\lang Russian
\jour Funktsional. Anal. i Prilozhen.
\vol 14
\yr 1980
\pages 53--54
    \endref
    
    
    
    
    
    
    
 
    
    
 
    
    
    
    

    
    
    
    \ref
\key JoMc
\paper Finite odometer factors of rank one $\Bbb Z^d$-actions
\by A. S. A. Johnson and  D. M. McClendon 
\paperinfo preprint,  arXiv:2306.09477
\endref


\ref\key Ha
\paper A measure theoretical proof of the Connes-Woods theorem on AT-flows
\by T. Hamachi
\jour Pacific J. Math. 
\vol 154 \yr 1992
\pages 67--85
\endref


\ref\key HaSi\paper
On nonsingular Chacon transformations
\by T. Hamachi and C. E. Silva
\jour  Illinois J. Math.,\vol 44\yr 2000\pages 868--883
\endref


\ref\key Ki
\by J. King
\paper The commutant is the weak closure of the powers, for rank-1 transformations
\jour Ergodic Theory \& Dynam. Systems
\vol 6\yr 1986\pages  363--384
\endref

    
    
   \ref\key LiSaUg 
    \by S. Lightwood, A. {\c S}ahin and I. Ugarcovici
    \paper The structure and the spectrum of Heisenberg odometers
    \jour Proc. Amer. Math. Soc. 
    \vol 142 \yr 2014
    \pages 2429--2443
    \endref
    
    
  \ref\key Ru   
  \paper An example of a measure preserving map with minimal self-joinings, and applications
\by D. J. Rudolph
\jour J. Analyse Math. 
\vol 35 
\yr 1979
\pages 97--122
\endref

\ref\key Yu
\by H. Yuasa
\paper Uniform sets for infinite measure-preserving systems
\jour J. Anal. Math.
\vol 120
\yr 2013
\pages 333--356
\endref
    
    
\endRefs
\enddocument